\documentclass[preprint,sort&compress,12pt]{elsarticle}
\usepackage{bbm}
\usepackage{amsmath}
\usepackage{mathrsfs}
\usepackage{stmaryrd}
\usepackage{mathrsfs}
\usepackage{bm}
\usepackage{amsmath}
\usepackage{amsfonts}
\usepackage{amssymb}
\usepackage{amsthm}
\usepackage{lineno}
\usepackage[RGB]{xcolor}

\newcommand{\mm}{\mathrm}
\newcommand{\ml}{\mathcal}

\newcommand{\be}{\begin{equation}}
\newcommand{\bea}{\begin{equation}\begin{aligned}}
\newcommand{\beas}{\begin{equation*}\begin{aligned}}
\newcommand{\eeas}{\end{aligned}\end{equation*}}
\newcommand{\eea}{\end{aligned}\end{equation}}
\newcommand{\ee}{\end{equation}}

\renewcommand{\div}{{\rm div }}

\begin{document}
\begin{frontmatter}
\title{Asymptotic Behaviors of Global Solutions to the Two-Dimensional
Non-resistive MHD Equations with Large Initial Perturbations}

\author[FJ]{Fei Jiang}\ead{jiangfei0591@163.com}
\author[sJ]{Song Jiang}\ead{jiang@iapcm.ac.cn}
\address[FJ]{College of Mathematics and
Computer Science, Fuzhou University, Fuzhou, 350108, China.}
\address[sJ]{Institute of Applied Physics and Computational Mathematics, Huayuan Road 6, Beijing 100088, China.}
\begin{abstract}
This paper is concerned with the asymptotic behaviors of global strong solutions to the incompressible non-resistive viscous
magnetohydrodynamic (MHD) equations with large initial perturbations in two-dimensional periodic domains in Lagrangian coordinates.
First, motivated by the odevity conditions imposed in [Arch. Ration. Mech. Anal. 227 (2018), 637--662],
we prove the existence and uniqueness of strong solutions under some class of large initial perturbations,
where the strength of impressive magnetic fields depends increasingly on the $H^2$-norm of the initial perturbation values of
both velocity and magnetic field. {Then, we establish time-decay rates of strong solutions. Moreover, we find that
 $H^2$-norm of the velocity decays faster than the perturbed magnetic field. Finally, by developing some new analysis techniques,
 we show that the strong solution convergence in a rate of the field strength to the solution of the corresponding linearized problem
 as  the strength of the impressive magnetic field goes to infinity.
 In addition, an extension of similar results to the corresponding inviscid case with damping is presented. }
\end{abstract}
\end{frontmatter}


\newtheorem{thm}{Theorem}[section]
\newtheorem{lem}{Lemma}[section]
\newtheorem{pro}{Proposition}[section]
\newtheorem{cor}{Corollary}[section]
\newproof{pf}{Proof}
\newdefinition{rem}{Remark}[section]
\newtheorem{definition}{Definition}[section]

\section{Introduction}\label{introud}
\numberwithin{equation}{section}

We investigate the asymptotic behaviors of the global (-in-time) solutions to the following equations of incompressible
magnetohydrodynamic (MHD) fluids with zero resistivity:
\begin{equation}\label{1.1xx}
\begin{cases}
\rho  {v}_t+\rho {v}\cdot\nabla {v}+\nabla p- \mu\Delta v=\lambda M \cdot\nabla M/4\pi, \\
M_t+v\cdot\nabla M =M\cdot \nabla v,\\
\mathrm{div} v=\mathrm{div} M=0,
\end{cases}
\end{equation}
where the unknowns ${v}:={v}(x,t)$, $M:=M(x,t)$ and $p:=p(x,t)$ denote the velocity,  magnetic field and the sum of both magnetic and
kinetic pressures of MHD fluids respectively, and the three positive (physical) parameters $\rho$, $\mu$  and $ \lambda$
stand for the density, shear viscosity coefficient and permeability of vacuum, respectivly.

The global well-posedness of the system \eqref{1.1xx}, for which the initial data is a small perturbation around
a non-zero trivial stationary state (i.e., $v=0$, and $M$ is a non-zero constant vector $\bar{M}$, often called the impressive magnetic field),
has been widely investigated, see \cite{lin2015global,zhang2014elementary,ZTGS} and \cite{xu2015global} on the 2D and 3D
 Cauchy problems for \eqref{1.1xx} respectively, and see \cite{ren2016global} and \cite{TZWYJGw} on 2D and 3D initial-boundary value
 problems for \eqref{1.1xx} respectively. The existence of global solutions to the 2D Cauchy problem for \eqref{1.1xx}
 with large initial perturbations was obtained by Zhang under strong impressive magnetic fields \cite{ZTGS}. {As for the well-poseness
 of the 3D Cauchy and initial-boundary value problems for \eqref{1.1xx} with large initial perturbations, to our best knowledge,
 all available results are about the local (-in-time) existence, see \cite{fefferman2017local,HOCELETFCMDSFA,chemin2016local} for examples.
 We mention here that the corresponding compressible case has been also widely studied, see \cite{LYSYZG,LXLSNWDH,WJHWYF} and the references
 cited therein.} Since we are interested in the large-time behavior of a global solution and the asymptotic
 behavior of (a family of ) solutions with respect to the strength of the impressive magnetic field, we next briefly introduce the relevant
 progress on this topics, and our main results in this paper.

\subsection{Asymptotic behavior with respect to time}
It has been physically conjectured that in MHD fluids, the energy is dissipated at a rate that is independent of the resistivity \cite{CFCCRI}.
Hence, one can easily conclude that a non-resistive MHD fluid may still be dissipative. At present, this conclusion has been mathematically
verified for the global \emph{small perturbation solutions} of equations \eqref{1.1xx} around a non-zero trivial stationary state, in which the impressive magnetic field $M$ is given by
\begin{align}
\bar{M}=\bar{M}_N e_N.
\label{202008161713}
\end{align}
Here and in what follows $\bar{M}_N$ is a non-zero constant, $e_N$  a unit vector with the $N$-th component being $1$ and $N$ the spatial dimension.  In addition, we define $b:=M-\bar{M}_Ne_N$  and $\langle t\rangle^{-1} :=(1+t)^{-1}$.

The first mathematical verification result was probably given by Ren--Wu--Xiang--Zhang for the 2D Cauchy problem of \eqref{1.1xx}, they
established the following time-decay \cite{RXXWJHXZYZZF}:
\begin{align}
\langle t\rangle^{\frac{1+2k}{4}-\frac{\varepsilon}{2}}\|\partial_{x_2}^k(v,b)\|_{L^2(\mathbb{R}^2)}\leqslant c_I ,  \label{202008151719}
\end{align}
where $\varepsilon\in (0,1/2)$ is any given, $0\leqslant k\leqslant 2$, and the initial data $(v^0,b^0)$ of $(v,b)$
belongs to $H^{8}(\mathbb{R}^2)$. Here and in what follows, $c_I$ will denote a generic positive constant,
which may depend on the initial data $(v^0,b^0)$.
Then, Tan--Wang further obtained the almost exponential decay of solutions to the 2D/3D initial-boundary value problems \cite{TZWYJGw}:
\begin{align}
\langle  t\rangle^{n-2}\| (v,b)\|_{H^{2n+4}(\Omega)}\leqslant c_I,
\label{20200815}
\end{align}
where $n\geqslant 4$, the initial data $(v^0,b^0)$ belongs to $H^{4n}(\Omega)$, and $\Omega$ is a 2D/3D layer domain with finite height.
Later, Abidi--Zhang also got a decay rate of solutions to the 3D Cauchy problem \cite{ABIHZPOTG}:
\begin{align}
\label{20208142019}
\langle t\rangle^{1/4}\|(v,b)\|_{H^2(\mathbb{R}^3)}\leqslant c_I,
\end{align}
where $b^0=0$ and $v^0\in H^s$ with $s\in (3/2,3]$.
Very recently, under the assumption that the initial data are sufficiently smooth, Deng--Zhang further established
the following faster time-decay than \eqref{20208142019} \cite{deng2018large}:
\begin{align}
\label{2020606191026}
\langle t\rangle \|\nabla v\|_{L^2(\mathbb{R}^3)}+  {\langle t\rangle}^{1/2}\|(v,b)\|_{H^2(\mathbb{R}^3)}  \leqslant c_I.
\end{align}
In addition, Pan--Zhou--Zhu proved the existence of a unique global solution to the initial-boundary value problem of \eqref{1.1xx} defined
in a 3D periodic domain $\mathbb{T}^3$, and also obtained the following time-decay as a byproduct \cite{pan2018global}:
\begin{align}
\langle t\rangle^{(3-\sigma)/2}\| v \|_{H^{2(s-1)}(\mathbb{T}^3)}+\sum_{0\leqslant i\leqslant 1} \langle t\rangle^{(1-\sigma+2i)/2}\|\partial_3(v, b)\|_{H^{2(s-i)}(\mathbb{T}^3)}\leqslant c_I,
\label{202008031343}
\end{align}
where $0<\sigma<1$ and $(v^0,b^0)\in H^{2s+1}(\mathbb{T}^3) $ with $s \geqslant 5 $.

Motivated by the existence result on large perturbation solutions to the 2D Cauchy problem of \eqref{1.1xx} in \cite{ZTGS} and the time-decay
\eqref{2020606191026} for small perturbation solutions, we are interested in the time-decay of large perturbation solutions.
 {To the best of our knowledge, there is no result about time-decay rates of large perturbation solutions.
In this paper, we establish the time-decay rates of solutions with some class of large initial perturbations. }
Our results can be roughly described as follows.

 First, we prove the existence of a unique strong solution to the 2D initial-boundary value problem for \eqref{1.1xx} with periodic
 boundary conditions under some class of large initial perturbations, where the strength of impressive magnetic fields increasingly
 depend on the $H^2$-norm of  {the initially perturbed values of the velocity and magnetic field,
 and the initially perturbed data should satisfy some odevity conditions imposed by Pan--Zhou--Zhu in \cite{pan2018global}. Then, we show that
 the global solution enjoys the following decay in time: }
 \begin{align}
\label{202060619dsaf1026}u
\langle t\rangle^{3/2} \|  v\|_{H^1(\mathbb{T}^2)} +  {\langle t\rangle}  \| v \|_{H^2(\mathbb{T}^2)}+  {\langle t\rangle}^{1/2}\| b \|_{H^2(\mathbb{T}^2)}  \leqslant c_I,
\end{align}
where the initial datum $(b^0,v^0)$ belongs to $H^2$. We refer the reader to Theorem  \ref{202008161120} for the details (or see Theorems \ref{201904301948} and \ref{202006211824saf} for the version in Lagrangian coordinates). {We should point out here that by virtue of
\eqref{202060619dsaf1026}, the velocity in $H^2$-norm decays faster than the perturbation magnetic field, while in \eqref{2020606191026},
the $H^2$-norm of the velocity enjoys the same decay rate as the one of the perturbation magnetic field.}

\subsection{Asymptotic behavior with respect to the strength of impressive magnetic fields}
It is well-known that a non-resistive MHD fluid,  {the motion of which is described by the system \eqref{1.1xx},}
exhibits elastic characteristics. In particular, a MHD fluid strains when stretched and will quickly return to its original rest state
by the magnetic tension once the stress is removed \cite{JFJSOMITIN}. This means that the magnetic tension will have stabilizing effects
in the motion of MHD fluids. Moreover, the larger the strength of an impressive magnetic field is, the stronger this stabilizing effect will be,
see the inhibition phenomenon of flow instabilities by magnetic fields \cite{JFJSOUI,JFJSSETEFP,WYTIVNMI}. It is worth to mention here that Bardos--Sulem--Sulem used hyperbolicity of \eqref{1.1xx} with $\mu=0$ to establish an interesting global existence result of classical
solutions with small initial data in the H\"older space $H^s(\mathbb{R}^3)$ \cite{BCSCSPLL}; also see \cite{CYLZGW,HLBXLYP} for the case
of Sobolev spaces \cite{CYLZGW,HLBXLYP}). {Remark that such a result is not known for the system \eqref{1.1xx}
in 3D when the magnetic field is absent (the 3D incompressible Euler equations). }

In \cite{ZTGS} Zhang also found an interesting mathematical result that {the solution of the nonlinear system \eqref{1.1xx}
converges to the solution of some linear equations obtained from \eqref{1.1xx} under the use of the stream function
as $\bar{M}_2\to\infty$. More precisely, let $M=(\partial_2 ,-\partial_1 )^{\mm{T}}(\psi-\bar{M}_2x_1)$, then the system \eqref{1.1xx}
reduces to the following system: }
\begin{equation}
\label{1sfa.1xasfdsax}
\begin{cases}
\rho  {v}_t+\rho {v}\cdot\nabla {v}+\nabla p- \mu\Delta v-\bar{M}_2\lambda\partial_2(\partial_2\psi,-\partial_1\psi)^{\mm{T}} /4\pi\\
=\lambda(\partial_2\psi,-\partial_1\psi)\cdot
\nabla(\partial_2\psi,-\partial_1\psi)^{\mm{T}} /4\pi, \\
\psi_t+v\cdot\nabla \psi +\bar{M}_2v_1=0,\\
\mathrm{div} v=0,\\
(v,\psi)|_{t=0}=(v^0,\psi^0),
\end{cases}
\end{equation}
where the subscript $\mm{T}$ denotes the transposition, $\partial_i:=\partial_{x_i}$ and $x_i$ is the $i$-th component of $x\in \mathbb{R}^2$.
Thus, Zhang proved that the solution $(v,\psi)$ of \eqref{1sfa.1xasfdsax} converges to the solution $(v^{\mm{L}},\psi^{\mm{L}})$ of the
following linear pressureless equations as $\bar{M}_2\to \infty$:
\begin{equation}\label{1sfa.1xasafdsfdsax}
\begin{cases}
\rho  {v}_t^{\mm{L}} - \mu\Delta v^{\mm{L}}-\bar{M}_2\lambda\partial_2(\partial_2\psi^{\mm{L}},-\partial_1\psi^{\mm{L}})^{\mm{T}} /4\pi
=0, \\
\psi_t^{\mm{L}}+ \bar{M}_2v_2^{\mm{L}}=0,\\
(v^{\mm{L}},\psi^{\mm{L}})|_{t=0}=(v^0,\psi^0).
\end{cases}
\end{equation}
However, no convergence rate in the strength $\bar{M}_2$ is given in \cite{ZTGS}.
 {In this paper we shall prove a similar result in Lagrangian coordinates,
 and further provide a convergence rate in $\bar{M}_2$ of the global solution to \eqref{1.1xx} }  by developing some new analysis techniques.
 More precisely, we shall show that the difference between the solution of \eqref{1.1xx} in Lagrangian coordinates and
 {the solution of the corresponding linearized system}
 can be bounded from above by $c_I \bar{M}_2^{-1/2}$, see Theorem \ref{201912041028} for details.

Roughly speaking,  {the proof of our two results mentioned above is based on a key observation that the deviation function $\eta$
of MHD fluid particles enjoys the estimate }
\begin{align} \nonumber
 \|\bar{M}_2 \partial_2 \eta\|_0^2\leqslant c_I,
 \end{align}
where $c_I$ depends on the initial total mechanical energy. {The above estimate can be extended to the case of  higher-order derivatives of $\partial_2\eta$.}
Moreover, if $\eta$ additionally satisfies the odevity conditions, then we formally have an important inverse relation:
\begin{align}
\label{202102042101}
\nabla \eta \propto  c_I/\bar{M}_2,
\end{align}
see Section \ref{2020102012102} for a detailed discussion.
This relation intuitively not only reveals that the (nonlinear) solutions of \eqref{1.1xx} in Lagrangian coordinates can be approximated
by the (linear) solutions of the corresponding linearized equations for $\bar{M}_2 \gg c_I$, {but also  provides
a convergence rate in $\bar{M}_2$.}
Since the nonlinear solutions can be approximated by the linear solutions, we naturally expect the existence of strong solutions
  under some class of large initial perturbations as in \cite{ZTGS}. In fact, in \cite{ZTGS} Zhang first obtained the (linear) solution of \eqref{1sfa.1xasafdsfdsax}, then proved the existence of the small error solution $(v-v^{\mm{L}},\psi-\psi^{\mm{L}})$ as $\bar{M}_2 \gg c_I$,
  and finally got the large solution $(v,\psi)$ by adding the linear solution and the small error solution together.
  It is worth to mention that the relation \eqref{202102042101} allows us to directly establish the existence of solutions
  under some class of large initial perturbations by one-step procedure, rather than Zhang's three-step procedure.

 {We mention that recently,} some authors studied the case of inviscid, non-resistive MHD fluids with zero resistivity,
i.e., the viscosity term in the system \eqref{1.1xx} is replaced by the velocity damping term $\kappa \rho v$ with $\kappa$ being
the damping coefficient. Wu--Wu--Xu first proved the existence of a unique global solution with algebraic time-decay
to the 2D Cauchy problem, provided that the initial perturbation $(v^0,b^0)$ is small in $H^n(\mathbb{R}^2)$ with $n$ sufficiently large \cite{WJHWYFXXJG}.
Recently, Du--Yang--Zhou also obtained the existence of a unique global solution with exponential time-decay to the initial-boundary
value problem in a strip domain $\Omega$, provided that the initial perturbation $(v^0,b^0)$ around some non-trivial equilibrium
is small in $H^7(\Omega)$ \cite{DYYWZYOSJMA}.
Motivated by \cite{WJHWYFXXJG,DYYWZYOSJMA}, we can extend our aforementioned results in this paper on the asymptotic behavior of solutions
in the viscous case to the inviscid case with the damping term $\kappa\rho v$, and show that for the inviscid
case with damping, the decay in time is exponentially fast, just as in \cite{DYYWZYOSJMA},
while the convergence rate in $\bar{M}_2$ as $\bar{M}_2\to\infty$ of a strong solution of the original nonlinear system to the solution
of the corresponding linear system is in the form of $c_I \bar{M}_2^{-1}$, which is faster than that for the viscous case,
see Theorem \ref{2019043019sdfa48} for the details.

Finally, we mention that the asymptotic behaviors of solutions with respect to other parameters, such as the Mach and Alfv\'en numbers,
in MHD fluids have been also extensively investigated, see, for example, \cite{CBHQCSS} and the references cited therein.

The rest of this paper is organized as follows: In Section \ref{202008161246} we introduce our main results including
the existence of a unique strong solution with some class of large initial data to the 2D equations \eqref{1.1xx} in a periodic domain
in Lagrangian coordinates,  {and the time-decay of the strong solution, and the convergence rate as $\bar{M}_2\to\infty$ of the
strong solution as well as the extension to the inviscid case with damping term, i.e., Theorems \ref{201904301948}, \ref{202006211824saf}, \ref{201912041028} and \ref{2019043019sdfa48}, the proofs of which are given in Sections \ref{20190sdafs4301948}--\ref{201912asdfas062021},
        respectively. Finally, in Section \ref{201912062021}, we provide the proof of the local well-posedness for the equations
        of viscous, non-resistive MHD fluids and inviscid, non-resistive MHD fluids with damping, respectively.}

\section{Main results}\label{202008161246}
In this section we describe the main results in details. To begin with, we reformulate the equations \eqref{1.1xx} in Lagrangian coordinates.  {Recalling that \eqref{1.1xx} is considered with in a 2D periodic domain, we see, without loss of generality, that it suffices to consider the periodic domain $\mathbb{T}^2$ with $\mathbb{T}:= \mathbb{R}/\mathbb{Z}$.}

Let $(v,M)$ be the solution of the 2D system of equations \eqref{1.1xx}, and the flow map $\zeta$ be the solution to
\begin{equation}
\label{201806122101}
            \begin{cases}
\partial_t \zeta(y,t)=v(\zeta(y,t),t)&\mbox{ in }\mathbb{T}^2 \times\mathbb{R}^+ ,
\\
\zeta(y,0)=\zeta^0(y)&\mbox{ in }\mathbb{T}^2,
                  \end{cases}
\end{equation}
where $\zeta^0(y)$ satisfies $\det \nabla \zeta^0=1$ and ``$\det$" denotes the determinant.

Since $v$ is divergence-free, then
\begin{align}
\det \nabla \zeta =1
\label{202008161709}
\end{align}
as well as $\det \nabla \zeta^0=1$. Thus, we define $\ml{A}^{\mm{T}}:=(\nabla \zeta)^{-1}:=(\partial_j \zeta_i)^{-1}_{2\times 2}$.
In particular, by virtue of \eqref{202008161709},
\begin{align}  {\mathcal{A}}=\left(\begin{array}{cc}
            \partial_2\zeta_2 &  - \partial_1\zeta_2\\
                - \partial_2\zeta_1&   \partial_1\zeta_1
                 \end{array}\right). \nonumber
\end{align}

 We temporarily introduce some differential operators involving $\mathcal{A}$, which will be used later. The differential operators $\nabla_{\ml{A}}$, $\mm{div}_\ml{A}$ and $\Delta_{\ml{A}}$ are defined by $\nabla_{\ml{A}}f:=(\ml{A}_{1k}\partial_kf,
\ml{A}_{2k}\partial_kf)^{\mm{T}}$, $\mm{div}_{\ml{A}}(X_1,X_2)^{\mm{T}}$ $:=\ml{A}_{lk}\partial_k X_l$,
 $\mm{curl}_\mathcal{A}f:=\partial_1 \mathcal{A}_{1k}\partial_kf_2- \mathcal{A}_{2k}\partial_kf_1$    and $\Delta_{\mathcal{A}}f:=\mm{div}_{\ml{A}}\nabla_{\ml{A}}f$
for a scalar function $f$ and a vector function $X:=(X_1,X_2)^{\mm{T}}$, where $\ml{A}_{ij}$ denotes the $(i,j)$-th entry of the matrix $\mathcal{A}$. It should be remarked that  we have used the Einstein convention of summation over repeated indices, and $\partial_k=\partial_{y_k}$.
In addition, thanks to \eqref{202008161709}, we have
\begin{equation}
\partial_k  \mathcal{A}_{ik} =0.
\label{201909261909}
\end{equation}

Let $\mathbb{R}^+=(0,\infty)$, $\nu= \mu/{\rho}  $ and
\begin{equation*}
( u ,B,q)(y,t)=( v,M,p/\rho)(\zeta(y,t),t) \;\;\mbox{ for }\; (y,t)\in \mathbb{T}^2 \times\mathbb{R}^+ .
\end{equation*}
By virtue of the equations \eqref{1.1xx} and \eqref{201806122101}$_1$, the evolution equations for $(\zeta,u,q)$ in Lagrangian
coordinates read as follows.
\begin{equation}\label{01dsaf16asdfasfdsaasf}
                              \begin{cases}
\zeta_t=u , \\[1mm]
 u_t+\nabla_{\ml{A}}q- \nu \Delta_{\ml{A}}u=
 \lambda B\cdot\nabla_{\mathcal{A}}B/4\pi\rho, \\[1mm]
B_t-B\cdot\nabla_{\mathcal{A}} u=0,\\
\div_\ml{A}u=0,\\[1mm]\div_\ml{A}B=0.
\end{cases}
\end{equation}

 We can  derive from \eqref{01dsaf16asdfasfdsaasf}$_3$ the differential version of magnetic flux conservation \cite{JFJSOMITIN}:
\begin{equation} \nonumber
\ml{A}_{jl}B_j=\ml{A}_{jl}^0B_j^0,
 \end{equation}
 which yields
\begin{eqnarray}
 \label{0124}  B=\nabla\zeta \ml{A}^{\mm{T}}_0 B^0.
\end{eqnarray}
Here and in what follows, the notation $f_0$ as well as  $f^0$ denote the value of the function $f$ at $t=0$.
If we assume
\begin{equation}
\label{201903081437}
\ml{A}^{\mm{T}}_0 B^0=\bar{M} \;\;  \mbox{ (i.e., }B^0= \partial_{\bar{M}}\zeta^0\mbox{)},
 \end{equation}
 where $\bar{M}$ is defined by \eqref{202008161713} with $N=2$, then
 \eqref{0124} reduces to
 \begin{eqnarray}
\label{0124xx}   B=\partial_{\bar{M}}\zeta  .
\end{eqnarray}
Here we should point out that $B$ given by \eqref{0124xx} automatically satisfies \eqref{01dsaf16asdfasfdsaasf}$_3$
and \eqref{01dsaf16asdfasfdsaasf}$_5$. Moreover, from \eqref{0124xx} we see that the magnetic tension in Lagrangian coordinates
has the relation
\begin{equation}\nonumber
 B\cdot \nabla_{\ml{A}} B= \partial_{\bar{M}}^2\zeta.
 \end{equation}

Let $I$ denote a $2\times 2$ identity matrix, $m^2= \lambda  \bar{M}_2^2 /4\pi\rho   $, $\eta=\zeta-y$
and
$$\tilde{\mathcal{A}}=\left(\begin{array}{cc}
            \partial_2\eta_2 &  - \partial_1\eta_2\\
                - \partial_2\eta_1&   \partial_1\eta_1
                 \end{array}\right). \nonumber $$
Consequently, under the assumption \eqref{201903081437}, the equations \eqref{01dsaf16asdfasfdsaasf} are equivalent to
the following system:
\begin{equation}\label{01dsaf16asdfasf}
                              \begin{cases}
\eta_t=u , \\[1mm]
 u_t+\nabla_{\ml{A}}q- \nu \Delta_{\ml{A}}u=
  m^2 \partial_2^2\eta, \\[1mm]
\div_\ml{A}u=0
\end{cases}
\end{equation}
where $B=m(\partial_2\eta+e_2)$ and $\mathcal{A}=\tilde{\mathcal{A}}+I$.

For the well-posedness of \eqref{01dsaf16asdfasf} defined in $\mathbb{T}^2$, we impose the initial condition:
\begin{equation}\label{01dsaf16asdfasfsaf} (\eta, u)|_{t=0}=(\eta^0,  u^0)\;\;\mbox{ in }\;\mathbb{T}^2.
\end{equation}

 Before stating our main results, we introduce some notations which will be frequently used throughout this paper.
\begin{enumerate}
  \item[(1)] Basic notations:  $\overline{I_\infty}:=\mathbb{R}^+_0:=[0,\infty)$, $I_T:=(0,T)$ for $0<T\leqslant \infty$, $\overline{I_T}:=[0,T]$ for $T\in \mathbb{R}^+$,  $\Omega_T:=\mathbb{T}^2\times I_T$,
$\int:=\int_{(-1,1)^2}$, $(w)_{\mathbb{T}^2}:=\int w\mm{d}y$, $\alpha=(\alpha_1,\alpha_2)$ denotes the multi-index with respect to the variable $y$.
  \item[(2)] Simplified Banach spaces:
\begin{align}
& L^r:=L^r (\mathbb{T}^2)=W^{0,r}(\mathbb{T}^2),\;\; {H}^i:=W^{i,2}(\mathbb{T}^2),\;\; H:=H^1,\nonumber \\
&  H^i_\sigma:=\{u\in H^i~|~\mm{div}u=0\},\
  H^{i+1}_{1}:=\{\eta\in H^{i+1} ~|~\det\nabla(\eta +y )=1\},\nonumber \\
  & H^{i+1}_2:=\{u\in H^{i}~|~\partial_2u\in  H^{i}\},\   {H}^{i+1}_{1,2}:=H^{i}_{1}\cap H^{i+1}_2,\nonumber  \\
  &{ \underline{X}:=\{w\in X\cap L^2~|~(w)_{\mathbb{T}^2}=0\}, } \nonumber
  \end{align}
 {where $X$ denotes a Banach space, $1< r\leqslant \infty$ and $i \geqslant 0$ are integers. }
\item[(3)]
Simplified spaces of functions with values in a Banach space (please refer to Sections 1.3.10--1.3.12 in \cite{NASII04}
for the properties of functions with values in a Banach space):
 {\begin{align}
&\mathfrak{C}^0(\overline{I_T} ,H^i):=\{u\in L^2(\Omega_T) ~|~u\in C( {I_T}\setminus \mathfrak{Z},H^i)
\mbox{ for some zero-measurable set }\mathfrak{Z}\subset I_T\},  \nonumber\\
&\mathcal{U}_{i+1,T}:=\{u\in C^0(\overline{I_T}, H^{2i+2})~|~u_t\in C^0(\overline{I_T}, H^{2i}),\
(u,u_t)\in L^2(I_T,H^{2i+3}\times H^{2i+1})\},\label{202102111149} \\
&\underline{\mathcal{U}}_T:=\{u\in \mathcal{U}_{1,T}~|~(u)_{\mathbb{T}^2}=0\}.  \nonumber
\end{align}  }
 \item[(4)] Simplified function classes: for integer  $i \geqslant 1$,
\begin{align}
& H^i_{*}:=\{\xi\in H^i~|~ \xi(y)+y   : \mathbb{R}^2 \to \mathbb{R}^2 \mbox{ is a } C^1\mbox{-diffeomorphism mapping}\},\nonumber\\
&
  \underline{H}^{i,*}_{1,T}:=\{\eta\in C^0(\overline{I_T}, \underline{H}^{i}_{1})~|~ \eta(t) \in H^{i}_{*}\mbox{ for each }t\in \overline{I_T} \},\nonumber\\
& {C}^0_{B,\mm{weak}}({I_T} ,L^2):= \mathfrak{C}^0({I_T}, L^2)\cap L^\infty(I_T, L^2)\cap C^0(\overline{I_T}, L^2_{\mm{weak}}), \nonumber\\
& {\mathfrak{C}^0( {I_T},{H}^{i+1}_2 ):=\{\eta\in C(\overline{I_T} ,{H}^{i}) ~|~\partial^\alpha\partial_2\eta\in   {C}^0_{B,\mm{weak}}({I_T} ,L^2)\mbox{ for any }|\alpha|=i \},  } \nonumber \\
& {\underline{\mathfrak{H}}^{i+1,*}_{1,2,T}:=\{\eta\in \mathfrak{C}^0({I_T}, H^{i+1}_2 )~|~ \eta(t) \in \underline{H}^{i}_{1}\cap H^{i}_{*}\mbox{ for each }t\in \overline{I_T}\}, } \nonumber\\
&  { {\mathfrak{U}}_T^{i+1}:=  \{u\in C^0(\overline{I_T}, H^{i}) ~|~\partial^\alpha u\in {C}^0_{B,\mm{weak}}( {I_T} ,L^2)
\mbox{ for any }|\alpha|=i+1 \}, }  \nonumber\\
 &  { \underline{\mathfrak{U}}_T^{i+1}:=\{u\in {\mathfrak{U}}_T^{i+1}~|~(u)_{\mathbb{T}^2}=0\}. }  \nonumber
\end{align}
\item[(5)] Simplified norms: for integers $i\geqslant 0$ and $n\geqslant i$,
\begin{align}
&\|\cdot \|_{i} :=\|\cdot \|_{H^i(\mathbb{T}^2)},\;\; \|\nabla^i \cdot\|^2_0:=\sum_{|\alpha|=i}\| \partial^\alpha \cdot\|^2_0,
\;\;  \|\cdot \|_{i+1,2} :=\sqrt{\|\cdot \|_{i}^2+\|\nabla^i\partial_2\cdot \|_{0}^2},\nonumber \\
                 & {\mathfrak{E}_{n,i}(t):= \|\partial_2^i(\nabla \eta, u,m\partial_2\eta)(t)\|_{n-i}^2, \quad
                                    \mathfrak{E}_{n,i}^0:= \mathfrak{E}_{n,i}(0). } \nonumber
               \end{align}
\item[(6)]  {General} constants: ${c}_i$ ($1\leqslant i \leqslant 3$) and ${c}_i^{\kappa}$ ($1\leqslant i\leqslant 4$)
are fixed constants which may depend on the parameters $\nu$ and $\kappa$ respectively, but not on $m$. If not stated explicitly,
$c_0$, $c$, $c^\kappa$, $C$ and $C^\kappa$ will denote generic positive constants, which may vary from one place to another. Moreover,\end{enumerate}
\begin{itemize}
  \item $c_0$ is independent of any parameter;
  \item $c$ and $c^\kappa$ may depend on $\nu$ and $\kappa$ respectively (\emph{but not on $m$});
  \item  $C$ depends on $\nu$ and $\sqrt{\mathfrak{E}_{2,0}}$, and increases with respect to $\sqrt{\mathfrak{E}_{2,0}}$. In particular, $C$ only depends on $\nu$  and the norm $\|u^0\|_2$ for the case $\eta^0=0$;
\item  $C^\kappa$  depends   on $\kappa$ and $\|( u^0,\eta^0,m\partial_2\eta^0)\|_{4}$, and increases with respect to  $\|( u^0,\eta^0,m\partial_2\eta^0)\|_{4}$.
\end{itemize}
 {In addition, $A\lesssim_0 B$, $A\lesssim B$ and $A\lesssim_\kappa B$ mean that $A\leqslant c_0B$, $A\leqslant cB$ and  $A\leqslant c^\kappa B$, respectively.}

\subsection{Existence of global solutions}\label{2020102012102}

Before stating the global existence result of solutions to the initial-value problem \eqref{01dsaf16asdfasf}--\eqref{01dsaf16asdfasfsaf}
in some classes of large data under strong magnetic fields, let us fist mention the heuristic idea, which leads us to study this topic.

First, multiplying \eqref{01dsaf16asdfasf}$_2$ with $u$ in $L^2$, we obtain the basic energy identity:
\begin{align}
\label{ssebdaiseqinM0846dfgssgsdxx}
\frac{1}{2}\frac{\mm{d}}{\mm{d}t}\left( \|u\|_0^2+ \|m \partial_2 \eta\|_0^2\right)
+ {\nu} \|\nabla_{\mathcal{A}} u\|_0^2= 0,
\end{align}
which implies
\begin{align}
\label{202008221356}
 {\|u\|_0^2 }+ \|m\partial_2 \eta\|_0^2
+ {2\nu} \int_0^t \| \nabla_{\mathcal{A}} u\|_0^2\mm{d}\tau= {\|u^0\|_0^2 }+ \|m\partial_2 \eta^0\|_0^2=:I_0.
\end{align}
We call $I_0$ the initial total mechanical energy, which includes the kinetic energy, and the perturbation magnetic energy
that could be regarded as the potential energy.
We easily see from \eqref{202008221356} that $\|\partial_2 \eta\|_0\to 0$  as $m\to\infty$  for fixed $I_0$.
 {This basic relation motivates us to expect that}
 the deformation quantity $\nabla \eta$ may be small, when $m$ is sufficiently large. Fortunately,
 {this is indeed the case for $(\eta,u)$ satisfying the additional odevity conditions
 imposed by Pan--Zhou--Zhu in \cite{pan2018global}, see \eqref{202007301500} for details. }

We rewrite \eqref{01dsaf16asdfasf}$_2$--\eqref{01dsaf16asdfasf}$_3$ as a nonhomogeneous system of the Stokes equations:
\begin{equation}\label{s0106pnnnn}
 \begin{cases}
u_t+\nabla q-\nu\Delta  u =\mathfrak{F} , \\[1mm]
\div u =\mm{div}_{\tilde{\mathcal{A}}}u,
\end{cases}
\end{equation}
where we have defined $\mathfrak{F}:= m^2\partial_2^2\eta +\mathfrak{N}$, $\mathfrak{N}:=\mathcal{N}^\nu-\nabla_{\tilde{\ml{A}}}q$, $\mathcal{N}^\nu:= \partial_l(\mathcal{N}^\nu_{1,l}, \mathcal{N}^\nu_{2,l})^{\mm{T}}$ and
$\mathcal{N}^\nu_{j,l} := \nu (\mathcal{A}_{kl}\tilde{\mathcal{A}}_{km} + \tilde{\mathcal{A}}_{ml} )\partial_mu_j$.
 {This formally reveals that }
the system \eqref{s0106pnnnn} can be approximated by the corresponding linear system, if $\nabla \eta$ is sufficiently small.
Since the linear system admits a global solution, the nonlinear system \eqref{s0106pnnnn} may also admit a global solution
in some classes of large data under the strong magnetic fields. This result read as follows.
\begin{thm}
\label{201904301948}
There are positive constants $c_1\geqslant 4$, $c_2>0$ and a sufficiently small constant $c_3\in (0,1]$, such that for any
$(\eta^0,u^0)\in (\underline{H}^{3}_{1}\cap H^3_*)\times \underline{H}^{2}$ and $m$ satisfying the incompressible condition
$\mm{div}_{\mathcal{A}^0}u^0=0$, the odevity conditions
\begin{align}
(\eta^0_1,  u^0_1)(y_1,y_2)= (\eta^0_1,  u^0_1)(y_1,-y_2)\;\mbox{ and }\;
(\eta^0_2,u^0_2)(y_1,y_2)=- (\eta^0_2, u^0_2)(y_1,-y_2),
\label{202007301500}
\end{align} and the condition of strong magnetic fields
\begin{equation}
\label{201909281832} m\geqslant \frac{1}{c_3}\max\left\{\left( c_1 \mathfrak{E}_{2,0}^0 e^{c_2 \mathfrak{E}_{2,1}^0}\right)^{1/4},
c_1 \mathfrak{E}_{2,0}^0 e^{c_2 \mathfrak{E}_{2,1}^0}\right\},
\end{equation}
then the initial value problem \eqref{01dsaf16asdfasf}--\eqref{01dsaf16asdfasfsaf} admits a unique global strong solution
$(\eta, u,q)\in \underline{H}^{3,*}_{1,\infty}\times \underline{\mathcal{U}}_\infty\times C^0(\mathbb{R}^+_0,\underline{H}^2)$.
Moreover, the solution $(\eta,u)$ enjoys the stability estimate:
\begin{align}
& \mathfrak{E}_{2,0}(t) +  \int_0^t(\| u\|_3^2+  \|m\partial_2 \eta \|_2^2)\mm{d}\tau\lesssim\mathfrak{E}_{2,0}^0
e^{c_2 \mathfrak{E}_{2,1}^0} (1+\mathfrak{E}_{2,1}^0) \quad\mbox{ for any }t\geqslant 0,\label{20190safd5041053}
\end{align}
where $\| \eta \|_3+\|u \|_2+\|m\eta \|_{3,2}\lesssim_0 \sqrt{\mathfrak{E}_{2,0}}$. In addition,
\begin{align}
& \|\nabla q\|_1 \lesssim_0 \|\partial_2 u\|_1^2
  +m^2 (   \|\partial_2^3 \eta\|_0 \|\partial_2 \eta\|_2 +\|\partial_2^2 \eta\|_1^2 ), \label{20211191238} \\
&\|u\|_{i+1,2}\lesssim_0
                     \begin{cases}
                   \|\partial_2 u\|_0+\|\partial_2^2\eta\|_1\|\partial_2u\|_1  \ & \hbox{ for }i=0; \\
                    \|\partial_2 u\|_i & \hbox{ for }i=1\mbox{ and }2,
                     \end{cases} \label{2020111211028} \\
& \|\eta\|_{i+1,2} \lesssim_0 \|\partial_2 \eta\|_i,\ \|\eta\|_{L^\infty} \lesssim_0 \sqrt{\| \partial_2\eta\|_0\|\partial_2\eta\|_1}.
\label{202011119211101}
\end{align}
\end{thm}
\begin{rem}We can easily construct a family of $(\eta^0,u^0)$ satisfies all the assumptions in Theorem \ref{201904301948},
where $\eta^0\neq 0$ and $u^0\neq 0$. In fact, let
$\bar{\eta}=\bar{u}=(\sin x_1 \cos x_2,-\cos x_1\sin x_2)$. Because $\mm{div}\bar{\eta}=\mm{div}\bar{u}=0$, for sufficiently
small $\varepsilon$, there exists a function pair $(\eta^0 ,u^0 )$ enjoying the form
$(\eta^0 ,u^0 )= (\varepsilon\bar{\eta}+\varepsilon^2 \eta^{\mm{r}},\bar{u}+\varepsilon u^{\mm{r}})$,
where $(\eta^{\mm{r}},u^{\mm{r}})$ satisfies
$\| \eta^{\mm{r}}\|_3+\| u^{\mm{r}}\|_2\leqslant c_0 $,  \begin{align}
   \begin{cases}
   -\Delta \eta^{\mm{r}} +\nabla \beta_1=0, \\
   \mm{div} \eta^{\mm{r}}=  \mm{div}\big( (\bar{\eta}_1 +\varepsilon \eta^{\mm{r}}_1)( -\partial_2(\bar{\eta}_2 +\varepsilon\eta^{\mm{r}}_2),  \partial_1(\bar{\eta}_2+\varepsilon \eta^{\mm{r}}_2) )^{\mm{T}}\big), \\
     (\eta^{\mm{r}})_{\mathbb{T}^2}  = 0,
   \end{cases}
\end{align}
and
\begin{align} \begin{cases}
-\Delta u^{\mm{r}} +\nabla \beta_2=0, \\
\mm{div} u^{\mm{r}}=\varepsilon^{-1}\mm{div}_{\tilde{\mathcal{A}}^0}(u^0+\varepsilon u^{\mm{r}}),\\
({u}^{\mm{r}})_{\mathbb{T}^2}  = 0,
\end{cases}
\nonumber
\end{align}
 {for a proof of which we refer to \cite[Proposition 5.1]{JFJSWZhangwei} }.
It is easy to check that for sufficiently small $\varepsilon$, $(\eta^0 ,u^0 )$ is non-zero,
belongs to $ (\underline{H}^{3}_{1}\cap H^3_*)\times \underline{H}^{2} $, and satisfies $\mm{div}_{\mathcal{A}^0}u^0=0$ and \eqref{202007301500}.
We further take $m=\varepsilon^{-1}$ to immediately see that $(\eta^0 ,u^0 )$ and $m$ satisfy the condition of strong magnetic fields \eqref{201909281832} for sufficiently small $\varepsilon$. Furthermore, $\|(\nabla \eta^0, u^0,m\partial_2\eta^0)\|_{2}\leqslant c_0$
for some constant $c_0$ independent of $\varepsilon$ and $m$.
\end{rem}
\begin{rem}
Noting that the initial perturbation magnetic field ``$B^0-\bar{M}$'' is equal to $\bar{M}_2\partial_2\eta$, we see from \eqref{201909281832}
that the strength of the impressive magnetic field increasingly depends on the $H^2$-norm of the initial velocity and perturbation
magnetic field.
\end{rem}
\begin{rem}
In the above theorem, we have assumed the condition $(\eta^0)_{\mathbb{T}^2}=(u^0)_{\mathbb{T}^2}=0$.
If $((\eta^0)_{\mathbb{T}^2}, (u^0)_{\mathbb{T}^2})\neq 0$, we can define $\bar{\eta}^0:=\eta^0-(\eta^0)_{\mathbb{T}^2}$ and
$\bar{u}^0:=u^0-(u^0)_{\mathbb{T}^2}$. Then, by Theorem \ref{201904301948}, there exists a unique global strong solution $(\bar{\eta},\bar{u},\bar{q})$ to the initial value problem \eqref{01dsaf16asdfasf}--\eqref{01dsaf16asdfasfsaf} with initial data $(\bar{\eta}^0,\bar{u}^0)$. It is easy to verify that
$(\eta,u,q):=(\bar{\eta}+ t(u)_{\mathbb{T}^2}+ (\eta^0)_{\mathbb{T}^2},\bar{u}+(u)_{\mathbb{T}^2},\bar{q})$
is just the unique strong solution of \eqref{01dsaf16asdfasf}--\eqref{01dsaf16asdfasfsaf} with initial data $({\eta}^0, {u}^0)$.
\end{rem}
\begin{rem}\label{202008191942}
Since $(\eta^0, u^0)$ satisfies the odevity conditions \eqref{202007301500}, the strong solution $(\eta,u)$ of  \eqref{01dsaf16asdfasf}--\eqref{01dsaf16asdfasfsaf} with an associated pressure function $q$
  enjoys the same odevity conditions as $(\eta^0,u^0)$ does, i.e.,
\begin{align}
(\eta_1,  u_1)(y_1,y_2,t)= (\eta_1,  u_1)(y_1,-y_2,t)\;\mbox{ and }\;
(\eta_2,u_2)(y_1,y_2,t)=-(\eta_2,  u_2)(y_1,-y_2,t).
\label{20208192055}
\end{align}
\end{rem}

Now, we briefly describe the proof idea of Theorem \ref{201904301948}.
Motivated by \eqref{202008221356}, we naturally except that for given value $\mathfrak{E}_{2,0}^0$,
the solution $(\eta,u)$ enjoys the estimate
 {\begin{align}
\label{20200asdf8211944}
\|(u,m\partial_2 \eta)\|_2\leqslant {C}\quad\mbox{ for sufficiently large }\; m.
\end{align} }
Thus, we want to derive the \emph{a priori} estimate of $(\eta,u)$ like
\begin{align}
\label{20200asdf8sadfa211944}
\|(u,m\partial_2 \eta)\|_2\leqslant   {C}/2
\end{align}
under the \emph{a priori} assumption \eqref{20200asdf8211944}.

 {However, if we follow the above idea, we find that the assumption \eqref{20200asdf8211944} does not suffice
to establish the \emph{a priori} estimate \eqref{20200asdf8sadfa211944}, expect further requiring the additional assumption: }
\begin{equation}
\label{aprpiasdfsaose1}
 \|\nabla \eta\|_2^2   \leqslant  C.
\end{equation}
More precisely, we can conclude that there are constants $K$ (depending on $ \mathfrak{E}_{2,0}^0$, $\mathfrak{E}_{2,1}^0$
and $m$) and $\delta$, such that
\begin{equation}
\label{aprpiose1}
\sup_{0\leqslant t\leqslant T}  \mathfrak{E}_{2,0}(t)  \leqslant{K}^2/ 4,
\end{equation}
if
\begin{align}
&\sup_{0\leqslant t\leqslant T}( \|(\nabla \eta,m\partial_2 \eta)(t)\|_2^2+\|\partial_2 u(t)\|_1^2 )
\leqslant K^2\;\;\mbox{ for any given }\; T>0\label{aprpiosesnewxxxx}
\end{align}
and
\begin{equation}
\label{aprpiose1snewxxxxz}
\max\{K^{1/2}, K^2\}/m\in (0,\delta]\footnote{From \eqref{aprpiosesnewxxxx} we easily get $K/m\leqslant\delta$.
We should point out that the term $K^{1/2}$ in \eqref{aprpiosesnewxxxx} can be replaced by $K$
to establish Theorem \ref{201904301948} with ``$1/2$" in place of ``$1/4$" in \eqref{201909281832}.  Here we further choose $K^{1/2}$
in order to make sure that $m^{-1}$ does not appear on the right hand of \eqref{202005201033}. }.\end{equation}

The above \emph{a priori} stability estimate, together with a local well-posedness result on \eqref{01dsaf16asdfasf}--\eqref{01dsaf16asdfasfsaf},
immediately yields Theorem \ref{201904301948}.  The detailed proof will be presented in Section \ref{20190sdafs4301948}.  In addition, the proof of the existence of a unique local solution will be provided in Section \ref{202009111946}.
It should be remarked that the odevity conditions \eqref{20208192055} play an important role in the derivation of the \emph{a priori} stability
estimate, see the key Lemma \ref{201805141072dsafa}.

\subsection{Asymptotic behaviors of solutions}

Now, we turn to stating the asymptotic behaviors of the global solution given in Theorem \ref{201904301948}. To begin with,
we state the result of the asymptotic behavior with respect to the time.
\begin{thm}
\label{202006211824saf}
Let $(\eta,u,q)$ be the global solution of \eqref{01dsaf16asdfasf}--\eqref{01dsaf16asdfasfsaf} established in Theorem \ref{201904301948},
then
\begin{align}
&  \sum_{i=1}^2\Big(\langle t\rangle^{i}\|\partial_2^i\eta (t)\|_{3-i}^2
+\langle t\rangle^{i+1}\|\partial_2^{i}(u,m\partial_2 \eta )(t) \|_{2-i}^2\nonumber\\
& + \int_0^t\big( \langle \tau\rangle^{i+1}\| \partial_2^{i} u \|_{3-i}^2
 + \langle \tau\rangle^i \| m  \partial_2^{i+1}\eta  \|_{2-i}^2 \big) \mm{d}\tau\Big) \nonumber \\
 & +\langle t\rangle \| (u,m\partial_2 \eta )(t) \|_{2}^2  + \int_0^t  \langle \tau\rangle  \|   u  \|_{3}^2   \mm{d}\tau\leqslant C   \label{202005201033}
\end{align}
and
\begin{align}
&\sum_{i=1}^2 \left( \langle t\rangle^{i+1} \|u(t)\|_{3-i}^2 +   \langle t\rangle^{i+1}\int_0^t
e^{\nu (\tau-t)/2}\| u \|_{4-i}^2     \mm{d}\tau\right)  \leqslant  (1+m^2)C\quad\mbox{for any }t\geqslant 0. \label{20200safa5201033x}
\end{align}
\end{thm}
\begin{rem}
In view of \eqref{2020111211028}, \eqref{202011119211101},
\eqref{202005201033} and Poinc\'are's inequality, we get from \eqref{202005201033} that
\begin{align}
&\langle t\rangle^{5/2}\|m\eta (t)\|_{L^\infty(\mathbb{T}^2)}^2 +\sum_{i=1}^2\left(\langle t\rangle^{i}\|\eta (t)\|_{4-i,2}^2
 +\int_0^t\langle \tau\rangle^i \| m  \partial_2\eta  \|_{3-i,2}^2   \mm{d}\tau\right.\nonumber\\
&\left. + \langle t\rangle^{i+1}\|( u,m\partial_2 \eta )(t) \|_{3-i,2}^2+ \int_0^t  \langle \tau\rangle^{i+1}
\|u \|_{4-i,2}^2\mm{d}\tau\right)\leqslant C, \qquad\forall\, t\geqslant 0.
\end{align}
\end{rem}

From \eqref{202005201033} we immediately see that $u$ enjoys the same time-decay rate as $m\partial_2 \eta$ does.
Similar results were obtained for the global solution with small perturbation, see \eqref{202008151719}--\eqref{202008031343}.
Next, we briefly explain the basic idea how to get the faster time-decay of the velocity \eqref{20200safa5201033x}.

It is well-known that any solution of the following homogeneous Stokes equations
\begin{equation}\label{s0106safasdfaapnnnn}
 \begin{cases}
u_t+\nabla q-\nu\Delta  u =0 , \\[1mm]
\div u =0
\end{cases}
\end{equation}
 {decays exponentially in time. Thus, for the nonhomogeneous case in \eqref{s0106pnnnn}, the algebraic time-decay of $\|u\|_{3-i}$
 depends on $\mathfrak{F}$ where $i=1$ and $2$. Since the decay rate of the linear term
 $\|\partial_2^2\eta\|_{2-i}$ is $\langle t\rangle^{-(i+1)}$,
    we naturally except that $\|u(t)\|_{3-i}$ also decays in the rate of $\langle t\rangle^{-(i+1)}$.}  Fortunately, by employing carefully
     estimates, one sees that the nonlinear term $\mathfrak{N}$ does not prevent us from obtaining the desired decay rate, see Section \ref{202008231326}.

Recalling that the solution $\eta$ in Theorems \ref{201904301948} satisfies
\begin{align}
&\zeta:= \eta(y, t)+y   : \mathbb{R}^2\to\mathbb{R}^2 \mbox{ is a }C^1\mbox{ diffeomorphism mapping}, \label{2312018031adsadfa21601xx}
\end{align}
one can easily recover the decay result in Theorem \ref{202006211824saf} from Lagrangian coordinates to the one in Eulerian coordinates:
\begin{thm}
\label{202008161120}
Let $(v,\eta,q)$ be the global solution given in Theorem \ref{201904301948}, $\zeta^0=\eta^0+y$, $\zeta^{-1}_0$ denotes
the inverse function and
$$(v,M,p)=(u,m\partial_2 \zeta,  q)|_{y=\zeta^{-1}}.$$
Then $(v,M,p)$ belongs to $C^0(\mathbb{R}_0^+,\underline{H}^2_\sigma\times H^2 \times  \underline{H}^2 )$,
and is a unique global strong solution of the initial value problem \eqref{1.1xx} with initial data
$(v,M)|_{t=0}=(u^0,m\partial_2 \zeta^0)|_{y=\zeta^{-1}_0}$ in $\mathbb{T}^2$. Moreover, the solution enjoys the decay estimates
\begin{align}
& \sum_{i=0}^2\left(\langle t\rangle^{i+1}\| (v,b)(t)\|_{2-i}^2+ \int_0^t \big( \langle \tau\rangle^{i+1}\| v\|_{3-i}^2
+ \langle \tau\rangle^i \| b\|_{2-i}^2\big)  \mm{d}\tau \right)\leqslant C   \nonumber
\end{align}
and
\begin{align}
\sum_{i=1}^2 \left( \langle t\rangle^{i+1} \| v(t)\|_{3-i}^2 +   \langle t\rangle^{i+1}\int_0^t
e^{\nu (\tau-t)/2}\| v \|_{4-i}^2     \mm{d}\tau\right)   \leqslant (1+m^2)C, \nonumber
\end{align}
where $b:=M-\bar{M}$.
\end{thm}
\begin{rem}
Since the periodic cell of $\mathbb{T}^2$ is bounded, we can also establish a result of almost exponential decay, where
the decay rate is faster than Tan--Wang's result \eqref{20200815} in \cite{TZWYJGw} under the same regularity on initial data.
In fact, if the initial data $(\eta^0,u^0)$ in Theorem \ref{201904301948} is in $H^{n+1}\times H^n$ with $n\geqslant 3$, and $m$ satisfies \eqref{201909281832} with $\mathfrak{E}_{n,0}^0$ in place of $\mathfrak{E}_{2,0}^0$, then there exists a unique classical solution
$(v,M,p)\in C^0(\mathbb{R}_0^+,\underline{H}^n_\sigma\times H^n\times\underline{H}^n)$ to the initial value problem \eqref{1.1xx}
with initial data $(v,M)|_{t=0}=(u^0,m\partial_2 \zeta^0)|_{y=\zeta_0^{-1}}$. Moreover,
\begin{align}
&  \sum_{i=0}^n\left(\langle t\rangle^{i+1}  \|  ( v,b) \|_{n-i}^2+ \int_0^t  \langle \tau\rangle^{i+1}\| v\|_{n+1-i}^2\mm{d}\tau\right)+  \sum_{i=1}^n \int_0^t \langle \tau\rangle^i \| b\|_{n-i}^2   \mm{d}\tau \leqslant c_I,  \nonumber \\
&
  \sum_{i=1}^n  \left(\langle t\rangle^{i+1}  \| v\|_{n+1-i}^2  +   \langle t\rangle^{i+1}\int_0^t
e^{\nu (\tau-t)/2}\| v \|_{n+2-i}^2  \mm{d}\tau  \right) \leqslant (1+m^2)c_I,\nonumber
\end{align}
where $b:=M-\bar{M}$ and the constant $c_I$ depends on $\nu$ and $\mathfrak{E}_{n,0}^0$.
\end{rem}

Next, we further state the asymptotic behavior of solutions with respect to $m$. Noting that the inhomogeneous term $\mathfrak{N}$
in $\mathfrak{F}$ includes $\partial_i\eta_j$ for $1\leqslant i$, $j\leqslant 2$, and
$\|\nabla \eta\|_{L^\infty}\lesssim_0 \|\partial_2\eta\|_2\leqslant K/2m$ by \eqref{aprpiose1} and \eqref{2020043013014}, we formally see
that $\mathfrak{N}\to 0$ as $m\to \infty$ for fixed $\mathfrak{E}_{2,0}^0$. Thus, the solution $(\eta,u)$ established in
Theorem \ref{201904301948} converges in the rate $m^{-1/2}$ to the solution $(\eta^{\mm{L}},u^{\mm{L}})$ of the corresponding linearized
system as $m\to \infty$. More precisely, we have
 \begin{thm}\label{201912041028}
Let $(\eta,u,q)$ be the global solution of \eqref{01dsaf16asdfasf}--\eqref{01dsaf16asdfasfsaf} given in Theorem \ref{201904301948}.
\begin{enumerate}[(1)]
   \item Then, one can use the initial data of $(\eta,u)$ to construct a function pair
   $(\eta^{\mm{r}}, u^{\mm{r}})\in  \underline{H}^3\times \underline{H}^2$, such that the following linear pressureless
   initial-value problem
   \begin{equation}\label{202001070914}
 \begin{cases}
\eta_t^{\mm{L}}=u^{\mm{L}} , \\[1mm]
u_t^{\mm{L}} -\nu\Delta  u^{\mm{L}}= m^2\partial_2^2\eta^{\mm{L}}  , \\[1mm]
\div u^{\mm{L}} =0,\\
(\eta^{\mm{L}},u^{\mm{L}})|_{t=0}=(\eta^0+\eta^{\mm{r}}, u^0+u^{\mm{r}})
\end{cases}
\end{equation}
admits a unique strong solution
$(\eta^{\mm{L}}, u^{\mm{L}} )\in C^0(\mathbb{R}_0^+, \underline{H}^{3}_\sigma)\times\underline{\mathcal{U}}_\infty $. Moreover,
\begin{itemize}
  \item  $(\eta^{\mm{L}}, u^{\mm{L}} )$ also satisfies the odevity conditions \eqref{20208192055} as $(\eta ,u)$ does;
  \item  the function pair $(\eta^{\mm{r}}, u^{\mm{r}})$ satisfies
\begin{align}
&\mm{div}(u^0+u^{\mm{r}})=\mm{div}(\eta^0+\eta^{\mm{r}})=0,\label{202011032122}  \\
  &\| u^{\mm{r}} \|_2\lesssim_0 \|\partial_2 \eta^0\|_2\|  u^0\|_2,\;\mbox{ and }\;
\|\partial_2^j\eta^{\mm{r}}\|_{3-j}\lesssim_0\|\partial_2\eta^0\|_2\|\partial_2^j\eta^0\|_{3-j},\;\;\; j=0,1.\label{202011032123}
\end{align}
\end{itemize}
 \item  Let $(\eta^{\mm{d}}, u^{\mm{d}} )= (\eta-\eta^{\mm{L}} , u-u^{\mm{L}} ) $, then for any $t\in \mathbb{R}^+_0$,
 \begin{align}
&   \|\eta^{\mm{d}}(t)\|_3^2+ \int_0^t \| ( u^{\mm{d}},  m \partial_2 \eta^{\mm{d}})\|_2^2\mm{d}\tau\leqslant C m^{-1},
\label{201905041053xxx} \\
& \sum_{i=1}^2\bigg(\langle t\rangle^{i}\|\partial_2^i\eta^{\mm{d}}(t)\|_{3-i}^2 +\langle t\rangle^{i+1}
  \|\partial_2^{i}( u^{\mm{d}},m\partial_2 \eta^{\mm{d}} )(t) \|_{2-i}^2 \nonumber\\
&  \qquad \  +\int_0^t (\langle \tau\rangle^{i+1}\| \partial_2^{i} u^{\mm{d}} \|_{3-i}^2+\langle  \tau\rangle^i
 \| m\partial_2^{i+1}\eta ^{\mm{d}} \|_{2-i}^2 )\mm{d}\tau\bigg) \leqslant C m^{-1}  \label{2020sfa05201033},
\end{align}
where the error function
$(\eta^{\mm{d}}, u^{\mm{d}})$ enjoys the estimates \eqref{2020111211028} and \eqref{202011119211101} with $(\eta^{\mm{d}},u^{\mm{d}})$
in place of $(\eta,u)$.
 \item
If  $\eta^0$ further satisfies the additional regularity
\begin{align}\partial_1(\eta^0,u^0)\in H^3\times H^2,
\label{2020111222152}
\end{align} then for any $t\geqslant 0$,
\begin{align}
&\langle t\rangle\|(u^{\mm{d}},m \partial_2 \eta^{\mm{d}})(t)\|_2^2+ \int_0^t \langle \tau \rangle\|   u^{\mm{d}} \|_3^2\mm{d}\tau
\leqslant C m^{-1}(1+\|\partial_1 \eta^0\|_3^2+\|\partial_1 u^0\|_2^2). \label{202011032108}
\end{align}
\end{enumerate}
\end{thm}

{We can follow the idea of deriving the estimate \eqref{202005201033} to establish Theorem \ref{201912041028}, the proof of
which will be presented in Section \ref{2020011sdf92326}. Here we explain why one has to modify $(\eta^0,u^0)$
to be initial data of $(\eta^{\mm{L}},u^{\mm{L}} )$ in \eqref{202001070914}$_4$, and why one imposes the additional regularity
\eqref{2020111222152} in order to get \eqref{202011032108}. }
\begin{enumerate}
  \item[(1)]
Since the initial data $u^{\mm{L}}|_{t=0}$ is divergence-free, i.e., $\mm{div}(u^{\mm{L}}|_{t=0})=0$, one has to adjust
the initial data $u^0$ as in \eqref{202001070914}$_4$.
  \item[(2)]
 { If the initial data $\eta^0$ of $\eta$ is directly used to be initial data of the corresponding linear problem, then
 we see that $\mm{div}\eta^{\mm{L}}=\mm{div}\eta^0$, and a time-decay of $\partial_2\eta^{\mm{d}}$ in \eqref{2020sfa05201033}
 can not be expected unless $\mm{div}\eta^0=0$. Hence, we have to modify $\eta^0$ as in \eqref{202001070914}$_4$,
  so that the obtained new initial data ``$\eta^0+\eta^{\mm{r}}$'' is also divergence-free. }
\item[(3)] Subtracting \eqref{202001070914} from \eqref{01dsaf16asdfasf}--\eqref{01dsaf16asdfasfsaf},
one obtains
\begin{equation}\label{01dsaf16safafasdfasfxx}
\begin{cases}
\eta_t^{\mm{d}}=u^{\mm{d}}, \\[1mm]
u_t^{\mm{d}}+\nabla q -\nu\Delta    u^{\mm{d}}- m^2
\partial_2^2\eta^{\mm{d}}=\mathfrak{N}, \\[1mm]
\div u^{\mm{d}}=-\mathrm{div}_{\tilde{\mathcal{A}}} {u},\\
(\eta^{\mm{d}},u^{\mm{d}})|_{t=0}=-(\eta^{\mm{r}},u^{\mm{r}}).
\end{cases}
\end{equation}
\quad
Let $\alpha=(\alpha_1,\alpha_2)$, an application of $\partial^\alpha$ with $|\alpha|=2$ to \eqref{01dsaf16safafasdfasfxx}$_2$ yields
\begin{equation}
\partial^\alpha u_t^{\mm{d}}+\partial^\alpha\nabla q -\nu\Delta \partial^\alpha  u^{\mm{d}}- m^2
\partial_2^2\partial^\alpha\eta^{\mm{d}}=\partial^\alpha\mathfrak{N}\qquad (|\alpha|=2).
\label{2020120717757}
\end{equation}
 We multiply \eqref{01dsaf16safafasdfasfxx}$_2$ by $\partial^\alpha\eta^{\mm{d}}$ and $\partial^\alpha u^{\mm{d}}$ in $L^2$,
respectively, and integrate to obtain the following two energy identities concerning derivatives of $(\eta^{\mm{d}},u^{\mm{d}} )$:
\begin{align}
&\frac{\mm{d}}{\mm{d}t}\int \left(\partial^\alpha\eta^{\mm{d}}\cdot \partial^\alpha u^{\mm{d}}+ \frac{\nu}{2}
|\partial^\alpha\nabla\eta^{\mm{d}}|^2\right)\mm{d}y+  \int |m \partial_2 \partial^\alpha \eta^{\mm{d}}|\mm{d}y-\int |\partial^\alpha u^{\mm{d}}|^2\mm{d}y \nonumber \\
&= \int \partial^\alpha\mathfrak{N}\cdot \partial^\alpha\eta^{\mm{d}} \mm{d}y
+\int \partial^\alpha q \partial^\alpha\mm{div}\eta^{\mm{d}} \mm{d}y=:\mathfrak{I}_1,  \label{20211292118}
\end{align}
and
\begin{align}
&\frac{1}{2}\frac{\mm{d}}{\mm{d}t} \int ( |\partial^\alpha u^{\mm{d}}|^2 +
|m\partial^\alpha\partial_2\eta^{\mm{d}}|^2)\mm{d}y +\nu\int |\nabla \partial^\alpha u^{\mm{d}}|^2\mm{d}y\nonumber \\
& = \int \partial^\alpha \mathfrak{N}\cdot \partial^\alpha u^{\mm{d}}
+\int  \partial^\alpha q \partial^\alpha\mm{div}u^{\mm{d}}\mm{d}y=:\mathfrak{I}_2.
\label{202011292118}
\end{align}
 {It seems that the integral terms $\mathfrak{I}_2$ for $\alpha_2\neq0$ and $\mathfrak{I}_1$ could provide the convergence rate $m^{-1}$
by directly employing the estimate $\|m\partial_2\eta\|_2\leqslant C$. } This idea, can not be directly applied to $\mathfrak{I}_2$ for $\alpha_2=0$
due to the integral term $\int\partial_1^3\eta_2\partial_2u \cdot\partial_1^3 u^{\mm{d}}\mm{d}y$ hidden in $\mathfrak{I}_2$.
 To circumvent this difficult term, we rewrite it as follows.
\begin{align}
\int \partial_1^3\eta_2  \partial_2u \cdot  \partial_1^3 u^{\mm{d}}\mm{d}y=
\int \partial_1^3\eta_2^{\mm{d}}  \partial_2u \cdot  \partial_1^3 u^{\mm{d}}\mm{d}y+\int \partial_1^3\eta_2^{\mm{L}}  \partial_2u \cdot  \partial_1^3 u^{\mm{d}}\mm{d}y. \label{20202111292107}
\end{align}
 {Noting that $\partial_1^3\eta_2^{\mm{d}}$ can provide a convergence rate $m^{-1}$ by using \eqref{20211292118} and the energy
identity \eqref{202011292118} with $|\alpha|=1$,} thus the first term on the right hand of \eqref{20202111292107} also provide a convergence rate $m^{-1}$. In addition, it is easy to formally derive from the odevity condition of $\eta_2^{\mm{L}}$ and the linear problem \eqref{202001070914} that
\begin{align}\|\partial_1^3\eta_2^{\mm{L}}\|_0\lesssim_0 \|\partial_1^3\partial_2\eta_2^{\mm{L}}\|_0
\leqslant C m^{-1}(1+\|\partial_1 \eta^0\|_3^2+\|\partial_1 u^0\|_2^2),
\label{202011292140}
\end{align}
thus we need the regularity condition \eqref{2020111222152} to make the derivation procedure of
  \eqref{202011292140} sense. Consequently, $\mathfrak{I}_2$ also implies the convergence rate $m^{-1}$ for $\alpha_2=0$
  under the additional condition \eqref{2020111222152}.
\end{enumerate}
\subsection{Extension to the inviscid case with damping}\label{202008091956}

We now describe how to extend the aforementioned results to the inviscid case with damping.
The equations of incompressible inviscid MHD fluids with zero resistivity and a low-order damping read as follows.
\begin{equation}\label{1.1xsfadax}
\begin{cases}
\rho  {v}_t+\rho {v}\cdot\nabla {v}+\nabla p- \kappa\rho v=\lambda M \cdot\nabla M/4\pi, \\
M_t+v\cdot\nabla M =M\cdot \nabla v,\\
\mathrm{div} v=\mathrm{div} M=0,
\end{cases}
\end{equation}
where $\kappa >0$ is the damping coefficient. We mention that the well-posedness for the idea MHD system with a velocity damping
has been widely investigated, see \cite{TZWY2013JDE,PRHZK2009JDE,WWKYT2001JDE,STCTBWDH2012CPDE} for examples.

 { Similarly to \eqref{01dsaf16asdfasf}, we rewrite \eqref{1.1xsfadax} in the following form of
Lagrangian coordinates:  }
\begin{equation}\label{01dsaf16asdfsfdsaasf}
                              \begin{cases}
\eta_t=u , \\[1mm]
 u_t+\nabla_{\ml{A}}q+ \kappa u=
  m^2 \partial_2^2\eta, \\[1mm]
\div_\ml{A}u=0,
\end{cases}
\end{equation}
with initial data
\begin{equation}\label{01dsafsfda16asdfasfsaf}
(\eta, u)|_{t=0}=(\eta^0,  u^0)\;\; \mbox{ in }\mathbb{T}^2.
\end{equation}
Then we have the following results, which can be regarded as an extension of Theorems \ref{201904301948}, \ref{202006211824saf} and \ref{201912041028}.
\begin{thm}
\label{2019043019sdfa48}
Let $\kappa$ be a positive constant. There are constants $ {c}_1^\kappa\geqslant 4$ and $ {c}_2^\kappa\in (0,1]$, such that
for any $(\eta^0,u^0)\in (\underline{H}^{5}_{1,2}\cap H^4_*)\times \underline{H}^4$ and $\kappa$ satisfying $\mm{div}_{\mathcal{A}^0}u^0=0$,
the odevity conditions \eqref{202007301500} and the condition of strong magnetic fields
\begin{equation}
\label{2019092sadf81832} m\geqslant \frac{1}{c_3^\kappa}\max\left\{ \mathfrak{K}^{1/2} ,\mathfrak{K} ^2\right\},
\end{equation} where ${\mathcal{A}^0}$ denotes the initial data of $\mathcal{A}$ and
$$
\mathfrak{K}:= \sqrt{ {c^\kappa_1} \|(\eta^0,u^0, m \partial_2\eta^0) \|_4^2e^{c^\kappa_2\| (\eta^0,u^0,m \partial_2\eta^0)\|_3^2 }},$$
 the initial value problem \eqref{01dsaf16asdfsfdsaasf}--\eqref{01dsafsfda16asdfasfsaf} admits a unique global classical solution
 $(\eta, u,q)\in\underline{\mathfrak{H}}^{5,*}_{1,2,\infty}\times\underline{\mathfrak{U}}_\infty^4\times ({C}^0(\mathbb{R}^+_0,\underline{H}^3)
 \cap \mathfrak{C}^0(\mathbb{R}^+,H^4)) $.
Moreover, the solution $(\eta,u)$ enjoys
\begin{enumerate}[(1)]
  \item Decay estimate:
\begin{align}
 &\sup_{t\geqslant 0}\big\{ e^{  c_3^\kappa \min\{1,m\} t}(\| (\eta,u)(t)\|_4^2+\| m\eta (t)\|_{5,2}^2)\big\}  \nonumber \\
& \qquad +\int_0^\infty e^{c_3^\kappa\min\{1,m\} \tau}(\| u\|_4^2+\| m \eta \|_{5,2}^2) \mm{d}\tau \lesssim_\kappa  \mathfrak{K}^2
  \label{202008071638}
\end{align}
and
\begin{align}
& \sum_{0\leqslant i\leqslant 1}\left(\sup_{t\geqslant 0}\big\{ \langle t\rangle^i(\|(\eta,u)(t)\|_4^2 +\|m\eta (t)\|_{5,2}^2)\big\}
+ \int_0^\infty\langle  \tau\rangle^i( \|  u\|_4^2+\|m\eta \|_{5,2}^2)\mm{d}\tau\right)\nonumber \\
& \qquad + \int_0^t\langle  \tau\rangle \|\eta \|_{5,2}^2\mm{d}\tau\leqslant C^\kappa ,
  \label{202911sfa9fasdfd52245}
\end{align}
where we remark that the decay rates in \eqref{202911sfa9fasdfd52245} do not depend on $m$ for fixed
$\|(\eta^0,u^0, m\partial_2\eta^0)\|_4^2$.
\item Stability around the solution $(\eta^{\mm{L}},u^{\mm{L}})$ of the linear problem:
 \begin{align}
&  \sup_{t\geqslant 0}\big\{e^{c_3^\kappa\min\{1,m\} t}(\|(\eta^{\mm{d}},u^{\mm{d}})(t)\|_4^2+\| m\eta^{\mm{d}}(t)\|_{5,2}^2)\big\} \nonumber \\
& \qquad +\int_0^\infty e^{ c_3^\kappa  \min\{1,m\} \tau }(\| u^{\mm{d}} \|_4^2
+ \| m \eta^{\mm{d}}\|_{5,2}^2) \mm{d}\tau \leqslant C^\kappa m^{-2} \label{2020212052001}
\end{align}
and
 \begin{align}
& \sup_{t\geqslant 0}\big\{ \langle t\rangle^i(\|(\eta^{\mm{d}}, u^{\mm{d}})(t)\|_4^2+\|m\eta^{\mm{d}}(t)\|_{5,2}^2)\big\} \nonumber \\
& \qquad  +\int_0^\infty\langle  \tau\rangle^i (\| u^{\mm{d}}\|_{4}^2
+(1+m^2)\| \eta^{\mm{d}}\|_{5,2}^2 )\mm{d}\tau \leqslant C^\kappa  m^{-2},\quad\; i=1,2,\label{202012052058}
\end{align}
\end{enumerate}
where the error function $(\eta^{\mm{d}},u^{\mm{d}}):=(\eta-\eta^{\mm{L}},u-u^{\mm{L}})$ and $(\eta^{\mm{L}},u^{\mm{L}})\in \mathfrak{C}^0(\mathbb{R}^+,{H}^5_2)\times\underline{\mathfrak{U}}_\infty^4 $ is the unique classical solution of the following linear
pressureless initial-value problem:
\begin{equation}
\label{01dsaf16asdfsfsdafafdsaasf}
\begin{cases}
\eta_t^{\mm{L}}=u^{\mm{L}} , \\[1mm]
 u_t^{\mm{L}}  - \kappa u^{\mm{L}}=
  m^2 \partial_2^2\eta^{\mm{L}}, \\[1mm]
\div u^{\mm{L}}=0,\\
(\eta^{\mm{L}},u^{\mm{L}})|_{t=0}=(\eta^0+\eta^{\mm{r}}, u^0+u^{\mm{r}}),
\end{cases}
\end{equation}
for some $(\eta^{\mm{r}}, u^{\mm{r}})\in  \underline{H}^5\times \underline{H}^4$ satisfying  \eqref{202011032122} and
\begin{align}
&\|u^{\mm{r}} \|_4\lesssim_0 \|\partial_2 \eta^0\|_4\|u^0\|_4\;\mbox{ and }\;
\|\partial_2^i\eta^{\mm{r}}\|_{4}\lesssim_0\|\partial_2\eta^0\|_{3+i}\|\nabla\partial_2^i\eta^0\|_{3},\;\;\; i=0,1.
\label{20asdff2011032123}
\end{align}
\end{thm}
\begin{rem}
We have assumed $(\eta^0)_{\mm{T}^2}=(u^0)_{\mm{T}^2}=0$ in Theorem \ref{2019043019sdfa48}. If $((\eta^0)_{\mm{T}^2},(u^0)_{\mm{T}^2})\neq0$,
then the unique solution $(\eta,u,q)$ of \eqref{01dsaf16asdfsfdsaasf}--\eqref{01dsafsfda16asdfasfsaf} enjoys the following form:
$$ (\eta,u,q)=(\bar{\eta}+\kappa^{-1}(u^0)_{\mathbb{T}^2}(1-e^{-\kappa t})+ (\eta^0 )_{\mathbb{T}^2},  \bar{u}+(u^0)_{\mathbb{T}^2}e^{-t},\bar{q}),$$
where $(\bar{\eta},\bar{u},\bar{q})$ is the unique solution, established in Theorem \ref{2019043019sdfa48}, of the following initial-value problem
\begin{equation}
\begin{cases}
\bar{\eta}_t=\bar{u} , \\[1mm]
\bar{u}_t+\nabla_{\bar{\ml{A}}}\bar{q}+\kappa\bar{u}=
m^2 \partial_2^2\bar{\eta},\\[1mm]
\div_{\bar{\ml{A}}}\bar{u}=0,
\end{cases}
\nonumber
\end{equation}
with initial data
$(\bar{\eta}, \bar{u})|_{t=0}=(\eta^0-(\eta^0)_{\mathbb{T}^2},  u^0-(u^0)_{\mathbb{T}^2})$ in $\mathbb{T}^2$.
\end{rem}

The key idea in the proof of Theorem \ref{2019043019sdfa48} is similar to that in the proof of Theorems \ref{202006211824saf}--\ref{201912041028}.
 {But there are remarkable differences between the decay results for the inviscid and viscous cases, which
will be further explained in the proof process of Theorem \ref{2019043019sdfa48},} see all the footnotes in Section \ref{201912asdfas062021}.

\subsection{Verification of preserving the odevity of solutions}
We end this section by verifying the assertion in Remark \ref{202008191942}.  To this end, let $(\eta,u,q)$ be a strong solution of
\eqref{01dsaf16asdfasf}--\eqref{01dsaf16asdfasfsaf}, $\psi=(\eta_1,-\eta_2)(y_1,-y_2,t )$, $w=( u_1,-u_2)(y_1,-y_2,t)$ and
$p=( q_1, q_2)(y_1,-y_2,t)$. {Due to the uniqueness of strong solutions, we see that
 to get the desired conclusion, it suffices to verify that $(\psi,w)$ is also a strong solution of \eqref{01dsaf16asdfasf}. }
 It is obvious that $(\psi,w)$ satisfies \eqref{01dsaf16asdfasf}$_1$. Next, we show that $(\psi,w)$ also satisfies \eqref{01dsaf16asdfasf}$_2$
 and \eqref{01dsaf16asdfasf}$_3$.

Defining \begin{align} {\mathcal{A}}=\left(\begin{array}{cc}
            \partial_2\eta_2+1 &  - \partial_1\eta_2\\
                - \partial_2\eta_1&   \partial_1\eta_1+1
                 \end{array}\right)\mbox{ and } \mathcal{B}=\left(\begin{array}{cc}
            \partial_2\psi_2+1 &  - \partial_1\psi_2\\
                - \partial_2\psi_1&   \partial_1\psi_1+1
                 \end{array}\right), \nonumber
\end{align}
then
\begin{align}
& (\mathcal{B}_{11} ,\mathcal{B}_{22})= ({\mathcal{A}}_{11} ,{\mathcal{A}}_{22})|_{y_2=-y_2 }\mbox{ and } (\mathcal{B}_{12}, \mathcal{B}_{21})= -({\mathcal{A}}_{12} ,{\mathcal{A}}_{21})|_{y_2=-y_2 }.
\label{202002122032}
\end{align}

 Let $\chi $ be a function, and $\tilde{\chi}=\chi(y_1,-y_2) $. Then, by \eqref{202002122032},
\begin{align}\nabla_{{\ml{B}}} \tilde{\chi}   = ({\ml{A}}_{1i} ,- {\ml{A}}_{2i})^{\mm{T}}  \partial_i \chi |_{y_2=-y_2 }
\label{20200212203safdsaf2}
\end{align}
 and
$$ \mm{div}_{\ml{B}} w = \mm{div}_{ {\ml{A}}}u |_{y_2=-y_2 }=0, $$
from which we see that $(\eta,u)$ satisfies \eqref{01dsaf16asdfasf}$_3$.

By \eqref{202002122032} and \eqref{20200212203safdsaf2}, we further have
\begin{align}\Delta_{ {\ml{B}}}\tilde{\chi} =\Delta_{\ml{A}} {\chi} |_{y_2=-y_2 }.\label{2020safa02122032}
\end{align}
Thus, the identities \eqref{20200212203safdsaf2} and \eqref{2020safa02122032} yield
$$
\begin{aligned}
  \partial_tw_i +\nabla_{ {\mathcal{B}}}p - \nu\Delta_{ {\ml{B}}}w_i =
      \begin{cases}
   (\partial_t u_1+ {\mathcal{A}}_{1i}\partial_i q- \nu\Delta_{{\ml{A}}}u_1 )|_{y_2=-y_2 }  & \hbox{ for }i=1, \\
-  (\partial_t u_2+{\mathcal{A}}_{2i}\partial_i q- \nu \Delta_{{\ml{A}}}u_2 )|_{y_2=-y_2 }  & \hbox{ for }i=2,
      \end{cases}
\end{aligned}$$
and
$$ \partial_2^2 (\psi_1,\psi_2)=  \partial_2^2 (\eta_1,-\eta_2)|_{y_2=-y_2 } .$$
Hence, we see that $(\psi,w)$ satisfies \eqref{01dsaf16asdfasf}$_2$ by the above two identity.  This proves the property of preserving
the odevity of strong solutions.

\section{Proof of Theorem \ref{201904301948} } \label{20190sdafs4301948}

This section is devoted to the proof of Theorem \ref{201904301948}. First we derive some basic (\emph{a priori}) estimate for $(\eta,u,q)$
under the \emph{a priori} assumption \eqref{aprpiosesnewxxxx} associated with the smallness condition \eqref{aprpiose1snewxxxxz} in Subsection
\ref{202061410014}, then further establish the stability estimate \eqref{20190safd5041053} in Subsection \ref{202006141016}, and finally,
 introduce a local well-posedness result for \eqref{01dsaf16asdfasf}--\eqref{01dsaf16asdfasfsaf} and complete the proof of
 Theorem \ref{201904301948} by a standard continuity argument in Subsection \ref{20206141015}.

\subsection{Energy estimates} \label{202061410014}
Let $(\eta,u,q)$ be a solution
of the initial-value problem \eqref{01dsaf16asdfasf}--\eqref{01dsaf16asdfasfsaf} defined on $\Omega_T$  for any given $T>0$,
where $(\eta^0,u^0) $ belongs to $\underline{H}^{3}_{1} \times \underline{H}^2$, and satisfies  $\mm{div}_{\mathcal{A}^0}u^0=0$
and the odevity conditions \eqref{202007301500}. We recall here that the solution automatically satisfies  $(\eta)_{\mathbb{T}^2}=(u)_{\mathbb{T}^2}=0$ and the odevity conditions in \eqref{20208192055} if $(\eta^0,u^0)$ does.
We further assume that $(\eta,u,q)$ and $K$
satisfy  $(q)_{\mathbb{T}^2}=0$, \eqref{aprpiosesnewxxxx} and \eqref{aprpiose1snewxxxxz},
  where $K\geqslant 1$ will be defined by \eqref{201911262060} and $\delta\in (0,1]$ is a sufficiently small constant.
  It should be noted that the smallness of $\delta$ only depends on the parameter $\nu$. In addition, by Young's inequality,
  one easily finds from \eqref{aprpiosesnewxxxx} and \eqref{aprpiose1snewxxxxz} that
\begin{equation}\label{aprpiosesnewxxxx1x}
 \sup_{0\leqslant t \leqslant T} \|\partial_2\eta(t)\|_2  \leqslant \frac{K}{m}\leqslant \frac{1}{3m}(2K^{1/2}+ K^2)\leqslant \delta.
\end{equation}

Before deriving the energy estimates for $(\eta,u)$, we introduce some basic inequalities and establish some preliminary estimates of $(\eta,u)$
by using the following two lemmas.
\begin{lem}\label{202005091306}
 We have the following basic inequalities:
\begin{enumerate}[(1)]
   \item
 Generalized Poinc\'are's inequalities:
\begin{align}
& \| f\|_0 \lesssim_0 \|\partial_2 f\|_0\;\;\mbox{ for any }f\in \overline{H}^1\mbox{ satisfying }\int_{-1}^1 f(y_1,y_2)\mm{d}y_2=0  , \label{201005041609}\\
&\label{20160614fdsa19asfda57x}
 \| f\|_i \leqslant \tilde{c}_i \|\nabla^i f\|_0\;\;\mbox{ for any }f\in \underline{H}^i ,
\end{align}
where  and the positive constant $\tilde{c}_i$ only depends on $i$.
 \item Interpolation inequality: for any $f\in H^2_2$,
\begin{align}
&\|f\|_{L^\infty }^2 \lesssim_0  \|f\|_0\|(f,\partial_1 f)\|_0+ \|f\|_0^{1/2} \|(f,\partial_1 f)\|_0^{1/2}\|\partial_2 f \|_{0}^{1/2} \|\partial_2(f,\partial_1f) \|_{0}^{1/2}.
\label{202004221saffad412}
\end{align}
In particularly, by Young's inequality,
\begin{align}
&\|f\|_{L^\infty } \lesssim_0    \|\partial_2(f,\partial_1f) \|_{0} +\|(f,\partial_1 f)\|_0 .
\label{202004221412}
\end{align}
\item Product estimate: for any $(f, g ,h)\in H^1_2\times {H}  \times L^2$,
\begin{equation}
\int | fgh|\mm{d}y \lesssim_0 \sqrt{\|f\|_{0}\|f\|_{1,2}}  \|(g,\partial_1 g)\|_0\|h\|_0
\lesssim_0 \|(f,\partial_2 f)\|_0 \|(g,\partial_1 g)\|_0\|h\|_0. \label{20200508}
\end{equation}
\item There is a constant $\delta\in (0,1]$, such that for any $(\xi-y)\in H$ satisfying $\|\nabla (\xi-y)\|_{L^\infty}\leqslant \delta$,
\begin{align}
\|\nabla_{\mathcal{B}}f\|_0\lesssim_0  \| \nabla f\|_0\lesssim_0 \|\nabla_{{\mathcal{B}}} f\|_0\mbox{ for any }\nabla f\in L^2,
\label{202010261646}
\end{align}
where $\ml{B}:=(\nabla \xi+I)^{-\mm{T}}$.
\end{enumerate}
\end{lem}
\begin{pf}
(1) The inequalities \eqref{201005041609} and \eqref{20160614fdsa19asfda57x} are obvious to get
by virtue of classical Poinc\'are's inequality:
  $$\left\|g\right\|_{L^r(\Omega)}^r \leqslant  \tilde{c}_{r,\Omega}\left(\|\nabla g\|_{L^r(\Omega)}^r+
 \left|\int_\Omega g\mm{d}y\right|^r \right)\mbox{ for any } g\in W^{1,r}(\Omega),$$
where $\Omega=\mathbb{T}^n$ with $n\geqslant 1$, $r\geqslant 1$ is a constant and the constant $\tilde{c}_{r,\Omega}$ depends only on $r$ and $\Omega$.

(2) Now, we turn to the derivation of \eqref{202004221saffad412}. Here and in what follows, we denote
$$\|f\|_{L^p_{y_i}}:= \left(\int_{-1}^1 |f|^p\mm{d}y_i\right)^{1/p}\;\; \mbox{ for }f\in L^p.$$

Since $C^2(\mathbb{T}^2)$ is dense in $H^2_2$, it suffices to prove \eqref{202004221saffad412} for $f\in C^2(\mathbb{T}^2)$.
Noting that, for any $g\in C^0(\mathbb{T}^2)$,
$$ \sup_{y_1\in(-1,1) } g(y_1,y_2)\mbox{ is a measurable function defined on }(-1,1),$$
we use the Fubini theorem and the one-dimensional interpolation inequality (see \cite[Theorem]{NLOE}) to deduce that
\begin{align}
\sup_{s\in(-1,1) }|\phi(s)|^2\lesssim_0  \|\phi(s)\|_{L^2(-1,1)}\| (\phi,\phi')(s)\|_{L^2(-1,1)}\;\;\mbox{ for }\phi\in H^1(-1,1). \label{202008021649}
\end{align}
Therefore,
$$
\begin{aligned}
\|f\|_{L^\infty }^2= &\sup_{y_1\in(-1,1) }\left(\sup_{y_2\in(-1,1) }|f(y_1,y_2)|^2\right)& \\
\lesssim_0 &
\sup_{y_1\in(-1,1) }\left( \| f\|_{L^2_{y_2} }\|\partial_2 f \|_{L^2_{y_2} }+\| f\|_{L^2_{y_2} }^2\right)  \\
\leqslant  &  \left\|\sup_{y_1\in(-1,1) }
|f|^2\right\|_{L^1_{y_2} }^{1/2}\left\|\sup_{y_1\in(-1,1) }|\partial_2 f|^2 \right\|_{L^1_{y_2} }^{1/2}+\left\|
\sup_{y_1\in(-1,1) }|f|^2\right\|_{L^1_{y_2} }   \\
\lesssim_0   &\|f\|_0^{1/2} \|(f,\partial_1 f)\|_0^{1/2}\|\partial_2 f \|_{0}^{1/2}
\|\partial_2(f,\partial_1f) \|_{0}^{1/2}+\|f\|_0\|(f,\partial_1 f)\|_0,
\end{aligned}
$$
which gives \eqref{202004221saffad412}.

(3) Finally, we prove \eqref{20200508}. Recalling that $C^1(\mathbb{T}^2)$ is dense in $H_2^1$ and $H$, it suffices to prove
\eqref{20200508} for $(f,g)\in C^1(\mathbb{T}^2)$. Let $(f,g)\in C^1(\mathbb{T}^2)$, we have
 $$
\left\|\sup_{y_2\in (-1,1)}|f |\right\|_{L^2_{y_1}} \lesssim   \sqrt{\|f\|_0\|f\|_{1,2}},\quad
\sup_{y_1\in (-1,1)}\left\| g \right\|_{L^2_{y_2}}  \lesssim    \sqrt{\|g\|_0\|(g,\partial_1 g)\|_0},
$$
which, together with the Fubini theorem and H\"older's inequality, implies
 \begin{align}
\int \left| f g h\right|\mm{d}y\lesssim_0 &
\int_{\mathbb{T} }\sup_{y_2\in (-1,1)}|f| \int_{\mathbb{T} }|gh|\mm{d}y_2\mm{d}y_1  \nonumber \\
\lesssim_0&\sup_{y_1\in (-1,1)} \|g\|_{L^2_{y_2}}
\int_{\mathbb{T}} \sup_{y_2\in (-1,1)}|f |   \|h\|_{L^2_{y_2}}\mm{d}y_1\nonumber \\
\lesssim_0&\left\|\sup_{y_2\in (-1,1)}|f |\right\|_{L^2_{y_1}}\sup_{y_1\in (-1,1)} \|g\|_{L^2_{y_2}}\|h\|_0 \nonumber \\
\lesssim_0 &  \sqrt{\|f\|_0\|  f\|_{1,2}}    \|_{0}\|(g,\partial_1 g)\|_0\|h\|_0. \nonumber
\end{align}
Therefore, an application of Young's inequality yields \eqref{20200508}.

(4) The equivalent estimate \eqref{202010261646} holds obviously.
\hfill$\Box$
\end{pf}
 \begin{lem}
\label{201805141072dsafa}
Under the condition \eqref{aprpiosesnewxxxx1x} with sufficiently small  $\delta\in (0,1]$, we have
 \begin{enumerate}[(1)]
 \item  Estimates for $\eta$:
\begin{align}
& \| \partial_i \eta_j \|_{L^\infty}\lesssim _0
                                           \begin{cases}
                                   \|\partial_2\eta\|_2   & \hbox{ for }i=1\mbox{ and }1\leqslant j\leqslant 2, \\
  \|\partial_2^2\eta\|_1   & \hbox{ for }i=2, \\
                                                  \end{cases}
\label{2020043013014}\\
& \|\eta\|_{i+1,2}\lesssim_0 \|\partial_2\eta\|_i \hbox{ for }0\leqslant i\leqslant 2, \label{20202111191231}\\
& \|\eta\|_{L^\infty} \lesssim_0 \sqrt{\| \partial_2\eta\|_0\|\partial_2\eta\|_1}. \label{2020112af12454}
\end{align}
 \item  Estimates for $u$:
\begin{align}
&\|u\|_{i+1,2}\lesssim_0
                     \begin{cases}
                   \|\partial_2 u\|_0+\|\partial_2^2\eta\|_1\|\partial_2u\|_1  \ & \hbox{ for }i=0, \\
                    \|\partial_2 u\|_i & \hbox{ for }i=1\mbox{ and }2,
                     \end{cases}\label{202006121323}\\
&\|\partial_i u \|_{L^\infty}\lesssim_0
\begin{cases}
\|\partial_2u\|_2   & \hbox{ for }i=1, \\
  \|\partial_2^2u\|_1  & \hbox{ for }i=2.
                            \end{cases}  \label{202006101611}
\end{align}
\item  Product estimates:
\begin{align}
&   \int | \partial^\alpha\eta gh|\mm{d}y \lesssim_0
                                                      \begin{cases}
                                       \| \partial_2 \eta\|_{| \alpha|} \|g\|_1\|h\|_0 & \hbox{ for }\alpha_2= 0, \\
    \|\partial_2\partial^\alpha \eta\|_{0} \|g\|_1\|h\|_0 &  \hbox{ for }\alpha_2\neq 0,
                                                      \end{cases}
\label{202safas00508}\\[1mm]
&   \int | \partial^\beta u gh|\mm{d}y \lesssim_0
                                                      \begin{cases}
                                       \| \partial_2 u\|_{| \beta|} \|g\|_1\|h\|_0 & \hbox{ for }\beta_2= 0, \\
    \|\partial_2 \partial^\beta u\|_{0} \|g\|_1\|h\|_0 &  \hbox{ for }\beta_2\neq 0
                                                      \end{cases}\label{202safas00508x}
\end{align}
for any $(g,h)\in   H \times L^2$, where $\alpha$ and $\beta$  satisfy  $1\leqslant   |\alpha|\leqslant 2$ and  $|\beta|= 2$.
\end{enumerate}
 \end{lem}
\begin{rem}
Thanks to \eqref{aprpiosesnewxxxx1x} and \eqref{2020043013014}, one has
 \begin{align}
\|\nabla \eta\|_{L^\infty}\lesssim_0 \|\partial_2\eta\|_2\lesssim_0 \delta,  \label{20206111731}
\end{align}
which, together with \eqref{202010261646}, yields the following equivalent estimate: for sufficiently small $\delta$, it holds that
\begin{align}
  & \|\nabla_{\mathcal{A}}f\|_0\lesssim_0  \| \nabla f\|_0\lesssim_0 \|\nabla_{{\mathcal{A}}} f\|_0\;\;\mbox{ for any }\nabla f\in L^2.   \label{20200501}
 \end{align}
\end{rem}
\begin{pf}
(1) By \eqref{202004221412} we have that for $1\leqslant i, j\leqslant 2$,
\begin{align}
\|  \partial_i \eta_j \|_{L^\infty}\lesssim_0
\|\partial_2\partial_i( \eta_j,\partial_1 \eta_j) \|_{0} +\|\partial_i(\eta_j,\partial_1 \eta_j)\|_0. \label{202005051503}
\end{align}

Since $ \eta_2(y_1,y_2)$ is a odd function with respect to $y_2$ for any given $y_1$, one sees that
\begin{align}
\int_{-1}^1\partial_1^k\eta_2\mm{d}y_2=0\;\;\mbox{ for }0\leqslant k\leqslant 2. \nonumber
\end{align}
In addition, $ \eta_2(y_1,y_2)$ is a periodic function with respect to $y_2$ for any given $y_1$, then
\begin{align}
\int_{-1}^1 \partial_2^l\partial_1^k\eta_2\mm{d}y_2=0\;\;\mbox{ for }1\leqslant k+l\leqslant 3\mbox{ and }1\leqslant l.  \nonumber
\end{align}

Using \eqref{201005041609} and the above two relations, we find that
\begin{align}
\|\partial_1^k\eta_2\|_0\lesssim_0 \|\partial_2\partial_1^k\eta_2\|_0\;\;\mbox{ for }0\leqslant k\leqslant 2
\label{20d2005041539}
\end{align}
and
\begin{align}
 \|\partial_2^l\partial_1^k\eta_j\|_0\lesssim_0\|\partial_2^{l+1}\partial_1^k \eta_j\|_{0}
 \;\;\mbox{ for }0\leqslant l+k\leqslant 2\mbox{ and }1\leqslant l. \label{202005121605}
\end{align}
Putting \eqref{202005051503}--\eqref{202005121605} together, we get
\begin{align}
& \|  \partial_i \eta_j  \|_{L^\infty}\lesssim _0
                                                  \begin{cases}
                                                     \|\partial_2\eta_j\|_2   & \hbox{ for }i=1\mbox{ and }j=2, \\
   \|\partial_2^2\eta_j\|_1   & \hbox{ for }i=2\mbox{ and }j=1,\ 2.
                                                  \end{cases}
\label{202004301301sfa4}
\end{align}
To obtain \eqref{2020043013014}, we next show that \eqref{202004301301sfa4} also holds for $(i,j)=(1,1)$.

Keeping in mind that $\det(I+\nabla \eta)=1$, one gets from Sarrus' rule that
\begin{align}
\mathrm{div}\eta=&\partial_1\eta_2\partial_2 \eta_1
-\partial_1\eta_1\partial_2 \eta_2 . \label{201903211118}
\end{align}
In particular,
\begin{align}
\partial_1\eta_1=\partial_1\eta_2\partial_2 \eta_1
-\partial_1\eta_1\partial_2 \eta_2-\partial_2\eta_2.\label{202005201201}
\end{align}
Applying $\nabla$ to \eqref{202005201201}, and then multiplying the resulting identity by $\nabla \partial_1\eta_1$ in $L^2$, we obtain
\begin{align}
\|\nabla \partial_1\eta_1\|_0^2=\int \nabla   (\partial_1\eta_2\partial_2 \eta_1
-\partial_1\eta_1\partial_2 \eta_2-\partial_2\eta_2)\cdot \nabla \partial_1 \eta_1\mm{d}y. \label{2020061534}
\end{align}

Making use of H\"older's inequality, \eqref{20160614fdsa19asfda57x}, \eqref{20200508}, \eqref{20d2005041539} and \eqref{202004301301sfa4}, we infer
from \eqref{2020061534} that
\begin{align}
\|  \partial_1\eta_1\|_1^2\lesssim_0 & \|\partial_1\eta_1\|_1(\|  \partial_2 \eta_2\|_1
+ \|\partial_2 \eta_2\|_{L^\infty} \|\partial_1\eta_1\|_1+ \|\partial_1 \eta_2\|_{L^\infty} \|   \partial_2\eta_1\|_1) \nonumber \\
& +\int |\nabla    \partial_1\eta_2||\partial_2 \eta_1| |\nabla \partial_1 \eta_1|\mm{d}y+ \int |\partial_1\eta_1 ||\nabla \partial_2 \eta_2|
 |\nabla \partial_1 \eta_1|\mm{d}y   \nonumber \\
\lesssim_0 & \|\partial_1\eta_1\|_1(\|  \partial_2 \eta_2\|_1 + \|\partial_1 \eta_1\|_1 \|\partial_2 \eta_2\|_{2} +(\|\partial_1 \eta_2\|_{L^\infty}
+ \|\partial_2\eta_2 \|_2)\|   \partial_2\eta_1\|_1). \nonumber
\end{align}
By \eqref{aprpiosesnewxxxx1x}, \eqref{202004301301sfa4} and Young's inequality, we further deduce from the above estimate that
\begin{align}
\label{202006101449}
\|  \partial_1\eta_1\|_1 \lesssim_0 \|  \partial_2 \eta\|_1 .
\end{align}
We immediately see from \eqref{202005051503} and \eqref{202006101449} that $\|\partial_1 \eta_1  \|_{L^\infty}\lesssim _0                                                      \|\partial_2\eta_1\|_2 $. This completes the proof of \eqref{2020043013014}.

Thanks to \eqref{2020043013014}, we can get from \eqref{202005201201} that
\begin{align}\label{202006sdfa101449}\|  \partial_1\eta_1\|_0 \lesssim_0 \|  \partial_2 \eta\|_0.\end{align}
Putting  \eqref{20d2005041539},  \eqref{202006101449} and \eqref{202006sdfa101449} together, one concludes
\begin{align}
&\| \nabla \eta  \|_{i-1}\lesssim_0 \| \partial_2 \eta\|_{i}\;\;\mbox{ for }1\leqslant i\leqslant 2.\label{202004301335sdfa}
\end{align}
Thus, from \eqref{20160614fdsa19asfda57x}, \eqref{20d2005041539}, \eqref{202006sdfa101449} and \eqref{202004301335sdfa}, the estimate \eqref{20202111191231} follows immediately.
Finally, the estimate \eqref{2020112af12454} is obvious by virtue of \eqref{202004221saffad412} and \eqref{20202111191231}.

(2) Since $u$ enjoys the same odevity and the same periodicity as $\eta$ does,
we obtain, similarly to \eqref{202005051503}--\eqref{202004301301sfa4}, that
\begin{align}
& \|\partial_1^ku_2\|_0\lesssim_0 \|\partial_2\partial_1^ku_2\|_0\;\mbox{ for }0\leqslant k\leqslant 2, \label{20d20050asda41539} \\
& \|\partial_2^l\partial_1^ku_j \|_0\lesssim_0 \|\partial_2^{l+1} \partial_1^k u_j  \|_{0}\;\mbox{ for }0
\leqslant l+k\leqslant 2\mbox{ and }1\leqslant l, \label{202005121605sdfa}\\
& \|\partial_1 u_1\|_{L^\infty} \lesssim_0 \|\partial_2\partial_1 (u_1,\partial_1 u_1)\|_{0} +\|\partial_1(u_1,\partial_1 u_1)\|_0
\label{202005081107}
\end{align}
and
\begin{align}
& \|  \partial_i u_j  \|_{L^\infty}\lesssim _0
                                                  \begin{cases}
                                                     \|\partial_2u_j\|_2   & \hbox{ for }i=1\mbox{ and }j=2, \\
  \|\partial_2^2u_j\|_1   & \hbox{ for }i=2\mbox{ and }j=1,\ 2.
                                                  \end{cases}
\label{2020043013sfda01sfa4}
\end{align}

By \eqref{01dsaf16asdfasf}$_3$ and the definition of $\tilde{\mathcal{A}}$, we find that
\begin{align}
 \mm{div}u=- \mm{div}_{\tilde{\mathcal{A}}}{u},
\label{202005011525}
\end{align}
where
\begin{align}   \tilde{\mathcal{A}}=\left(\begin{array}{cc}
            \partial_2\eta_2 &  - \partial_1\eta_2\\
                - \partial_2\eta_1&   \partial_1\eta_1
                 \end{array}\right). \label{20206101609}
\end{align}

If we apply the norm $\|\cdot\|_0$ to \eqref{202005011525}, we get
$$\|\partial_1 u_1\|_0\lesssim_0 \|\partial_2 u \|_0 +\|\partial_2\eta_1\|_{L^\infty} \|\partial_1 u_2\|_0+ \delta\|\partial_1 u_1\|_0,$$
which, together with \eqref{2020043013014} and \eqref{20d20050asda41539}, gives
\begin{align}
\label{202010241545}
\|\partial_1 u_1\|_0\lesssim_0 \|\partial_2 u \|_0+ \|\partial_2^2 \eta \|_1\|\partial_2 u \|_1
\end{align}
 for sufficiently small $\delta$.

In addition, following the same process as in the derivation of \eqref{202006101449}, we easily deduce from \eqref{202005011525} that
\begin{align}
\|\partial_1 u_1 \|_1\lesssim_0 &
\|\partial_2u \|_2 .   \label{202005092249}
\end{align}
Thus, the estimate \eqref{202006121323} follows from \eqref{aprpiosesnewxxxx1x}, \eqref{20160614fdsa19asfda57x}, \eqref{20d20050asda41539},
\eqref{202010241545} and \eqref{202005092249}.

Plugging \eqref{202005092249} into \eqref{202005081107}, we get $\|\partial_1 u_1\|_{L^\infty}\lesssim \|\partial_2 u\|_2$,
which, together with \eqref{2020043013sfda01sfa4}, yields \eqref{202006101611}.

(3) The estimate \eqref{202safas00508} follows from \eqref{20200508},  \eqref{202005121605} and \eqref{202004301335sdfa}.
Similarly, in view of \eqref{20200508},    \eqref{20d20050asda41539},  \eqref{202005121605sdfa} and \eqref{202005092249},
we see that \eqref{202safas00508x} holds. This completes the proof of Lemma \ref{201805141072dsafa}.
\hfill $\Box$
\end{pf}

Now, we proceed to derive some basic energy estimates for $(\eta,u,q)$.
\begin{lem}\label{201612132242nx}
Under the condition \eqref{aprpiosesnewxxxx1x} with sufficiently small $\delta$, we have
\begin{align}
&\|\nabla q\|_i \lesssim_0
\begin{cases}
\|(u,\partial_2 u)\|_0\|\nabla u\|_1+  m^2\|\partial_2^2\eta \|_0\|\partial_2\eta\|_2    & \hbox{ for }i=0 ; \\
 \|\partial_2^2 u\|_0\|\nabla u\|_1\;\mbox{ or }\; \|\partial_2 u\|_1^2  & \\
  +m^2 (   \|\partial_2^3 \eta\|_0 \|\partial_2 \eta\|_2 +\|\partial_2^2 \eta\|_1^2 ) & \hbox{ for }i=1,
 \end{cases} \label{2017020614181721}  \\[1mm]
 & \|\nabla q_t\|_1 \lesssim_0 \|\partial_2 u\|_2 \|\nabla u\|_2  (1+\| \nabla \eta\|_2+\|\nabla u\|_1 )
   + m^2\|\partial_2\eta\|_{2}( \|\partial_2\eta\|_2\|\nabla u\|_2 +\|\partial_2u\|_2).\label{2020103safda02130}
\end{align}
\end{lem}
\begin{rem} The  estimate  \eqref{2020103safda02130} will be used in the derivation of the error estimate \eqref{202011032108} in  Theorem \ref{201912041028}, see the last term in \eqref{202011220947}.
\end{rem}
\begin{pf}
(1) By \eqref{201909261909} and \eqref{01dsaf16asdfasf}$_3$, we have
\begin{align} \mm{div}_{\mathcal{A}} u_t = -\mm{div}_{\mathcal{A}_t} u = -\mm{div}({\mathcal{A}}_t^{\mm{T}} u).  \label{202005101008}
\end{align}
Multiplying \eqref{01dsaf16asdfasf}$_2$ by $\nabla_{\mathcal{A}} q$ in $L^2$ and using \eqref{201903211118}, \eqref{202005101008}, we
integrate by parts and recall the fact
\begin{align}
\mm{div}_{\mathcal{A}}\Delta_{\mathcal{A}}u=0 \label{202011211548}
\end{align}
to infer that
\begin{align}
\|\nabla_{\mathcal{A}}q\|_0^2 =
&    \int ( {\mathcal{A}}_t^{\mm{T}} u)\cdot \nabla q\mm{d}y+
m^2 \int \partial_2^2 \eta \cdot \nabla_{\tilde{\mathcal{A}}}q\mm{d}y\nonumber \\
& +m^2 \int \partial_2(\partial_1\eta_2\partial_2 \eta_1
-\partial_1\eta_1\partial_2 \eta_2 )   \partial_2 q\mm{d}y.\label{n0101nn928}\end{align}
 Exploiting \eqref{20200508} and \eqref{202safas00508}, we deduce from the above identity that
{\begin{align}
 \|\nabla_{\mathcal{A}} q\|_0^2\lesssim_0 &  (\|(u, \partial_2 u)\|_0 \| \mathcal{A}_t \|_1 +
m^2 (\|\nabla \eta\|_{L^\infty} \| \partial_2^2 \eta \|_0  + \| \partial_2^2 \eta \|_0\|\partial_2 \eta\|_{2} ))\|\nabla q\|_0.\nonumber \label{202safda005101002}
\end{align}  }
Finally, with the help of the above estimate, and \eqref{01dsaf16asdfasf}$_1$, \eqref{2020043013014} and \eqref{20200501}, we obtain
  \eqref{2017020614181721} with $i=0$.

(2) Applying $\mm{div}_{\mathcal{A}}$ to \eqref{01dsaf16asdfasf}$_2$, using then \eqref{01dsaf16asdfasf}$_1$, \eqref{202011211548}
 and the first identity in \eqref{202005101008}, we arrive at
\begin{equation}
 \Delta  q =   f,  \label{202006111729}
\end{equation}
where
$$
\begin{aligned}
f:= &  (\partial_2u_2\partial_1-\partial_1u_2\partial_2)  u_1+ (\partial_1u_1\partial_2-\partial_2u_1\partial_1)  u_2 - (\mm{div}_{\tilde{\ml{A}}}\nabla_{\ml{A}}q
+ \mm{div} \nabla_{\tilde{\ml{A}}}q) \\
&+ m^2(\partial_2^2 \mm{div}\eta+ (\partial_2\eta_2\partial_1-\partial_1\eta_2\partial_2) \partial_2^2 \eta_1+ (\partial_1\eta_1\partial_2-\partial_2\eta_1\partial_1) \partial_2^2 \eta_2) .
\end{aligned}$$

Multiplying \eqref{202006111729} by $\Delta q$ in $L^2$, we use the regularity theory of elliptic equations to get
\begin{align}
& \| \nabla q\|_1^2\lesssim_0 \|\Delta q\|_0^2 \lesssim_0\int| f  \Delta q|\mm{d} y, \label{201711118xx41x}
\end{align}
where the integral term on the right hand side can be estimated as follows.
\begin{align}
\int| f  \Delta q|\mm{d} y\lesssim_0 & ( \|\partial_2^2 u\|_0\|\nabla u\|_1 \mbox{ or } \|\partial_2 u\|_1^2+ m^2 ( \| \partial_2^3\eta\|_0 (\|\nabla \eta\|_{L^\infty} +\|\partial_2 \eta\|_2)+\| \partial_2\eta\|_{L^\infty} \|\partial_2^2\eta\|_1 ) \nonumber \\
& + (1+\|\nabla \eta\|_{L^\infty})(\| \nabla  \eta \|_{L^\infty}
   \|\nabla^2q\|_0 + \|  \partial_2  \eta \|_2
   \|\nabla q\|_{1}) )\|\Delta q \|_0\nonumber \\ \lesssim_0 &   (\|\partial_2^2 u\|_0\|\nabla u\|_1 \mbox{ or } \|\partial_2 u\|_1^2+ m^2 (
\|\partial_2^3 \eta\|_0 \|\partial_2 \eta\|_2  +\|\partial_2^2 \eta\|_1^2)   +\delta
   \|\nabla q\|_1)\|\Delta q\|_0, \label{201909121931}
\end{align}
where we have used \eqref{202safas00508}, \eqref{202safas00508x} and \eqref{201903211118} in the first inequality,
and \eqref{2020043013014} and \eqref{20206111731} in the second inequality.
Consequently, combining \eqref{201711118xx41x} with \eqref{201909121931}, one obtains \eqref{2017020614181721} with $i=2$.

(3) Applying $\nabla$ to \eqref{01dsaf16asdfasf}$_2$ and multiplying the resulting identity by $\nabla u_t$ in $L^2$, we have
$$\int |\nabla u_t|^2\mm{d}y= \int \nabla ( m^2 \partial_2^2\eta +\nu \Delta_{\ml{A}}u- \nabla_{\ml{A}}q):\nabla u_t\mm{d}y.$$
It is easy to see from the above identity that
$$\|\nabla u_t\|_0\lesssim m^2\|\partial_2^2\eta\|_1 + \|\nabla u\|_2(1+\|\nabla \eta\|_2)+\|\nabla q\|_1,$$
which, combined with \eqref{2017020614181721} with $i=1$, implies that for sufficiently small $\delta$,
 \begin{align}
 \|\nabla u_t\|_0 \lesssim_0    m^2\|\partial_2^2\eta\|_1+ \|\nabla u\|_2(1+\|  \nabla \eta\|_2+\|\nabla u\|_1) . \label{202010302129}
 \end{align}

 An application of $\partial_t$ to \eqref{202006111729} yields $\Delta  q_t = f_t$. Thus, multiplying this equation with $\Delta q_t$
 in $L^2$ and applying the regularity theory of elliptic equations, we get
 \begin{align} & \| \nabla q_t\|_1^2\lesssim_0 \|\Delta q_t\|_0^2 \lesssim_0 \int| f_t \Delta q_t|\mm{d} y.\label{20171111asfdsa8xx41x} \end{align}
 Furthermore, it is easy to verify that
 \begin{align}&\int| f_t \Delta q_t|\mm{d} y\lesssim_0   (m^2   \|\partial_2 \eta\|_{2} \| \partial_2 u\|_2+\|\partial_2 u\|_2
 \|\nabla u_t\|_0\nonumber \\ & \qquad \qquad \qquad\quad +  \|\partial_2u\|_2\|\nabla q\|_1 +\delta\|\nabla q_t\|_1)
 \|\Delta q_t\|_0.\label{2017asfdsa11118xx41x} \end{align}
 Consequently, if we making use of \eqref{2017020614181721} with $i=1$, and \eqref{202010302129}, \eqref{2017asfdsa11118xx41x}
 and Young's inequality, we obtain \eqref{2020103safda02130} from  \eqref{20171111asfdsa8xx41x}.  \hfill$\Box$
\end{pf}

\begin{lem}\label{qwe201612132242nn}
Under the conditions of \eqref{aprpiosesnewxxxx}--\eqref{aprpiose1snewxxxxz} with sufficiently small $\delta$, we have
\begin{align}
\frac{\mm{d}}{\mm{d}t}\|\nabla u(t)\|_0^2+ \nu\|u(t)\|_2^2 \lesssim & (1+m^2)\| m\partial_2^3 \eta\|_0^2, \quad t\geqslant 0
\label{qwessebdaiseqinM0846dsadfafgssgsd}
 \end{align}
and
\begin{align}
\frac{\mm{d}}{\mm{d}t} \|\nabla^2 (u, m\partial_2\eta )(t)\|_0^2 & + \nu\| u(t)\|_3^2
\lesssim F_1:=   \|\partial_2^2 \eta\|_1^2 \|m\partial_2 \eta\|_2^4+ \|\nabla \eta\|_2^2
 \|\partial_2^2 u\|_1^2\nonumber \\
&  +\|\nabla \eta\|_2\|\partial_2 u\|_2(
 m^2\|\partial_2^2 \eta\|_1\|\partial_2 \eta\|_2+\|\partial_2 u\|_1^2), \quad t\geqslant 0.
\label{qwessebdaiseqinM0846dfgssgsd}
 \end{align}
  \end{lem}
\begin{pf}
Let $\alpha$ satisfy $1\leqslant |\alpha|\leqslant 2$.
Applying $ \partial^{ \alpha} $ to \eqref{s0106pnnnn}$_1$ yields
$$ \partial^{ \alpha}(u_t+\nabla q-\nu\Delta  u-m^2\partial_2^2\eta
)=\partial^{ \alpha}(\mathcal{N}^\nu-\nabla_{\tilde{\ml{A}}}q ).$$
If we multiply the above identity by $\partial^{\alpha} u$ in $L^2$, integrate by parts and make use of \eqref{202005011525},
we get
\begin{align}
& \frac{1}{2}\frac{\mm{d}}{\mm{d}t} \|\partial^\alpha u\|_0^2 + \nu\|\nabla \partial^\alpha u \|_0^2 = I_1+I_2-
                                     \begin{cases}
                     \displaystyle  m^2 \int \partial_2^2 \eta\cdot  \partial^{2\alpha}u\mm{d}y  & \hbox{ for }|\alpha|=1; \\[0.7em]
        \displaystyle  \frac{1}{2}\frac{\mm{d}}{\mm{d}t}
\|m\partial_2  \partial^\alpha\eta\|_0^2  & \hbox{ for }|\alpha|=2,
                                          \end{cases}  \label{qwe201808100154628xxx}
\end{align}
where
\begin{align}
&I_1:= - \int   \partial^{\alpha } \mathcal{N}^\nu_{j,l}  \partial_l \partial^{\alpha} u_j\mm{d}y,\ I_2:= \int (\partial^{\alpha^-} \nabla_{\tilde{\mathcal{A}}} q  \cdot \partial^{\alpha^+} {u}-\partial^\alpha  q   \partial^\alpha \mm{div}_{\tilde{\mathcal{A}}} u ) \mm{d}y ,\\
\label{202004182208}
&\alpha^-:=
             \begin{cases}
              (\alpha_1-1,\alpha_2) & \hbox{ for }\alpha_2=0; \\
   (\alpha_1,\alpha_2-1)         & \hbox{ for }\alpha_2\geqslant 1
             \end{cases}
\mbox{ and }\alpha^+:=
             \begin{cases}
              (\alpha_1+1,\alpha_2) & \hbox{ for }\alpha_2=0; \\
   (\alpha_1,\alpha_2+1)         & \hbox{ for }\alpha_2\geqslant 1.
             \end{cases}
\end{align}
Next, we consider two cases.

(1) {\it Case $|\alpha|=1$.}

In view of \eqref{202safas00508} and \eqref{20206111731}, we find that
\begin{align}
I_1\lesssim  (1+\|\nabla \eta\|_{L^\infty})\|\ \partial_2 \eta\|_{2}\| u\|_2 \|\nabla^2 u\|_0\lesssim
\delta \|   u\|_2^2.\label{202sdfaf04101345}
\end{align}
Recalling \eqref{aprpiosesnewxxxx} and \eqref{aprpiose1snewxxxxz}, we see that
\begin{align}
 \|\partial_2 \eta\|_2 (\| \partial_2u\|_1+\|(\nabla \eta, m\partial_2 \eta)\|_2)+(\|\nabla \eta\|_2 +\|\nabla \eta\|_2^2 )/m\lesssim_0 \delta.
\label{202006141452}
\end{align}
The integral $I_2$ can be estimated as follows.
\begin{align}
I_2\lesssim&
\| \partial_2\eta\|_2 \|\nabla  q\|_0\|\nabla u\|_1 \nonumber \\
\lesssim& \| \partial_2\eta\|_2  \| u\|_2 ( \| \partial_2 u \|_1\|u\|_2+ m^2  \|\partial_2^2\eta \|_0\|\partial_2\eta\|_2)\nonumber \\
\lesssim &\delta ( \| m\partial_2^2  \eta \|_0^2 +\|  u\|_2^2),\label{20200614145dfs2}
   \end{align}
where we have used \eqref{2020043013014} and  \eqref{202safas00508} in the first inequality, \eqref{202006121323} with $i=0$ and \eqref{2017020614181721} with $i=0$ in the second inequality, and \eqref{202006141452} and  Young's inequality in the last inequality.

Finally, putting \eqref{202sdfaf04101345} and \eqref{20200614145dfs2} into \eqref{qwe201808100154628xxx},  and utilizing
\eqref{20160614fdsa19asfda57x}, \eqref{202005121605} and Young's inequality, we get \eqref{qwessebdaiseqinM0846dsadfafgssgsd}.

(2) \emph{Case $i=2$.}

We integrate by parts and use \eqref{202safas00508x} to deduce
\begin{align}
&\int \partial_1^2\tilde{\mathcal{A}}_{12}\partial_1 u\cdot \partial_1^2\partial_2u\mm{d}y
= -\int \partial_1^2\tilde{\mathcal{A}}_{12}\partial_1\partial_2 u\cdot \partial_1^2u\mm{d}y
+\int \partial_1\partial_2\tilde{\mathcal{A}}_{12}\partial_1 (\partial_1  u\cdot \partial_1^2u)\mm{d}y \nonumber \\
& \qquad \lesssim \|\nabla^3 \eta\|_0\|\nabla \partial_2^2 u\|_0 \|\partial_1^2 u\|_{1}+\|\nabla^2\partial_2\eta\|_0
 ( \|\nabla^2 u\|_1^2+ \|\partial_1 u\|_{L^\infty} \|\nabla^3  u\|_0). \label{2020410asads1345}
\end{align}
Thus, the integral $I_1$ can be estimated as follows.
\begin{align}
I_1\lesssim  &((1+\|\nabla \eta\|_{L^\infty})(\|\nabla \eta\|_{L^\infty}\|\nabla^3 u\|_0 + \|\nabla^2 \partial_2 \eta\|_0(\|\nabla^2 u\|_1+\|\partial_1 u\|_{L^\infty})+\|\nabla^3\eta\|_0( \|\partial_2  u\|_{L^\infty} \nonumber \\
&+\|\nabla \partial_2^2 u\|_0)) + (\| \partial_2\eta\|_{L^\infty}+ \|\nabla^2 \partial_2\eta\|_0)(\|\nabla\partial_2 \eta\|_{1}\|\partial_1  u\|_{L^\infty} +\|\nabla^2 \eta\|_{1}\|\nabla  u\|_{L^\infty} )) \|\nabla^2 u\|_1\nonumber \\
\lesssim &\delta \|  u\|_3^2+    \|\nabla \eta\|_2
\|\partial_2^2 u\|_1 \|u\|_3,\label{20204101345}
\end{align}
where we have used \eqref{202safas00508} and \eqref{2020410asads1345} in the first inequality,
and \eqref{202004221saffad412}, \eqref{202006101611}, \eqref{20206111731} and \eqref{202006141452} in the second inequality.

Similarly, the second integral $I_2$ can be estimated as follows, using \eqref{2017020614181721} with $i=1$.
\begin{align}
I_2\lesssim&( \| \partial_2\eta\|_2 \|u\|_3+\| \nabla \eta\|_2 \|\partial_2 u\|_2 )\|\nabla  q\|_1 \nonumber
\\ \lesssim &(\|\partial_2\eta\|_2\| u\|_3+\|\nabla \eta\|_2\|\partial_2 u\|_2)(\|\partial_2 u\|_1^2 +m^2\|\partial_2^2 \eta\|_1\|\partial_2 \eta\|_2) , \label{202004302052}
\end{align}
Consequently, plugging the above two estimates into \eqref{qwe201808100154628xxx}, and using Young's inequality, \eqref{aprpiosesnewxxxx1x}, \eqref{20160614fdsa19asfda57x} and \eqref{202006141452}, we get \eqref{qwessebdaiseqinM0846dfgssgsd}. \hfill$\Box$
\end{pf}
\begin{lem}
\label{qwe201612132242nxsfssdfsxx}
 Under the conditions of \eqref{aprpiosesnewxxxx}--\eqref{aprpiose1snewxxxxz} with sufficiently small $\delta$, we have \begin{align}
 & \frac{\mm{d}}{\mm{d}t}\left(\frac{\nu}{2}\|\nabla^{3}\eta\|_0^2  +   \sum_{|\alpha|=2} \left( \int \partial^\alpha \eta\cdot \partial^\alpha u\mm{d}y +\nu F_2^\alpha \right)\right)+  \| m \partial_2 \eta \|_2^2 \leqslant
    \| \nabla^2 u \|_{ 0}^2+  c \delta  \| u\|_3^2,
  \label{qweLem:03xx}
 \end{align}
where  $F_2^\alpha: =  \int \partial^\alpha \mathcal{A}_{12} \partial_2 \eta \cdot \partial_1\partial^\alpha\eta    \mm{d}y$.
\end{lem}
\begin{pf}
Applying $ \partial^{\alpha}$ ($|\alpha|=2$) to \eqref{s0106pnnnn}$_1$, and multiplying then the resulting identity by $\partial^{\alpha}\eta$ in $L^2$,
 one sees that
\begin{align}
 &\frac{\mm{d}}{\mm{d}t}  \left(\frac{\nu}{2}\|\nabla \partial^\alpha \eta \|_0^2+ \int \partial^\alpha \eta \cdot\partial^\alpha u\mm{d}y \right)
  +  \|m\partial_2 \partial^\alpha \eta \|_0^2 = \|\partial^\alpha u\|_0^2+I_3+I_4, \label{qwe2010154628}
\end{align}
where
$$I_3:=-\int  \partial^\alpha {\mathcal{N}}^\nu_{j,l}\cdot\partial_l\partial^\alpha\eta_j\mm{d}y,\quad
I_4:=  \int(  \partial^{\alpha^-} \nabla_{\tilde{\ml{A}}}q \cdot   \partial^{\alpha^+}  \eta+ \partial^\alpha q\partial^\alpha \mm{div} \eta)\mm{d}y,$$
and $\alpha^-$ and $\alpha^+$ are defined by \eqref{202004182208}.

Noting that
$$
\begin{aligned}
  \int    \partial^\alpha \mathcal{A}_{12} \partial_2 u\cdot \partial_1\partial^\alpha\eta    \mm{d}y  = \frac{\mm{d}}{\mm{d}t}  F_2^\alpha -    \int  \partial_t( \partial^\alpha \mathcal{A}_{12}\partial_1 \partial^\alpha\eta)\cdot \partial_2 \eta  \mm{d}y,
\end{aligned} $$
arguing similarly to the derivation of \eqref{20204101345} and \eqref{202004302052}, we find that
\begin{align}
|I_3|\leqslant &  c (\|\partial_2 \eta\|_2^2\|u\|_3+\|\nabla \eta\|_2( \|\partial_2^2 \eta\|_1 \|u\|_3(1+\|\nabla \eta \|_2)+ \|\partial_2 \eta\|_2 \|\partial_2u\|_2 )) - \frac{\mm{d}}{\mm{d}t} F_2^\alpha  , \label{202005safas0119091}\\
|I_4|\lesssim &\|\partial_2 \eta\|_2 \|\nabla \eta\|_{2} (\|\partial_2 u\|_1^2 +m^2    \|\partial_2^2 \eta\|_1  \|\partial_2 \eta\|_2 ). \label{2020050119091}
\end{align}
Inserting \eqref{202005safas0119091} and \eqref{2020050119091} into \eqref{qwe2010154628}, and using  \eqref{aprpiosesnewxxxx1x}, \eqref{20160614fdsa19asfda57x}, \eqref{202006141452} and Young's inequality, we obtain \eqref{qweLem:03xx}. \hfill$\Box$
 \end{pf}
\begin{lem} \label{202005121727}
 Under the assumptions of \eqref{aprpiosesnewxxxx}--\eqref{aprpiose1snewxxxxz}  with sufficiently small $\delta$, we have
 \begin{align}
 \frac{\mm{d}}{\mm{d}t} \|\nabla^{2-i}\partial_2^i(u, m\partial_2\eta )\|_0^2 + \nu\|\partial_2^i u \|_{{3-i}}^2
  \lesssim   \begin{cases}
 m^4\|\partial_2\eta\|_2^2 ( \|\partial_2^3 \eta\|_0^2\|\partial_2\eta\|_2^2 +\|\partial_2^2\eta\|_1^4)& \mbox{for }i=1, \\
 F_3 & \mbox{for }i=2,
            \end{cases}
\label{qwesse6dfgssgsd}
 \end{align}
 and
\begin{align}
& \frac{\mm{d}}{\mm{d}t}\left(\frac{\nu}{2}\| \nabla^{3-i}\partial_2^i\eta\|_0^2+  \sum_{|\alpha|=2-i} \left( \int
\partial^{\alpha}\partial_2^{i}\eta \cdot\partial^{\alpha}\partial_2^{i} u\mm{d}y +\nu F_4^\alpha  \right)+(i-1)F_5\right) +   \|m \partial_2^{i+1} \eta \|_{2-i}^2\nonumber \\
& \leqslant \|  \partial_2^i  u \|_{2- i}^2+\begin{cases}
c(\|  \partial_2^2 \eta\|_1\|\partial_2\eta\|_2\|\partial_2u\|_2 + \delta  \| \partial_2 u\|_2^2) & \mbox{for }i=1, \\
F_6
+\|\partial_2^2\eta\|_1\|\partial_2\eta\|_2\|\partial_2^2u\|_1  +\delta\|\partial_2 u\|_1^2 &\mbox{for }i=2,
            \end{cases}
  \label{qwe3xx}
 \end{align}
where
$$
\begin{aligned}
&F_3 := \| \partial_2^3 \eta\|_0 \|\partial_2^2u\|_0 \|u\|_2 \|\partial_2u\|_2  + \| \partial_2^3 \eta\|_0^2 (\|\partial_2u\|_2^2+
 \|m\partial_2\eta\|_2^4+m^2\|\partial_2\eta\|_2\|\partial_2u\|_2 ) \\
& \qquad\  + \|m\partial_2^2\eta\|_1^2( \|\partial_2^3\eta\|_0
\|\partial_2u\|_2+ \|m\partial_2\eta\|_2^2\|\partial_2^2\eta\|_1^2) ,\\
&F_4^\alpha:=
               \begin{cases}
         0 & \hbox{for }\alpha_1=0, \\[0.7em]
   \displaystyle     \int     \partial^\alpha \partial_2 \tilde{\mathcal{A}}_{12}
   \partial_2 \eta \cdot \partial_1\partial_2\partial^\alpha\eta  \mm{d}y & \hbox{for }\alpha_1=1 ,
               \end{cases}\\
&F_5:= \int (\partial_2^2(\mathcal{A}_{k1}\tilde{\mathcal{A}}_{k1} +\tilde{\mathcal{A}}_{11})\partial_1 \eta\cdot \partial_1\partial_2^2\eta + \partial_2^2(\mathcal{A}_{k2}\tilde{\mathcal{A}}_{k1}+\tilde{\mathcal{A}}_{12} )\partial_1 \eta\cdot \partial_2^3\eta))\mm{d}y, \\
& F_6:= m^2\|  \partial_2^2 \eta\|_1^2 ( \|\partial_2^3 \eta \|_0\| \partial_2 \eta\|_2 +\|\partial_2^2\eta\|_1^2).
\end{aligned}
$$
\end{lem}
\begin{pf}
Let $\alpha$ satisfy $|\alpha|=2-i$. Similar to the derivation of \eqref{qwe201808100154628xxx} and \eqref{qwe2010154628}, one gets
\begin{align}
& \frac{1}{2}\frac{\mm{d}}{\mm{d}t}\left(    \|\partial^\alpha\partial_2^{i}  u\|_0^2
+ \|m\partial^\alpha\partial_2^{i+1}  \eta\|_0^2 \right)  +
  \nu \|\partial^\alpha\nabla\partial_2^i  u \|_0^2 =   I_5+ I_6 \label{qwe201808100154628xxx12}
\end{align}
and
\begin{align}
 &\frac{\mm{d}}{\mm{d}t}  \left(
\frac{\nu}{2}\|\nabla \partial^\alpha \partial_2^i  \eta \|_0^2  +\int\partial^\alpha \partial_2^i  \eta \cdot \partial^\alpha \partial_2^i u\mm{d}y \right)+ \|m\partial^\alpha \partial_2^{i+1}   \eta \|_0^2
\nonumber \\
 & = \|\partial^\alpha \partial_2^i u\|_0^2+I_7+I_8, \label{qwe2010154628121321}
\end{align}
where
$$
\begin{aligned}
&  I_5:= -\int   \partial^\alpha \partial_2^i {\mathcal{N}}^\nu_{j,l}  \cdot   \partial_l \partial^\alpha\partial_2^i  u_j\mm{d}y,\ I_6:=   \int  (\partial^{\alpha} \partial_2^{i-1} \nabla_{\tilde{\mathcal{A}}} q  \cdot \partial^{\alpha} \partial_2^{i+1} {u}-\partial^\alpha \partial_2^i  q   \partial^\alpha \partial_2^i \mathrm{div}_{\tilde{\mathcal{A}}} {u} )\mm{d}y  ,\\
& I_7:=-\int \partial^\alpha \partial_2^i {\mathcal{N}}^\nu_{j,l}  \cdot \partial_l\partial^\alpha \partial_2^i\eta_j\mm{d}y,\ I_8:=\int (\partial^{\alpha } \partial_2^{i-1}\nabla_{\tilde{\mathcal{A}}} q  \cdot \partial^{\alpha } \partial_2^{i+1} \eta+\partial^{\alpha} \partial_2^{i}q   \partial^{\alpha}\partial_2^{i}\mathrm{div} \eta )\mm{d}y.
\end{aligned}$$
Next, we estimate $I_5$--$I_8$ by considering the cases $i=1$ and $i=2$ respectively.

(1) \emph{Case $i=1$}

From \eqref{aprpiosesnewxxxx1x}, \eqref{202006101611} and \eqref{202safas00508x} we get
\begin{align}
&I_5 \lesssim   \|\partial_2\eta\|_2\| \partial_2 u\|_2^2\lesssim \delta  \| \partial_2 u\|_2^2. \label{202006151927}
\end{align}
The integral $I_6$ can be estimated in the following way.
\begin{align}
& I_6 \lesssim \|  \partial_2\eta\|_2 \|\partial_2 u\|_2\|\nabla q\|_1\lesssim
  m^2 \|  \partial_2\eta\|_2(   \|\partial_2^3 \eta\|_0 \|\partial_2 \eta\|_2 +\|\partial_2^2 \eta\|_1^2 ) \|\partial_2 u\|_2+ \delta\|\partial_2 u\|_2^2,\label{20200615saf1927}
\end{align}
where we have used \eqref{2020043013014} and \eqref{202006101611}--\eqref{202safas00508x} in the first inequality,
and \eqref{2017020614181721} with $i=1$, \eqref{202006141452} and  Young's inequality in the second inequality.

Noting that for $\alpha_1=1$,
$$
\begin{aligned}
 \int    \partial^\alpha \partial_2 {\mathcal{A}}_{12} \partial_2 u\cdot \partial_1\partial_2\partial^\alpha\eta  \mm{d}y =  \frac{\mm{d}}{\mm{d}t} F_4^\alpha -   \int   \partial_t(   \partial^\alpha\partial_2 \mathcal{A}_{12}\partial_1\partial_2 \partial^\alpha\eta)\cdot \partial_2 \eta  \mm{d}y,
\end{aligned}
$$
we argue, similarly to \eqref{202006151927} and \eqref{20200615saf1927}, to infer that
\begin{align}
&I_7 \leqslant c\|  \partial_2^2 \eta\|_1\|\partial_2\eta\|_2 \|  \partial_2 u\|_2-  \frac{\mm{d}}{\mm{d}t} F_4^\alpha ,  \label{20200615saf19s27} \\
& I_8\lesssim\|\partial_2^2\eta\|_1 \|  \partial_2\eta\|_2\|\nabla q\|_1 \lesssim\delta ( \|m \partial_2^2 \eta\|_1^2+ \|\partial_2 u\|_2^2). \label{20200615saf1927xxxxx}
\end{align}

Thanks to the  four estimates \eqref{202006151927}--\eqref{20200615saf1927xxxxx} and Young's inequality, we derive
\eqref{qwesse6dfgssgsd} and \eqref{qwe3xx} from  \eqref{qwe201808100154628xxx12} and \eqref{qwe2010154628121321} with $i=1$, respectively.

(2) \emph{Case $i=2$.}

By employing a partial integration, \eqref{202006101611} and \eqref{202safas00508x}, one can easily see that
\begin{align}
\int\partial_2^2\tilde{\mathcal{A}}_{12}\partial_1 u\cdot \partial_2^3u\mm{d}y
& = \int\partial_2^2\partial_1\eta_2  \partial_1 \partial_2 u\cdot \partial_2^2u\mm{d}y -\int\partial_2^3\eta_2 \partial_1( \partial_1 u\cdot \partial_2^2u)\mm{d}y  \nonumber  \\
& \lesssim  \|\partial_2^3\eta\|_0\|\partial_2^2u\|_1\|\partial_2 u\|_2 +\|\partial_2\eta\|_2\|   \partial_2^2 u\|_{1}^2 \nonumber .
\end{align}
Similarly, we can also obtain
\begin{align}
 \int\partial_2^2{\mathcal{A}}_{k2}\tilde{\mathcal{A}}_{k1}\partial_1 u\cdot \partial_2^3u\mm{d}y \lesssim  \|\partial_2^3\eta\|_0\|\partial_2^2u\|_1\|\partial_2 u\|_2 +\|\partial_2\eta\|_2\|   \partial_2^2 u\|_{1}^2 .\nonumber
\end{align}
Thanks to the above two estimates, we easily deduce that
\begin{align}
I_5 \lesssim &\|\partial_2^3\eta\|_0\|\partial_2^2u\|_1\|\partial_2 u\|_2  + \|\partial_2\eta\|_2\|   \partial_2^2 u\|_{1}^2 + \int\partial_2^2\tilde{\mathcal{A}}_{12}\partial_1 u\cdot \partial_2^3u\mm{d}y+ \int\partial_2^2\tilde{\mathcal{A}}_{k2}\tilde{\mathcal{A}}_{k1}\partial_1 u\cdot \partial_2^3u\mm{d}y\nonumber\\
\lesssim &\|\partial_2^3\eta\|_0\|\partial_2^2u\|_1\|\partial_2 u\|_2+ \|\partial_2\eta\|_2\|   \partial_2^2 u\|_{1}^2 .\nonumber
\end{align}

Noting that
\begin{align}
&\int (\partial_2^2(\mathcal{A}_{k1}\tilde{\mathcal{A}}_{k1} +\tilde{\mathcal{A}}_{11})\partial_1 u\cdot \partial_1\partial_2^2\eta + \partial_2^2(\mathcal{A}_{k2}\tilde{\mathcal{A}}_{k1}+\tilde{\mathcal{A}}_{12} )\partial_1 u\cdot \partial_2^3\eta))\mm{d}y\nonumber \\
&\qquad \leqslant \frac{\mm{d}}{\mm{d}t}F_5+c \|\partial_2^2\eta\|_1\|\partial_2\eta\|_2\|\partial_2^2u\|_1,
\end{align}
 making use of \eqref{2020043013014} and \eqref{202006101611}--\eqref{202safas00508x}, we can control $I_6$--$I_8$ as follows.
\begin{align}
  I_6   \lesssim   &  (\| \partial_2^3 \eta\|_0 \|\partial_2  u \|_2 +\| \partial_2  \eta\|_2 \|\partial_2^2  u \|_1  )(   \|\partial_2^2u\|_0 \|u\|_2  +m^2 ( \|\partial_2^3 \eta\|_0 \|\partial_2 \eta\|_2  +  \|\partial_2^2 \eta\|_1^2)  ),\nonumber \\
I_7 \leqslant &  - \frac{\mm{d}}{\mm{d}t}F_5+c
\|\partial_2^2\eta\|_1\|\partial_2\eta\|_2\|\partial_2^2u\|_1 , \nonumber \\
 I_8   \lesssim & ( \|\partial_2^3 \eta \|_0 \| \partial_2 \eta\|_2+\|\partial_2^2\eta\|_1^2)(  \|\partial_2 u\|_1^2+ m^2(  \|\partial_2^3 \eta\|_0 \|\partial_2 \eta\|_2 + \|  \partial_2^2 \eta\|_1^2 ) ) .\nonumber \end{align}

Utilizing the estimates for $I_5$--$I_8$, \eqref{aprpiosesnewxxxx1x}, and \eqref{201005041609}, \eqref{202006141452} and Young's inequality,
we get \eqref{qwesse6dfgssgsd} and \eqref{qwe3xx} from \eqref{qwe201808100154628xxx12} and \eqref{qwe2010154628121321} with $i=2$,
respectively.  \hfill$\Box$
\end{pf}

\subsection{Stability estimates} \label{202006141016}
With the energy estimates in Lemmas \ref{qwe201612132242nn}--\ref{202005121727} in hand, we are in a position to establish
the stability estimate \eqref{20190safd5041053}.

By virtue of \eqref{20206111731}, it is easy to see that
\begin{align}
&|F_2^\alpha|\leqslant \delta \|\nabla \eta\|_2^2 , \label{200206161733} \\
&|F_4^\alpha|\leqslant \delta \|\partial_2 \eta\|_2^2,  \label{2020061737}\\
&|F_5^\alpha|\leqslant \delta \|\partial_2^2\eta\|_1^2.  \label{202006sfa1737}
\end{align}
Thus, we can use \eqref{20160614fdsa19asfda57x}, \eqref{202006141452}, \eqref{200206161733} and Young's inequality to
derive from  \eqref{qwessebdaiseqinM0846dfgssgsd}  and \eqref{qweLem:03xx} the two-order energy inequality:
\begin{align}
&\frac{\mm{d}}{\mm{d}t}    {\mathcal{E}}+ c(\|   u\|_3^2+\| m \partial_2 \eta \|_2^2)
 \lesssim \|( \nabla \eta, u, m\partial_2\eta)\|_2^2(\| \partial_2  u  \|_2^2  +\|m\partial_2^2\eta\|_1^2)
\label{2019112safd62012}
\end{align}
for sufficiently small $\delta$, where
$$
\begin{aligned}
 {\mathcal{E}} :=&  c  \|\nabla^2 (u,m\partial_2  \eta)\|_0^2 +\frac{\nu}{2}\|\nabla^3 \eta\|_0^2+
  \sum_{  |\alpha| = 2} \left(\int \partial^\alpha \eta \cdot \partial^\alpha   u\mm{d}y+\nu F_2^\alpha \right),
\end{aligned}$$
satisfying
\begin{align}
\mathfrak{E}_{2,0} \lesssim   {\mathcal{E}} \lesssim \mathfrak{E}_{2,0}.
\label{20205091359}
\end{align}

Furthermore, if one utilizes \eqref{20160614fdsa19asfda57x}, \eqref{202006141452}, \eqref{2020061737} and Young's inequality, one gets from  \eqref{qwesse6dfgssgsd} and \eqref{qwe3xx} that
\begin{align}
&\frac{\mm{d}}{\mm{d}t}    {\mathcal{E}}_i+
c (\| \partial_2^i u\|_{3-i}^2+\|m \partial_2^{i+1}   \eta  \|_{2-i}^2)   \lesssim
                                        \begin{cases}
    0             & \hbox{ for }i=1;\\
 {F}_7
  & \hbox{ for }i=2,
\end{cases}
\label{20191sadf1262012}
\end{align}
where  $F_7:= c F_3+F_6+\|\partial_2^2\eta\|_1^2\|\partial_2\eta\|_2^2$ and
$$
\begin{aligned}
{\mathcal{E}}_i  := &c \|\nabla^{2-i}\partial_2^i (u,m\partial_2  \eta)\|_0^2 +\frac{\nu}{2}\|\nabla^{3-i}\partial_2^{i}  \eta\|_0^2\\
&+
\sum_{  |\alpha| = 2-i} \left(\int  \partial^\alpha \partial^i_2\eta \cdot \partial^\alpha\partial^i_2 u\mm{d}y
+\nu F^\alpha_4\right)+(i-1)F_5,
\end{aligned}$$
satisfying
\begin{align}
\mathfrak{E}_{2,i}\lesssim   {\mathcal{E}}_i \lesssim \mathfrak{E}_{2,i}.  \label{202006141102}
\end{align}

From \eqref{20191sadf1262012} with $i=1$ we find that
\begin{align}
& \mathfrak{E}_{2,1}(t) +
c  \int_0^t(\|   \partial_2 u\|_2^2+\|m \partial_2^2 \eta \|_1^2) \mm{d}\tau\lesssim \mathfrak{E}_{2,1}^0.
\label{2019112safdsafsaf6201sdaf2}
\end{align}
Thanks to \eqref{20205091359}, an application of Gronwall's lemma to \eqref{2019112safd62012} yields
$$\mathfrak{E}_{2,0}  \lesssim  \mathfrak{E}_{2,0}^0 e^{c \int_0^T (\|  \partial_2u\|_2^2+\|m\partial_2^2\eta\|_1^2)\mm{d}\tau}.$$

Combining \eqref{20205091359} and \eqref{2019112safdsafsaf6201sdaf2} with the above estimate, we see that there is $\delta_1\in (0,1]$,
 such that for any $\delta\leqslant \delta_1$,
\begin{align}
 \mathfrak{E}_{2,0}(t) +\int_0^t( \| \partial_2 u  \|_2^2+ \|m \partial_2^2 \eta  \|_1^2) \mm{d}\tau
 \leqslant  {c}_1 \mathfrak{E}_{2,0}^0 e^{ {c}_2 \mathfrak{E}_{2,1}^0  }/4\;\;\;\mbox{ for any }t\in I_T,\label{20190sfda5sadfa041053}
\end{align}
where $c_1\geqslant 4$. In addition, due to \eqref{2019112safdsafsaf6201sdaf2} and \eqref{20190sfda5sadfa041053}, we can obtain from \eqref{2019112safd62012} that
\begin{equation}
\int_0^t (\|  u\|_3^2+\|m \partial_2 \eta  \|_2^2)\mm{d}t \lesssim  \mathfrak{E}_{2,0}^0(1+\mathfrak{E}_{2,1}^0 e^{ {c}_2 \mathfrak{E}_{2,1}^0}). \label{202008011727}
\end{equation}

Finally, if we take
\begin{equation}
\label{201911262060}
K:=  \sqrt{ c_1 \mathfrak{E}_{2,0}^0 e^{c_2 \mathfrak{E}_{2,1}^0}}>0.
\end{equation}
 we complete the derivation of \eqref{aprpiose1} under the conditions \eqref{aprpiosesnewxxxx} and \eqref{aprpiose1snewxxxxz}
 with $\delta\leqslant \delta_2$.

\subsection{Proof of Theorem \ref{201904301948}} \label{20206141015}
We start with introducing a local (-in-time) well-posedness result for the initial value problem \eqref{01dsaf16asdfasf}--\eqref{01dsaf16asdfasfsaf} and a result concerning diffeomorphism mappings.
\begin{pro} \label{pro:0401nxdxx}
Let $(\eta^0,u^0)\in H^3\times H^2$ satisfy $\|(\nabla\eta^0,u^0)\|_2\leqslant B$ and
$\mm{div}_{\mathcal{A}^0}u^0=0$, where $B$ is a positive constant, $\zeta^0:= \eta^0+y$ and $\mathcal{A}^0$ is defined by $\zeta^0$.
Then there is a constant $\delta_2\in (0,1]$, such that for any $(\eta^0,u^0)$ satisfying
\begin{align}
\|\nabla \eta^0\|_{2,2} \leqslant \delta_2, \label{201912261426}
\end{align}
 there exist a local existence time $T>0$ (depending possibly on $B$, $\nu$, $m$ and $\delta_2$) and a unique local strong solution
$(\eta, u,q)\in C^0(\overline{I_T},{H}^{3} )\times \mathcal{U}_{1,T}\times C^0(\overline{I_T},\underline{H}^2)$
to the initial value problem \eqref{01dsaf16asdfasf}--\eqref{01dsaf16asdfasfsaf}, satisfying
$0<\inf_{(y,t)\in \mathbb{R}^2\times \overline{I_T}} \det(\nabla \eta+I)$ and $\sup_{t\in \overline{I_T}}\|\nabla \eta\|_{2,2}\leqslant 2\delta_2$.
\footnote{Here the uniqueness means that if there is another solution
$( \tilde{\eta}, \tilde{u},\tilde{q})\in C^0(\overline{I_T},{H}^{3} ) \times \mathcal{U}_{1,T}\times C^0(\overline{I_T},\underline{H}^2)$
 satisfying $0<\inf_{(y,t)\in \mathbb{R}^2\times \overline{I_T}} \det(\nabla \tilde{\eta}+I)$, then
 $(\tilde{\eta},\tilde{u},\tilde{q})=(\eta,u,q)$ by virtue of the smallness condition
 ``$\sup_{t\in \overline{I_T}}\|\nabla \eta\|_{2,2}\leqslant 2\delta_2$''.}
\end{pro}
\begin{pf}
The proof of Proposition \ref{pro:0401nxdxx} will be given in Section \ref{201912062021}.
\hfill $\Box$
\end{pf}
\begin{rem}
\label{202010252143}

If $(\eta^0,u^0)$ in Proposition \ref{pro:0401nxdxx} further satisfies $(\eta^0,u^0)\in \underline{H}^3_1\times \underline{H}^2$
and the odevity conditions \eqref{202007301500}, then
$(\eta,u)$ belongs to $ C^0(\overline{I_T},\underline{H}^3_1)\times \underline{\mathcal{U}}_T $ and satisfies the odevity conditions
 \eqref{20208192055} as $(\eta^0,u^0)$ does.
\end{rem}
\begin{pro}\label{pro:0401nasfxdxx}
There is a positive constant $\delta_3 $, such that for any $\varphi\in H^3$ satisfying
$\|\nabla \varphi\|_{2,2}\leqslant \delta_3$, we have (after possibly being redefined on a set of measure zero) $\det(\nabla \varphi+I)>1/2$ and
\begin{align}
&   \psi : \mathbb{R}^2\to \mathbb{R}^2 \mbox{ is a }C^1\mbox{ homeomorphism mapping},
\end{align}
where $\psi:=\varphi+y$.
\end{pro}
\begin{pf}
In view of \eqref{202004221saffad412}, we see that $\det(\nabla\varphi +I)>1/2$ and
$\|\nabla \varphi\|_{L^\infty}\lesssim_0\|\nabla \varphi\|_{2,2}\leqslant  \delta_3$.
 Thus, we easily verify that for sufficiently small $\delta_3$, $\psi : \mathbb{R}^2\to\mathbb{R}^2$ is a $C^1$ homeomorphism mapping,
 {please refer to \cite[Lemma 4.2]{JFJSOMITIN} for a detailed proof.}
\hfill $\Box$
\end{pf}

With the \emph{a priori} estimate \eqref{20190sfda5sadfa041053} (under the assumptions \eqref{aprpiosesnewxxxx} and \eqref{aprpiose1snewxxxxz}
with $\delta\leqslant \delta_2$) and Propositions \ref{pro:0401nxdxx} and \ref{pro:0401nasfxdxx} in hand, we can easily establish
Theorem \ref{201904301948}. We briefly give the proof below.

Let $m$ and $(\eta^0,u^0)\in (\underline{H}^3_{1}\cap H^3_*) \times \underline{H}^2$ satisfy the odevity conditions \eqref{202007301500}
and
\begin{align}
\max\{K^{1/2},K^2\}/m\leqslant\min\{\delta_1,  \delta_2/c_{0}, \delta_3/c_{0}\}=:c_3\leqslant 1,
\label{201912261425}
\end{align}
where $K$ is defined by \eqref{201911262060}, and  {the constant $c_0$ is the same as in \eqref{20202111191231}. }
Thus we see that $\eta^0$ satisfies \eqref{201912261426} by \eqref{20202111191231} and \eqref{201912261425}. Thus,
by virtue of Proposition \ref{pro:0401nxdxx} and Remark \ref{202010252143}, there exists a unique local solution $(\eta, u,q)$
of \eqref{01dsaf16asdfasf}--\eqref{01dsaf16asdfasfsaf} with a maximal existence time $T^{\max}$, satisfying
\begin{itemize}
  \item for any $T\in  I_{T^{\max}}$,
the solution $(\eta,u,q)$ belongs to $C^0(\overline{I_T},\underline{H}^3_1)\times \underline{\mathcal{U}}_T\times C^0( \overline{I_T},\underline{H}^2) $ and $$\sup_{t\in \overline{I_T}} \|\nabla \eta\|_{2,2}\leqslant 2\delta_2;$$
  \item $\limsup_{t\to T^{\max} }\|\nabla \eta( t)\|_{2,2} > \delta_2$ or $\limsup_{t\to T^{\max} }\|u( t)\|_2=\infty$, if $T^{\max}<\infty$.
\end{itemize}
In addition, the solution enjoys the odevity conditions \eqref{20208192055}.

Let
\begin{equation}
\nonumber
T^{*}=\sup\left\{ T \in I_{T^{\max}}~\left|~\mathfrak{E}_{2,0}(t) \leqslant K^2 \mbox{ for any }t\leqslant T\right.\right\}.
\end{equation}
Recalling the definition of $K$ and  the condition $c_1\geqslant 4$, we easily see that the definition of $T^*$ makes sense and $T^*>0$.

By \eqref{20202111191231} we have $\| \nabla \eta\|_{2,2} \leqslant \delta_3$ for all $t\in I_{T^*}$,
then $\eta(t)\in H^3_*$ for all $t\in I_{T^*}$ by Proposition \ref{pro:0401nasfxdxx}.
Thus, to obtain the existence of a global solution, it suffices to verify $T^*=\infty$. Now, we show this by contradiction.

Assume $T^*<\infty$. Keeping in mind that $T^{\max}$ denotes the maximal existence time and $K/m\leqslant \delta_2/c_{0}$ by virtue of
\eqref{201912261425}, we apply Proposition \ref{pro:0401nxdxx} to find that $T^{\max}>T^*$  and
\begin{equation}
\label{201911262202}
  \mathfrak{E}_{2,0}(T^*)    =K^2.
\end{equation}
Since $\max\{K,K^2\}/m\leqslant  \delta_1$ and $\sup_{0\leqslant t\leqslant T^*} \mathfrak{E}_{2,0}(t) \leqslant K^2,$
 we can still show that the solution $(\eta,u )$ enjoys the stability estimate \eqref{20190sfda5sadfa041053} with $T^*$ in place of $T$
 by the regularity of $(\eta,u,q)$. More precisely, one has $\sup_{0\leqslant t\leqslant T^*} \mathfrak{E}_{2,0}(t) \leqslant K^2/4,$
which contradicts with \eqref{201911262202}. Hence, $ T^*=\infty$, and thus $T^{\max}=\infty$.

Obviously, the global solution $(\eta,u)$ enjoys the stability estimate \eqref{20190safd5041053} by using \eqref{20190sfda5sadfa041053}
and \eqref{202008011727}, and the estimates \eqref{20211191238}--\eqref{202011119211101} are easily obtained from
 \eqref{20202111191231}--\eqref{202006121323} and \eqref{2017020614181721} with $i=1$.
 The uniqueness of the global solutions is obvious due to the uniqueness of the local solutions in Proposition \ref{pro:0401nxdxx}
 and the fact $\sup_{t\geqslant 0}\|\nabla \eta\|_{2,2}\leqslant 2 \delta_2$. This completes the proof of Theorem \ref{201904301948}.

\section{Proof of Theorem \ref{202006211824saf}}\label{2020dsafa11sdf92326}
 {We now proceed to the derivation of the time-decay estimates stated in Theorem \ref{202006211824saf}.
  It should be noted that, under the assumptions of Theorem \ref{202006211824saf}, the solution $(\eta,u,q)$ of \eqref{01dsaf16asdfasf}--\eqref{01dsaf16asdfasfsaf}, established in
Theorem \ref{201904301948}, satisfies} the differential inequalities in Lemma \ref{qwe201612132242nn}, \eqref{qwesse6dfgssgsd} and \eqref{20191sadf1262012} for a.e. $t>0$, and the estimates \eqref{202006141102} for any $t>0$.

\subsection{Decay estimates for $\partial_2\eta$}
We begin with the derivation of the decay estimates for $\eta$. Multiplying \eqref{qwessebdaiseqinM0846dfgssgsd}, \eqref{qwesse6dfgssgsd}
and \eqref{20191sadf1262012} by $\langle t\rangle$, $\langle t\rangle^{i+1}$ and $\langle t\rangle^i$ respectively, we obtain for a.e. $t>0$,
\begin{align}
&\frac{\mm{d}}{\mm{d}t} \left(\langle t\rangle \| \nabla^2(u,m\partial_2\eta )\|_0^2 \right)+\nu\langle t\rangle \|u\|_3^2\lesssim
\| \nabla^2(u, m\partial_2\eta)\|_0^2+\langle t\rangle  F_1
\label{2019112ssfaasfdfa62012}\\
&\frac{\mm{d}}{\mm{d}t} \left(\langle t\rangle^{i+1} \|\nabla^{2-i} \partial_2^i(u,m\partial_2\eta )\|_0^2 \right)
+\nu\langle t\rangle^{i+1}\| \partial_2^iu\|_{3-i}^2\nonumber \\
&\quad \lesssim\langle t\rangle^i \| \nabla^{2-i}\partial_2^i(u,m\partial_2\eta)\|_0^2   + \begin{cases}
 t^2m^4\|\partial_2 \eta\|_2^2 (\|\partial_2^3 \eta\|_0^2 \|\partial_2 \eta\|_2^2 +\|\partial_2^2 \eta\|_1^4 )& \hbox{ for }i=1, \\
\langle t\rangle^3  F_3
&\hbox{ for }i=2,
            \end{cases}
\label{2019112sdfa62012}
\end{align}
and
\begin{align}
\frac{\mm{d}}{\mm{d}t}   (\langle t\rangle^i {\mathcal{E}}_i)+
c  \langle t\rangle^i( \|\partial_2^i    u\|_{3-i}^2+\|m \partial_2^{i+1} \eta  \|_{2-i}^2)
   \lesssim \langle t\rangle^{i-1} {\mathcal{E}}_i  +
\begin{cases}
0       & \hbox{ for }i=1, \\
\langle t\rangle^2{F}_7& \hbox{ for }i=2.
\end{cases}
\label{201911262012}
\end{align}

Integrating \eqref{201911262012} with $i=1$ over $(0,t)$, and then using \eqref{201909281832},  \eqref{20190safd5041053}
and \eqref{202006141102}, we have
\begin{align}
& \langle t\rangle \mathfrak{E}_{2,1} +  \int_0^t\langle \tau\rangle  (\|\partial_2 u\|_2^2+\|m \partial_2^2  \eta \|_{1}^2 ) \mm{d}\tau
\nonumber \\
& \qquad \leqslant  \left(\sqrt{\mathfrak{E}_{2,0}^0 e^{c_2 \mathfrak{E}_{2,1}^0} } +\mathfrak{E}_{2,0}^0 e^{c_2 \mathfrak{E}_{2,1}^0}\right) (1+\mathfrak{E}_{2,1}^0)\;\; \mbox{ for any }t> 0.
\label{201911262fsadf012}
\end{align}

Denoting
$$ \begin{aligned}
\tilde{F}_7= &
  \|  u\|_3^2+\|\partial_2\eta\|_2^2( 1+ m^2+m^4\|\partial_2\eta\|_2^2),
\end{aligned}$$
one gets by \eqref{201909281832} and \eqref{20190safd5041053} that
\begin{align}
 \int_0^t \tilde{F}_7\mm{d}\tau\leqslant C.\nonumber
\end{align}
Applying Gronwall's lemma to \eqref{201911262012} with $i=2$, and utilizing \eqref{202006141102}, \eqref{201911262fsadf012} and
the above estimate, and the fact $F_7\lesssim (\|\partial_2^2\eta\|_1^2+ \|\partial_2^2 u\|_0^2+\|m\partial_2^3 \eta\|_0^2 )\tilde{F}_7$,
one can infer from \eqref{201911262012} with $i=2$ that
\begin{align}
\langle t\rangle^2\mathfrak{E}_{2,2}
+ \int_0^t\langle \tau\rangle^2 (\| \partial_2^2  u\|_1^2+\|m \partial_2^3 \eta \|_0^2  ) \mm{d}\tau
&\lesssim \mathfrak{E}_{2,2}^0+
\int_0^t\langle \tau\rangle\mathfrak{E}_{2,2} \mm{d}\tau  e^{c \int_0^t \tilde{F}_7\mm{d}\tau } \nonumber \\
& \leqslant C \quad\mbox{ for all }t> 0. \label{202061309fsad40}
\end{align}

Integrating \eqref{2019112ssfaasfdfa62012} and \eqref{2019112sdfa62012} over $(0,t)$,
using \eqref{201909281832}, \eqref{20190safd5041053},  \eqref{201911262fsadf012} and \eqref{202061309fsad40}, we conclude
\begin{align}
&  \langle t\rangle^{i+1}  \|  \partial_2^i( u,m\partial_2 \eta ) \|_{2-i}^2+ \int_0^t  \langle \tau\rangle^{i+1}\| \partial_2^i   u \|_{3-i}^2\mm{d}\tau
\leqslant  C\mbox{ for }0\leqslant i\leqslant 2.  \label{20202061309740}
\end{align}

\subsection{Decay estimates for $u$} \label{202008231326}
  We proceed to derive decay estimates of higher derivatives of the velocity. Multiplying \eqref{qwessebdaiseqinM0846dfgssgsd} by $e^{\nu t/2}$
  and using Young's inequality, we obtain
  \begin{align}
\frac{\mm{d}}{\mm{d}t}( e^{\nu t/2} \|\nabla^2  u\|_0^2) + c\int_0^t e^{\nu \tau/2} \|   u \|_3^2 \mm{d}\tau \lesssim
 e^{\nu \tau/2}(m^4\|  \partial_2^2  \eta \|_1^2 +F_1),\nonumber
 \end{align}
which implies
  \begin{align}
  \|u\|_2^2  +  \int_0^t e^{\nu (\tau-t)/2} \|   u \|_3^2 \mm{d}\tau \lesssim
 \|\nabla^2 u^0\|_0^2 e^{-\nu t/2} +\int_0^t e^{\nu (\tau-t)/2}(m^4\|\partial_2^2\eta\|_1^2 +F_1)\mm{d}\tau.
\label{qwessebdafgsfassgsd}
 \end{align}
In addition,
\begin{itemize}
  \item  for any given $a$, $\theta\geqslant 0$,
there is a positive constant ${c}$, depending only on $a$ and $\theta$, such that
\begin{align}
e^{-at}\leqslant  {c}  \langle  t\rangle^{-\theta}\mbox{ for any }t\geqslant 0.
\label{20208021526}
\end{align}
  \item for
any given $r_1>0$ and $r_2\in [0,r_1]$,
there is a positive constant ${c}$, depending only on $r_1$ and $r_2$, such that (see Lemma 2.5 in \cite{ODM172})
\begin{align}
\int_0^t \langle  t-\tau\rangle^{-r_1}\langle  \tau\rangle^{-r_2}\mm{d}\tau\leqslant {c}\langle t\rangle^{-r_2}. \label{2020081530}
\end{align}
  \item for any $\tau\in (0,t)$,
\begin{align}
\langle t\rangle\leqslant 4 \langle  t-\tau\rangle \langle  \tau\rangle.  \label{202008021527}
\end{align}
\end{itemize}
Therefore, we make use of \eqref{20190safd5041053}, \eqref{20202061309740} with $i=1$, \eqref{20208021526}--\eqref{202008021527} and H\"older's inequality to deduce from \eqref{qwessebdafgsfassgsd} that
\begin{align}
  \|  u\|_2^2 +  \int_0^t  e^{\nu (\tau-t)/2} \|   u \|_3^2 \mm{d}\tau
 \leqslant (1+m^2)C  \langle t\rangle^{-2} \label{202010261956}.
 \end{align}

Similarly to \eqref{202010261956}, one can utilize \eqref{20202061309740} with $i=3$, \eqref{20208021526} and \eqref{2020081530} to get from \eqref{qwessebdaiseqinM0846dsadfafgssgsd} that
\begin{align}
\|  u\|_1^2 + \int_0^t e^{\nu (\tau-t)/2} \|   u \|_2^2 \mm{d}\tau \leqslant (1+m^2)C \langle  t\rangle^{-3}.\nonumber
\end{align}
Finally, from \eqref{201911262fsadf012}--\eqref{20202061309740} the estimate \eqref{202005201033} follows,
while from the above two estimates we obtain \eqref{20200safa5201033x}.

\section{Proof of Theorem \ref{201912041028}}\label{2020011sdf92326}

This section is devoted to the proof of Theorem \ref{201912041028}.
 Let $(\eta^0,u^0)$ satisfy all the assumptions in Theorem \ref{201904301948} and $(\eta,u,q)$ be the solution constructed by Theorem \ref{201904301948}.
 Then, by the regularity theory of the Stokes problem, there exists a unique solution $(\eta^{\mm{r}},u^{\mm{r}},Q_1,Q_2)$ satisfying
\begin{align}
   \begin{cases}
   -\Delta \eta^{\mm{r}} +\nabla Q_1=0, \\
   \mm{div} \eta^{\mm{r}}=-\mm{div}\eta^0,\\
     (\eta^{\mm{r}})_{\mathbb{T}^2}  = 0
   \end{cases}
\label{203212051247}
   \end{align}
 and
\begin{align}\begin{cases}
-\Delta u^{\mm{r}} +\nabla Q_2=0, \\
\mm{div} u^{\mm{r}}=\mm{div}_{\tilde{\mathcal{A}}^0}u^0,\\
({u}^{\mm{r}})_{\mathbb{T}^2}  = 0, \label{2032120512471}
\end{cases}
\end{align}
where $\tilde{\mathcal{A}}^0:={\mathcal{A}}^0-I$. Moreover, $(\eta^{\mm{r}},u^{\mm{r}})$ satisfies  \eqref{202011032123}, also cf.
the derivation of \eqref{2017020614181721}.

Let $\tilde{\eta}^0=\eta^0+\eta^{\mm{r}}$ and $\tilde{u}^0=u^0+u^{\mm{r}}$. Thus, it is easy to see that $(\tilde{\eta}^0,\tilde{u}^0)$
belongs to $ \underline{H}^3_\sigma\times \underline{H}^2_\sigma$ and enjoys the odevity conditions as $(\eta^0,u^0)$ does.
Therefore, there exists a unique global solution $(\eta^{\mm{L}},u^{\mm{L}},q^{\mm{L}} )\in C^0(\mathbb{R}_0^+, \underline{H}^{3}_\sigma)\times \underline{\mathcal{U}}_\infty\times C^0(\mathbb{R}_0^+, \underline{H}^2)$
to the initial-value problem \eqref{202001070914}. Moreover, the solution enjoys the odevity conditions as $(\eta,u)$ does.

Similarly to \eqref{202005201033}, we employ \eqref{201909281832} and \eqref{202011032123} to see that the solution
$(\eta^{\mm{L}}, u^{\mm{L}})$ of the linearized problem enjoys the following estimate:
\begin{align}
&  \|  \nabla \eta^{\mm{L}}\|_2^2+ \sum_{i=1}^2\left(\langle t\rangle^{i}\|\partial_2^i\eta^{\mm{L}}\|_{3-i}^2
 +\int_0^t \langle \tau\rangle^i \| m  \partial_2^{i+1} \eta^{\mm{L}}   \|_{2-i}^2   \mm{d}\tau\right)\nonumber\\
&\quad + \sum_{i=0}^2\left(\langle t\rangle^{i+1}  \|  \partial_2^{i}( u^{\mm{L}},m\partial_2 \eta^{\mm{L}} ) \|_{2-i}^2+ \int_0^t  \langle \tau\rangle^{i+1}\| \partial_2^{i}   u^{\mm{L}} \|_{3-i}^2\mm{d}\tau \right)\lesssim  \mathfrak{E}_{2,0}^0 .  \label{2020052010sdfa33}
\end{align}

Let $(\eta^{\mm{d}}, u^{\mm{d}})= (\eta-\eta^{\mm{L}}, u-u^{\mm{L}})$, then the error function $(\eta^{\mm{d}}, u^{\mm{d}} )$ satisfies \eqref{01dsaf16safafasdfasfxx}.
It is easy to see from \eqref{01dsaf16safafasdfasfxx} that $(u^{\mm{d}})_{\mathbb{T}^2}=(\eta^{\mm{d}})_{\mathbb{T}^2}=0$ for any $t> 0$, since $(\eta^{\mm{r}})_{\mathbb{T}^2}  = ({u}^{\mm{r}})_{\mathbb{T}^2}=0$. Moreover, $(\eta^{\mm{d}},u^{\mm{d}})$  also enjoys the odevity conditions
as  $(\eta,u)$ does.

Recalling that $(\eta,u,q)$ is constructed by Theorem \ref{201904301948}, the solution $(\eta,u)$ satisfies all the estimates
in Lemma \ref{201805141072dsafa}. Hence, we can follow the arguments in the proof of Lemmas \ref{qwe201612132242nn}--\ref{202005121727}
with slight modifications to derive from \eqref{01dsaf16safafasdfasfxx} that
\begin{align}
& \frac{\mm{d}}{\mm{d}t}  \| \nabla  (u^{\mm{d}},m \partial_2 \eta^{\mm{d}})\|_0^2 +
 \nu\|  u^{\mm{d}} \|_2^2  \lesssim F^{\mm{d}}_1 , \label{201702061safasaf4181721asfdx} \\
&  \frac{\mm{d}}{\mm{d}t}\bigg( \frac{\nu}{2}\|\nabla^3 \eta^{\mm{d}}\|_0^2+ \sum_{|\alpha|=2} \int \partial^\alpha \eta^{\mm{d}} \cdot \partial^\alpha u^{\mm{d}}\mm{d}y +
\int  \partial_1^2\tilde{\mathcal{A}}_{k2}\partial_2 \eta\cdot \partial_{1}^3\eta^{\mm{d}} \mm{d}y \bigg)
\nonumber \\
 &\quad +\|m \partial_2\eta^{\mm{d}} \|_2^2 \lesssim \| u^{\mm{d}} \|_2^2+  F^{\mm{d} }_2 , \label{Lem:03sfdaxx}\\
&\frac{\mm{d}}{\mm{d}t}\|\nabla^{2-i} \partial_2^i(  u^{\mm{d}},m\partial_2   \eta^{\mm{d}})\|_0^2+\nu \|\partial_2^iu^{\mm{d}}\|_{{3-i}}^2 \lesssim  F^{\mm{d},i}_3, \label{qwesse6dfgssgsdfssd} \\
& \frac{\mm{d}}{\mm{d}t}\Big(\frac{\nu}{2}\| \nabla^{3-i}\partial_2^{i}\eta^{\mm{d}}\|_0^2+  \sum_{|\alpha|=2-i}\int \partial^{\alpha}\partial_2^{i}\eta^{\mm{d}}\cdot\partial^{\alpha}\partial_2^{i} u^{\mm{d}}\mm{d}y\Big)
\nonumber \\
&\quad + \|m \partial_2^{i+1} \eta^{\mm{d}} \|_{2-i}^2 \lesssim \|\partial_2^i u^{\mm{d}} \|_{2- i}^2 + F^{\mm{d},i}_4
 \label{qwesdfsf3xx}
\end{align}
for $i=1,2$, where
$$
\begin{aligned}
& F^{\mm{d}}_1:= \|\partial_2^2\eta\|_1^2 \| u\|_2^2+ \|\partial_2 \eta\|_2^2\|\partial_2 u\|_1^2
+ \|\partial_2\eta\|_2^2\|\nabla q \|_0^2\nonumber  \\
 &\qquad \quad + \|\nabla q \|_0  (\|\partial_2^2\eta\|_1 \| u\|_2+ \|\partial_2 \eta\|_2\|\partial_2 u\|_1) \nonumber ,\\
&F^{\mm{d}}_2:=
\|\partial_2\eta\|_2   (\|\nabla \eta^{\mm{d}}\|_2 \|\nabla u\|_2+ \|\nabla \eta \|_2 \|\nabla u^{\mm{d}} \|_2)
+ \|\partial_2\eta\|_2\|\nabla (\eta,\eta^{\mm{d}})\|_2\|\nabla q\|_1 , \nonumber \\
&F^{\mm{d},i }_3:=
                    \begin{cases}
                 \|\partial_2\eta\|_2^2(\| \partial_2 u\|_2^2 +\|\nabla q\|_1^2)& \hbox{ for }i=1; \\
  \|\partial_2^3\eta\|_0^2\|\partial_2 u\|_2^2+\|\partial_2^2\eta\|_1^2\|\partial_2u\|_2^2  + \|\partial_2\eta\|_2^2\|(\partial_2^2 u,\nabla  q) \|_1^2 &\\ + \|\nabla q\|_1(\| \partial_2^3 \eta\|_0 \|\partial_2  u \|_2 +\| \partial_2  \eta\|_2 \|\partial_2^2  u \|_1  )         & \hbox{ for }i=2,
                    \end{cases}\end{aligned}
$$
and
$$\begin{aligned}
&F^{\mm{d},i}_4:=
                    \begin{cases}
           \|\partial_2\eta^{\mm{d}}\|_2(\|  \partial_2^2 \eta\|_1 \|  \partial_2 u\|_2  +\|  \partial_2 \eta\|_2 \|  \partial_2^2 u\|_1) &\\
 + \|  \partial_2\eta\|_2 \|\partial_2^2(\eta, \eta^{\mm{d}})\|_1  \|\nabla q\|_1        & \hbox{ for }i=1; \\
 \|\partial_2^3\eta^{\mm{d}}\|_0
(\|\partial_2^2\eta  \|_1\|\partial_2u\|_2+\|\partial_2\eta\|_2\|\partial_2^2u\|_1)&\\+
 \|\nabla q\|_1( \|\partial_2^3 (\eta,\eta^{\mm{d}}) \|_0 \| \partial_2 \eta\|_2+\|\partial_2^2\eta\|_1^2)                & \hbox{ for }i=2.
                    \end{cases}  \nonumber
\end{aligned}
$$

If one integrates by parts, one gets from \eqref{202011032123}, \eqref{201702061safasaf4181721asfdx} and \eqref{Lem:03sfdaxx} that
\begin{align}
& \|\eta^{\mm{d}}\|_3^2+ \|(u^{\mm{d}},m \partial_2 \eta^{\mm{d}})\|_1^2+\int_0^t\|( u^{\mm{d}}, m\partial_2\eta^{\mm{d}})\|_2^2\mm{d}\tau
\nonumber\\
& \qquad \lesssim Cm^{-1}+\|\partial_2 \eta\|_2 \|\nabla \eta\|_2 \|\nabla \eta^{\mm{d}}\|_2 +   \int_0^t( F^{\mm{d}}_1 +F^{\mm{d} }_2)\mm{d}\tau .
\label{2019112safasdsaafd62012}
\end{align}
Recalling that $(\eta,u,q)$ satisfies \eqref{2017020614181721}, we make use of \eqref{20190safd5041053}, \eqref{2017020614181721} and \eqref{2020052010sdfa33} to get \eqref{201905041053xxx} from \eqref{2019112safasdsaafd62012}.

Following the arguments of deriving \eqref{2019112ssfaasfdfa62012}--\eqref{201911262012}, using \eqref{qwesse6dfgssgsdfssd}
and \eqref{qwesdfsf3xx}, we arrive at
\begin{align}
\frac{\mm{d}}{\mm{d}t}   (\langle t\rangle^i {\mathcal{E}}_i^{\mm{d}})+
c \langle t\rangle^i( \|\partial_2^i    u^{\mm{d}}\|_{3-i}^2+\|m \partial_2^{i+1} \eta^{\mm{d}}  \|_{2-i}^2)
  \lesssim \langle t\rangle^{i-1} {\mathcal{E}}_i^{\mm{d}} + \langle t\rangle^i  (  F^{\mm{d},i}_3 +F^{\mm{d},i}_4)
\label{201911262sadfsa012}
\end{align}
and \begin{align}
&\frac{\mm{d}}{\mm{d}t} \left(\langle t\rangle^{i+1} \|\nabla^{2-i} \partial_2^i(u^{\mm{d}},m\partial_2\eta^{\mm{d}} )\|_0^2 \right)
+\nu\langle t\rangle^{i+1}\| \partial_2^iu^{\mm{d}}\|_{3-i}^2\nonumber \\
&\qquad \lesssim\langle t\rangle^i  \| \nabla^{2-i}\partial_2^i(u^{\mm{d}},m\partial_2\eta^{\mm{d}})\|_0^2+
\langle t\rangle^{i+1} F_3^{\mm{d},i}\;\;\;\mbox{ for }i=1,\ 2,
\label{2019112asdfassdfa62012}
\end{align}
where
$$
\begin{aligned}
&{\mathcal{E}}_i^{\mm{d}}  := c \|\nabla^{2-i}\partial_2^i (u^{\mm{d}},m\partial_2  \eta^{\mm{d}})\|_0^2 +\frac{\nu}{2}\|\nabla^{3-i}\partial_2^{i}  \eta^{\mm{d}}\|_0^2+
\sum_{  |\alpha| = 2-i} \int  \partial^\alpha \partial^i_2\eta^{\mm{d}} \cdot \partial^\alpha\partial^i_2 u^{\mm{d}}\mm{d}y ,
\end{aligned}$$
and $ \mathcal{E}_i^{\mm{d}}  $ satisfies
$$ \|\partial_2^i(\nabla \eta^{\mm{d}}, u^{\mm{d}},m\partial_2\eta^{\mm{d}}) \|_{2-i}^2\lesssim   {\mathcal{E}}_i^{\mm{d}} \lesssim  \|\partial_2^i(\nabla \eta^{\mm{d}}, u^{\mm{d}},m\partial_2\eta^{\mm{d}}) \|_{2-i}^2. $$
Utilizing \eqref{202005201033}, \eqref{202011032123}, \eqref{2017020614181721} with $i=1$ and \eqref{2020052010sdfa33}, we easily obtain \eqref{2020sfa05201033} from \eqref{201911262sadfsa012} and \eqref{2019112asdfassdfa62012}.

Next, we further assume that $\eta^0$ satisfies the additional regularity condition $\partial_1( \eta^0,u^0)\in H^3\times H^2$.
Thus $\partial_1(\eta^{\mm{r}}, u^{\mm{r}})\in  \underline{H}^3\times \underline{H}^2$ satisfies the estimates:
 \begin{align} &\| \partial_1 u^{\mm{r}} \|_2\lesssim_0 \|\partial_1\partial_2 \eta^0\|_2\|  u^0\|_2+\|\partial_2 \eta^0\|_2\|\partial_1u^0\|_2 ,
\label{202011212202}\\
&\| \partial_1\partial_2 \eta^{\mm{r}} \|_2\lesssim_0 \| \partial_2 \eta^0\|_2\| \partial_1\partial_2\eta^0\|_2.\label{2020110321asdfas23}
\end{align}
Keeping in mind that $\partial_1(\tilde{\eta}^0,\tilde{u}^0)\in {H}^3\times {H}^2$, we can further get, using
\eqref{202011212202} and \eqref{2020110321asdfas23}, that
\begin{align}
 \|\partial_1^3 ( u^{\mm{L}}_2,m\partial_2\eta^{\mm{L}}_2)(t)\|_{0}^2+\nu \int_0^t\|\partial_1^3\nabla u^{\mm{L}}_2\|_0^2\mm{d}\tau
& =\|\partial_1^3 ( u^{\mm{L}}_2,m\partial_2\eta^{\mm{L}}_2)(t)\|_{0}^2|_{t=0} \nonumber \\
& \lesssim (1+\|\partial_1\eta^0\|_3^2+\|\partial_1 u^0\|_2^2). \label{2203x}
\end{align}

If we make use of the estimate
$$\begin{aligned}
&\int \partial_1^2\mathcal{A}_{k2}\partial_2u_k \partial_1^2q\mm{d}y=\frac{\mm{d}}{\mm{d}t}\int \partial_1^2\mathcal{A}_{k2}\partial_2  \eta_k\partial_1^2q \mm{d}y
-\int \partial_2 \eta_k(\partial_t\partial_1^2\mathcal{A}_{k2}\partial_1^2q+\partial_1^2\mathcal{A}_{k2}\partial_1^2q_t )\mm{d}y, \\
&-\int \partial_1^2\tilde{\mathcal{A}}_{12}\partial_2u\cdot \partial_1^3 u^{\mm{d}}\mm{d}y=\int \partial_1^3(\eta_2^{\mm{d}}+
\eta_2^{\mm{L}})\partial_2u\cdot \partial_1^3 u^{\mm{d}}\mm{d}y \mbox{ and }\|\partial_1^3\eta^{\mm{L}}\|_0
\lesssim_0\|\partial_1^3\partial_2 \eta^{\mm{L}}\|_0,
\end{aligned} $$
 we can deduce from \eqref{01dsaf16safafasdfasfxx} that
 \begin{align}&\frac{\mm{d}}{\mm{d}t}  \left(  \| \nabla^2 (u^{\mm{d}}, m \partial_2 \eta^{\mm{d}} )\|_0^2
 +\int\partial_1^2\mathcal{A}_{k2}\partial_2  \eta_k  \partial_1^2q\mm{d}y\right)+
 c\|  u^{\mm{d}} \|_3^2  \lesssim F^{\mm{d}}_5 , \label{201702061safa4181721asfdx} \end{align}
where
\begin{align}
&F^{\mm{d} }_5:= \|\partial_2^2\eta\|_2^2\| u\|_3^2
+ \|\partial_2\eta\|_2^2\|\partial_2 u\|_2^2 +\|  \partial_1^3 \eta^{\mm{d}}\|_0^2\|\partial_2^2u \|_1^2+
\|  \partial_1^3 \partial_2 \eta^{\mm{L}}_2\|_0^2\|\partial_2^2u \|_1^2 \nonumber \\
 &\quad \quad  +  \| \partial_2\eta\|_2^2  \|\nabla q\|_1^2 +
 \|\nabla q\|_1(\| \partial_2^2\eta\|_1\|u\|_3 +\| \partial_2\eta\|_2\|\partial_2 u\|_2  )
 +\|\partial_2\eta \|_2\|\nabla \eta\|_2\|\nabla q_t\|_1.  \label{202011220947}
\end{align}
Consequently, with the help of \eqref{20190safd5041053}, \eqref{202005201033}, and \eqref{201905041053xxx}, \eqref{2020103safda02130}
and \eqref{2203x}, we easily obtain \eqref{202011032108} from \eqref{201702061safa4181721asfdx}.
This completes the proof of Theorem \ref{201912041028}.

\section{Proof of Theorem \ref{2019043019sdfa48}}\label{201912asdfas062021}
This section is devoted to the proof of Theorem \ref{2019043019sdfa48}. The key ideas to establish Theorem \ref{2019043019sdfa48} are
 similar to those in the proof of Theorems \ref{201904301948}, \ref{202006211824saf} and \ref{201912041028}. We break up the proof
 of Theorem \ref{2019043019sdfa48} into two subsections.

\subsection{ {Existence and uniqueness of a global time-decay classical solution}}\label{201912asdfassfd062021}
    We start with the derivation of the  \emph{a priori} stability estimates
\eqref{202008071638} and \eqref{202911sfa9fasdfd52245} in Theorem \ref{2019043019sdfa48}. To this end, let $(\eta,u,q)$ be a solution
of the initial-value problem \eqref{01dsaf16asdfsfdsaasf}--\eqref{01dsafsfda16asdfasfsaf} defined on $\Omega_T$  for any given $T>0$,
where $(\eta^0,u^0)$ belongs to $\underline{H}^5_{1,2}\times \underline{H}^4$, and satisfies $\mm{div}_{\mathcal{A}^0}u^0=0$
and the odevity conditions \eqref{202007301500}. It should be remarked that the solution automatically satisfies  $(\eta)_{\mathbb{T}^2}=(u)_{\mathbb{T}^2}=0$ and the odevity conditions \eqref{20208192055}.
We further assume that $(\eta,u,q)$ and $K_\kappa$ satisfy $(q)_{\mathbb{T}^2}=0$,
\begin{align}
 \label{aprpiosesnewxxxxx}
 \sup_{0\leqslant t\leqslant T}\sqrt{\| ( \eta,u, m\partial_2 \eta)(t)\|_4^2} \leqslant K_\kappa\;\;\mbox{ for any given }T>0,
\end{align}
and
\begin{equation}
\label{aprpiosesnewxxasfedsasxxzx}
\max\left\{ K_\kappa^{1/2},K_\kappa^2\right\}/m\leqslant \delta ,
\end{equation}
  where $K_\kappa$ will be given in \eqref{202008071646} and $\delta$ is a sufficiently small constant.
  We should point out here that the smallness of $\delta$ depends only on $\nu$ and $\kappa$.

 By virtue of \eqref{aprpiosesnewxxxxx} and \eqref{aprpiosesnewxxasfedsasxxzx}, it is easy to see that
\begin{equation}
\label{aprpiosesnewxxxx1xasdfsa}
 \sup_{0\leqslant t\leqslant  T} \|  \partial_2 \eta (t)\|_4  \leqslant  \delta
\end{equation}
and
\begin{equation}
\label{aprp}
 \sup_{0\leqslant t\leqslant  T}( m^{-1}\| \eta  (t)\|_4 +  \|u(t)\|_4\|\partial_2 \eta (t)\|_4 )\lesssim_0 \delta.
\end{equation}
Now, we establish some basic energy estimates for $(\eta,u,q)$ which will be used later.
\begin{lem}\label{20161xx}
Under the condition \eqref{aprpiosesnewxxxx1xasdfsa} with sufficiently small $\delta$, one has
\begin{align}
&\|  \eta\|_{i+2,2} \lesssim_0 \|\partial_2 \eta\|_{i+1},\label{2020008061647} \\
&\|\nabla q\|_3 \lesssim
   \|  u\|_3^2+m^2 \|\partial_2 \eta\|_3\|\partial_2 \eta\|_4 , \label{2017asdfsaf21} \\
& \|  \chi\|_{i+1}\lesssim _0\|\nabla^i \mm{curl}_{\mathcal{A}}\chi\|_0\lesssim_0 \|\chi\|_{i+1},\quad 0\leqslant i\leqslant 3,
 \label{20170206141817asdfasdfsaf21}
\end{align}
where $\chi$ denotes $\eta$, or $\partial_2 \eta$, or $u$.
\end{lem}
\begin{pf}
 {We can argue in the same manner as in the derivation of \eqref{20202111191231} and \eqref{2017020614181721} with $i=1$ to deduce}
 \eqref{2020008061647} and \eqref{2017asdfsaf21}. To show \eqref{20170206141817asdfasdfsaf21}, keeping in mind that
$\Delta \chi=\nabla^\bot \mm{curl}\chi+\nabla \mm{div}\chi$ where $\mm{curl} \chi=\partial_1 \chi_2- \partial_2 \chi_1$ and
$\nabla^\bot :=(-\partial_2,\partial_1)$, we easily obtain
\begin{align}
\|\nabla \chi\|_i^2= \| \mm{curl}\chi\|_i^2+\| \mm{div}\chi\|_i^2.
\label{202008061631}
\end{align}

Recalling the product estimate
\begin{align}   \label{fasdfagestims}
&  \|fg\|_{j}\lesssim      \begin{cases}
                      \|f\|_{1}\|g\|_{1} & \hbox{ for }j=0,  \\
  \|f\|_j\|g\|_{2} & \hbox{ for }0\leqslant j\leqslant 2,  \\
                    \|f\|_{2}\|g\|_j+\|f\|_j\|g\|_{2}& \hbox{ for }3\leqslant j\leqslant 4 ,
                    \end{cases}
\end{align}
and noting that $(\eta,u)$ satisfies \eqref{201903211118} and \eqref{202005011525}, we find that
\begin{align}
\| \mm{div}\chi\|_i\lesssim_0   \|\eta\|_{4}\|\nabla \chi\|_i,\qquad 0\leqslant i\leqslant 3.
\label{20201122}
\end{align}
Thus, with the help of \eqref{20160614fdsa19asfda57x}, \eqref{aprpiosesnewxxxx1xasdfsa} and \eqref{2020008061647},
we get from \eqref{202008061631} and \eqref{20201122} that
\begin{align}
\|\chi\|_{i+1}\lesssim_0\|\nabla^i \mm{curl} \chi\|_0.
\label{202010122009}
\end{align}
In addition, by \eqref{aprpiosesnewxxxx1xasdfsa}, \eqref{2020008061647} and \eqref{fasdfagestims},
\begin{align}
\|\nabla^i \mm{curl} \chi\|_0\lesssim_0 \|\nabla^i \mm{curl}_{\mathcal{A}} \chi\|_0\lesssim_0 \|\chi\|_{i+1},
\end{align}
which, together with \eqref{202010122009}, yields \eqref{20170206141817asdfasdfsaf21}. \hfill $\Box$
\end{pf}
\begin{lem}\label{211qsadfwe201612132242nn}
Under the condition \eqref{aprpiosesnewxxxx1xasdfsa} with sufficiently small $\delta$, one has
\begin{align}
&  \frac{\mm{d}}{\mm{d}t}    \|\nabla^3 \mm{curl}_{\mathcal{A}} (u,m \partial_2 \eta) \|_0^2 + \kappa \|   u \|_4^2  \lesssim_0  \|  (u,m\partial_2\eta)\|_3\| ( u,m\partial_2\eta)\|_4^2 \label{211qwessebsadd}
 \end{align}
and
\begin{align}
 & \frac{\mm{d}}{\mm{d}t}\Big( \sum_{|\alpha|=3}\int\partial^\alpha \mm{curl}_{\mathcal{A}}\eta \partial^\alpha\mm{curl}_{\mathcal{A}}u\mm{d}y
  +\frac{\kappa}{2}\| \nabla^3 \mm{curl}_{\mathcal{A}}\eta\|_0^2  \Big) + \|m \eta \|_{5,2}^2 \lesssim_\kappa\|u\|_4^2 .
\label{211qsdafweLem:03xx}
\end{align}
\end{lem}
\begin{pf}
One can not directly apply $\partial^\beta$ to \eqref{01dsaf16asdfsfdsaasf}$_2$ to derive \eqref{211qwessebsadd} and \eqref{211qsdafweLem:03xx}
for the case $|\beta|=4$,  {since $q$ does not have the fifth-order derivatives.} Instead,
we apply $\partial^\alpha \mm{curl}_{\mathcal{A}}$ to \eqref{01dsaf16asdfsfdsaasf}$_2$ to get
\begin{equation}\label{01dsaf16asdsafdafsfdsaasf}
\partial^\alpha \mm{curl}_{\mathcal{A}} u_t+  \kappa \partial^\alpha \mm{curl}_{\mathcal{A}} u=
  m^2\partial^\alpha \mm{curl}_{\mathcal{A}} \partial_2^2\eta ,
\end{equation}
where $|\alpha|=3$
\footnote{We explain why to take $\alpha$ satisfying $|\alpha|=3$.
Let us consider the first term $ \int \partial^\alpha\mm{curl}_{\mathcal{A}_t}u\partial^\alpha\mm{curl}_{\mathcal{A}}u\mm{d}y$
on the right hand of \eqref{211qwe201808sdfsadfafaas100154628xxx}. Obviously,
$\int \partial^\alpha\mm{curl}_{\mathcal{A}_t}u\partial^\alpha\mm{curl}_{\mathcal{A}}u\mm{d}y
\lesssim_0 \|\nabla u\|_{L^\infty} \|\partial^\alpha u\|_0^2$. Since $H^2\hookrightarrow L^\infty$, obviously, we need at least
$|\alpha|\geqslant 3$. However, for $|\alpha|=3$,
it seems difficult to close the energy estimates by \eqref{211qwe201808sdfsadfafaas100154628xxx} and \eqref{211qwe20sadfa1808sdfaas100154628xxx}.
To overcome this difficulty, we adopt the two-layers energy method, i.e., we close the lower-order
energy estimate $\|\nabla (u,m\partial_2\eta)\|_2$ by a lower-order energy inequality, and then further close the highest-order energy estimate
by a highest-order energy inequality. Thus apparently, we need at least $|\alpha|\geqslant 4$.
Here we remark that the initial-value problem \eqref{01dsaf16asdfsfdsaasf} with small initial data $(\eta^0,u^0)
\in (\underline{H}^{4}_{1,2}\cap H^3_*)\times \underline{H}^3$ admits a unique global strong solution
$(\eta ,u,q)\in  \underline{\mathfrak{H}}^{4,*}_{1,2,\infty}\times C^0(\mathbb{R}_0^+,\underline{H}^3\times\underline{H}^3)$. }.
Multiplying \eqref{01dsaf16asdsafdafsfdsaasf} by $\partial^\alpha \mm{curl}_{\mathcal{A}}u$, resp.  $\partial^\alpha \mm{curl}_{\mathcal{A}}\eta$,
in $L^2$, we have
\begin{align}
& \frac{1}{2}\frac{\mm{d}}{\mm{d}t}  \|\partial^\alpha  \mm{curl}_{\mathcal{A}}(u, m \partial_2\eta)\|_0^2   +
  \kappa \|  \partial^\alpha  \mm{curl}_{\mathcal{A}} u \|_0^2 \nonumber \\
&=   \int \partial^\alpha  \mm{curl}_{\mathcal{A}_t}u\partial^\alpha  \mm{curl}_{\mathcal{A}}u\mm{d}y +m^2
\int (\partial^\alpha \mm{curl}_{\mathcal{A}}\partial_2\eta
\partial^\alpha \mm{curl}_{\mathcal{A}_t}\partial_2\eta\nonumber
 \\
&\quad-\partial^\alpha \mm{curl}_{\partial_2\mathcal{A}}\partial_2\eta
\partial^\alpha \mm{curl}_{\mathcal{A}}u-\partial^\alpha \mm{curl}_{\mathcal{A}}\partial_2\eta
\partial^\alpha \mm{curl}_{\partial_2\mathcal{A}}u)
\mm{d}y=:I_9, \label{211qwe201808sdfsadfafaas100154628xxx}
\end{align}
resp.
\begin{align}
&  \frac{\mm{d}}{\mm{d}t}\left(    \int  \partial^\alpha \mm{curl}_{\mathcal{A}} u\partial^\alpha \mm{curl}_{\mathcal{A}}\eta \mm{d}y +\frac{\kappa}{2}\| \partial^\alpha \mm{curl}_{\mathcal{A}}\eta\|_0^2  \right)
 +    \| m\partial^\alpha\mm{curl}_{\mathcal{A}} \partial_2\eta \|_0^2\nonumber \\
& =\int ( \partial^\alpha \mm{curl}_{\mathcal{A}} u \partial^\alpha \mm{curl}_{\mathcal{A}_t}\eta + (\partial^\alpha \mm{curl}_{\mathcal{A}}u)^2  + \partial^\alpha \mm{curl}_{\mathcal{A}_t} u \partial^\alpha \mm{curl}_{\mathcal{A}}\eta )\mm{d}y\nonumber \\
&\quad  -  m^2\int \partial^\alpha \mm{curl}_{\mathcal{A}}\partial_2 \eta  \partial^\alpha \mm{curl}_{\partial_2 \mathcal{A}}\eta  \mm{d}y  + \kappa
\int \partial^\alpha \mm{curl}_{\mathcal{A}}\eta  \partial^\alpha \mm{curl}_{\mathcal{A}_t}\eta \mm{d}y=:I_{10}, \label{211qwe20sadfa1808sdfaas100154628xxx}
\end{align}
where $I_{9}$ and $I_{10}$ can be bounded as follows.
\begin{align}
&I_{9}\lesssim_0 \|\nabla (u,m\partial_2\eta)\|_2\|\nabla (u,m\partial_2\eta)\|_3^2  ,\nonumber \\
& I_{10}\lesssim_\kappa (1+\|\nabla \eta\|_3
)  \|\nabla u\|_3^2 + m^2\|\nabla \eta\|_3\|\nabla \partial_2\eta\|_3^2+  \| \nabla u\|_3\|\nabla \eta\|_3^2.\nonumber
\end{align}
Consequently, inserting the above two estimates into \eqref{211qwe201808sdfsadfafaas100154628xxx} and
\eqref{211qwe20sadfa1808sdfaas100154628xxx} respectively, using \eqref{aprpiosesnewxxxx1xasdfsa}--\eqref{2020008061647}
and Young's inequality, we obtain \eqref{211qwessebsadd} and \eqref{211qsdafweLem:03xx}.
\hfill$\Box$
\end{pf}

\begin{lem}\label{211qwe201612132242nn}
Under the conditions \eqref{aprpiosesnewxxxx1xasdfsa}--\eqref{aprp} with sufficiently small $\delta$, we have
\begin{align}
& \frac{\mm{d}}{\mm{d}t}    \|\nabla^3( u, m   \partial_2 \eta )\|_0^2 + \kappa \|  u \|_3^2  \lesssim_0 \delta \|m \eta\|_{4,2}^2
\label{211qwessebsaddfdaiseqinM0846dfgssgsd}
 \end{align}
and
\begin{align}
 & \frac{\mm{d}}{\mm{d}t}\Big( \sum_{|\alpha|=3}  \int  \partial^\alpha  \eta \cdot \partial^\alpha  u\mm{d}y
 +\frac{\kappa}{2}\| \nabla^3  \eta\|_0^2 \Big) + \|m \eta \|_{4,2}^2 \lesssim_0 \|u \|_3^2  .
\label{211qweasdfLem:03xx}
\end{align}
\end{lem}
\begin{pf} It is easy to see from the proof of Lemma \ref{211qsadfwe201612132242nn} that
 {we can not directly establish the estimates}
\eqref{211qwessebsaddfdaiseqinM0846dfgssgsd} and \eqref{211qweasdfLem:03xx} by the $\mm{curl}_{\mathcal{A}}$-estimate method,
since we shall use some smallness properties which are hidden in some integral terms involving $\nabla_{\mathcal{A}}q$.

To exploit the smallness properties hidden in the pressure term $\nabla_{\mathcal{A}}q$, we apply $\partial^\alpha$ ($|\alpha|=3$)
to \eqref{01dsaf16asdfsfdsaasf}$_2$ to get
\begin{equation}\nonumber
\partial^\alpha u_t + \partial^\alpha\nabla_{ {\mathcal{A}}}  q+\kappa \partial^\alpha u=m^2\partial^\alpha  \partial_2^2\eta .
\qquad ( |\alpha|=3 )
\end{equation}
Thus, multiplying the above identity by $\partial^{\alpha} u$ and $\partial^{\alpha}\eta$ in $L^2$ respectively, we have
\begin{align}
& \frac{1}{2}\frac{\mm{d}}{\mm{d}t}   \|\partial^\alpha  (u,m \partial_2\eta)\|_0^2   +
  \kappa \|  \partial^\alpha  u \|_0^2   =I_{11}:= -\int  \partial^{\alpha}  \nabla_{  {\mathcal{A}}} q \cdot   \partial^{\alpha }  {u} \mm{d}y, \label{218safxxx}
\end{align}
and
\begin{align}
&  \frac{\mm{d}}{\mm{d}t}\left(    \int  \partial^\alpha    \eta \cdot \partial^\alpha u\mm{d}y +\frac{\kappa}{2}\| \partial^\alpha \eta\|_0^2  \right)  +   \|m \partial^\alpha\partial_2  \eta \|_0^2-\|\partial^\alpha  u\|_0^2\nonumber \\
& =I_{12}:=- \int \partial^{\alpha} \nabla_{  {\mathcal{A}}} q \cdot \partial^{\alpha }  \eta \mm{d}y , \label{21safa1qw8xxx}
\end{align}
respectively, where the right-hand sides can be bounded as follows, using \eqref{201903211118}, \eqref{202005011525} and \eqref{2017asdfsaf21}.
\begin{eqnarray*}
&& I_{11}\lesssim_0  \|  \nabla \eta\|_3\|\nabla u\|_2  (\|  u\|_3^2+ m^2 \|\partial_2 \eta\|_3\|\partial_2 \eta\|_4 ) ,  \\
&& I_{12}\lesssim_0 \|  \nabla \eta\|_3\|\nabla \eta\|_2(\|u\|_3^2+m^2 \|\partial_2 \eta\|_3 \|\partial_2 \eta\|_4 ) ,
\end{eqnarray*}
which, together with \eqref{218safxxx}, \eqref{21safa1qw8xxx} and \eqref{aprpiosesnewxxxx1xasdfsa}--\eqref{2020008061647}, implies \eqref{211qwessebsaddfdaiseqinM0846dfgssgsd} and \eqref{211qweasdfLem:03xx}.
\hfill$\Box$
\end{pf}

Now, we are in a position to show the \emph{a priori} stability estimate \eqref{202008071638}.
Firstly, using \eqref{20170206141817asdfasdfsaf21}, we derive from Lemmas \ref{211qsadfwe201612132242nn}--\ref{211qwe201612132242nn} that
\begin{align}
& \frac{\mm{d}}{\mm{d}t}  \mathcal{E}^\kappa +  c^\kappa(\|  u\|_4^2+\|m\eta \|_{5,2}^2 )  \lesssim_\kappa    \| (u, m \partial_2\eta)\|_3^2\| (u,m\partial_2\eta)\|_4^2 ,
\label{211qwefM0846dfgssgsd} \\
&  \frac{\mm{d}}{\mm{d}t}\mathcal{E}^\kappa_l+ c^\kappa(\|   u \|_3^2+\| m   \eta\|_{4,2}^2)\leqslant 0,
\label{202008061sfa938}
\end{align}
where
\begin{align}
&\mathcal{E}^\kappa:=c \|\nabla^3\mm{curl}_{\mathcal{A}} (u,m\partial_2 \eta)\|_0^2+\sum_{|\alpha|=3}\int\partial^\alpha \mm{curl}_{\mathcal{A}}\eta \partial^\alpha  \mm{curl}_{\mathcal{A}} u\mm{d}y +\frac{\kappa}{2}\| \nabla^3 \mm{curl}_{\mathcal{A}}\eta\|_0^2,
\label{202012051924} \\
&\mathcal{E}^\kappa_l:= c\|\nabla^3( u, m   \partial_2 \eta )\|_0^2 + \sum_{|\alpha|=3}  \int  \partial^\alpha  \eta \cdot \partial^\alpha  u\mm{d}y +\frac{\kappa}{2}\| \nabla^3  \eta\|_0^2  , \nonumber \\
& \|(\eta, u)\|_4^2+\|m\eta \|_{5,2}^2   \lesssim_\kappa \mathcal{E}^\kappa \lesssim_\kappa\|( \eta,u,m\partial_2\eta) \|_{4}^2,\nonumber \\
& \| (\eta, u)\|_3^2+\| m  \eta\|_{4,2}^2 \lesssim_\kappa \mathcal{E}^\kappa_l\lesssim_\kappa\|( \eta,u,m\partial_2\eta) \|_3^2 . \nonumber
\end{align}
Thus, we further obtain
\begin{align}
& \frac{\mm{d}}{\mm{d}t}e^{ c_3^\kappa\min\{1,m\}  t} \mathcal{E}^\kappa + c^\kappa e^{ c_3^\kappa\min\{1,m\} t} (\|  u\|_4^2
+\|m\eta \|_{5,2}^2 )\nonumber \\
&\qquad \lesssim_\kappa   e^{ c_3^\kappa\min\{1,m\} t} \| (u, m \partial_2\eta)\|_3^2\| (u,m\partial_2\eta)\|_4^2 ,
\label{211qwessebsaddfasfdadaiseqinM0846dfgssgsd} \\
& \frac{\mm{d}}{\mm{d}t} (\langle t\rangle^i\mathcal{E}^\kappa) + c^\kappa\langle t\rangle^i(\|  u\|_4^2+\|m\eta \|_{5,2}^2 ) \nonumber \\
&\qquad\lesssim_\kappa \langle t\rangle^{i-1} \mathcal{E}^\kappa+ \langle t\rangle^i \| (u, m \partial_2\eta)\|_3^2
\| (u,m\partial_2\eta)\|_4^2\;\;\;\mbox{ for any }i\geqslant 1, \label{211qwessgsd} \\
& \|(\eta,u)\|_3^2+\| m\eta\|_{4,2}^2 +\int_0^t(\|u\|_3^2+\|m\eta\|_{4,2}^2)\mm{d}\tau\lesssim_\kappa\|(\eta^0,u^0, m\partial_2\eta^0\|_3^2.
\label{202008061938}
\end{align}

Let $c_4^\kappa$ be the positive constant in \eqref{202012212151}. In view of
\eqref{211qwessebsaddfasfdadaiseqinM0846dfgssgsd} and \eqref{202008061938}, we find that
\begin{align}
 e^{c_3^\kappa\min\{1,m\} t}(\|(\eta,u)\|_4^2+\|m\eta\|_{5,2}^2) \leqslant {c^\kappa_1} \|(\eta^0,u^0, m\partial_2 \eta^0)
 \|_4^2e^{c^\kappa_2\|  (\eta^0,u^0, m\partial_2 \eta^0) \|_3^2 }/c_4^\kappa,  \label{202008061949}
\end{align}
where $c^\kappa_1\geqslant 4$.
From now on, we take
\begin{align}
K_\kappa=  \sqrt{ {c^\kappa_1}  \|  (\eta^0,u^0, m  \partial_2 \eta^0) \|_4^2e^{c^\kappa_2\|  (\eta^0,u^0, m  \partial_2 \eta^0) \|_3^2 }}.
\label{202008071646}
\end{align}
Thus, there is a $\delta_1^\kappa\in (0,1]$, such that for any $\delta\leqslant \delta_1^\kappa$,
\begin{align}
\|(\eta,u)\|_4^2+\|m\eta\|_{5,2}^2\leqslant  e^{c_3^\kappa\min\{1,m\} t}(\|(\eta,u)\|_4^2+\|m\eta\|_{5,2}^2)   \leqslant K_\kappa^2/c_4^\kappa.  \label{202949}
\end{align}
Besides, we get from \eqref{211qwessebsaddfasfdadaiseqinM0846dfgssgsd} and \eqref{202949} that
\begin{align}
& \int_0^t e^{ c_3^\kappa \min\{1,m\} \tau }(\| u\|_4^2+\| m \eta \|_{5,2}^2) \mm{d}\tau  \lesssim_\kappa  K_\kappa^2/c_4^\kappa,  \label{20asdfas2008061949}
\end{align}
which, together with \eqref{202949}, yields the desired stability estimate:
\footnote{We can get the exponential decay in the fourth-order energy inequality \eqref{20200safda071638}, since the energy functional
$\mathcal{E}^\kappa$ is equivalent to the dissipative functional $\mathcal{D}^\kappa$ in \eqref{211qwefM0846dfgssgsd}
where $\mathcal{D}^\kappa:=\|u\|_4^2+\|m\eta \|_{5,2}^2 $. Obviously, for the viscous fluid case, a similar equivalent relation
does not hold in the second-order energy inequality \eqref{2019112safd62012}, in which one sees that the second-order energy functional
${\mathcal{E}}$ only controls the second-order dissipative functional $\mathcal{D}:=\|u\|_3^2+\| m\partial_2\eta\|_2^2$. Similarly,
if we further derive the first-order energy inequality, then the first-order energy functional, denoted by ${\mathcal{E}}^1$,
also controls the the first-order dissipative functional. However, the first-order energy functional ${\mathcal{E}}^1$ can be controlled
by the second-order dissipative functional ${\mathcal{D}}$; and this relation admits us to expect at most algebraic decay rates
in the viscous fluid case.}
\begin{align}
 & e^{  c_3^\kappa \min\{1,m\} t}(\| (\eta,u)\|_4^2+\| m  \eta\|_{5,2}^2) +\int_0^t e^{c_3^\kappa\min\{1,m\} \tau}(\| u\|_4^2+\| m \eta \|_{5,2}^2) \mm{d}\tau \lesssim_\kappa  {K}^2_\kappa.     \label{20200safda071638}
\end{align}

Next, we proceed to derive algebraic time-decay stated in \eqref{202911sfa9fasdfd52245}. To begin with,
using \eqref{aprpiosesnewxxasfedsasxxzx} and \eqref{20asdfas2008061949}, we have
\begin{equation}
\int_0^t  \|   \eta \|_{5,2}^2  \mm{d}\tau \lesssim  {K}_\kappa. \label{202010sdaf81313}
\end{equation}

Applying Gronwall's lemma to \eqref{211qwessgsd} and using \eqref{202008061938}, we conclude
\begin{align}
& \langle t\rangle^i(\|( \eta,u)\|_4^2+\|m\eta \|_{5,2}^2)\nonumber \\
& \lesssim_\kappa \left(\|(\eta^0, u^0,m\partial_2\eta^0)\|_{4}^2+\int_0^t \langle \tau\rangle^{i-1}\|( \eta,u,m\partial_2\eta) \|_{4}^2\mm{d}\tau\right)
  e^{c^\kappa_2( \|( \eta^0,u^0, m  \partial_2 \eta^0  )\|_3^2}.
  \label{202911952245}
\end{align}
Thanks to \eqref{202008061938} and \eqref{202911952245}, we further deduce from \eqref{211qwessgsd} that
\begin{align}
&\int_0^t\langle  \tau\rangle^i(\|  u\|_4^2+\|m\eta \|_{5,2}^2 ) \mm{d}\tau \nonumber \\
&\lesssim_\kappa  (1+ ( \|( \eta^0,u^0, m  \partial_2 \eta^0  )\|_3^2)\nonumber \\
&\qquad \times \left(\|( \eta^0,u^0,m\partial_2\eta^0)\|_{4}^2+\int_0^t \langle \tau\rangle^{i-1}\|( \eta,u,m\partial_2\eta) \|_{4}^2\mm{d}\tau\right)  e^{c^\kappa_2\|( \eta^0, u^0, m  \partial_2 \eta^0  )\|_3^2}.
\label{202011081439}
 \end{align}
Consequently, we conclude from \eqref{aprpiosesnewxxasfedsasxxzx}, \eqref{20asdfas2008061949}, \eqref{202010sdaf81313}
and \eqref{202011081439} that for $i=0$ and $1$,
\begin{align}
 \int_0^t\langle  \tau\rangle^i\big( \| u\|_4^2+(1+m^2)\|\eta \|_{5,2}^2\big)\mm{d}\tau \leqslant C^\kappa.
\label{202sdfas00110814002}
 \end{align}
Finally, we get from \eqref{202911952245} and \eqref{202sdfas00110814002} that for any $i =1$ and $2$,
\begin{align}
& \langle t\rangle^i(\|( u, \eta)\|_4^2+\|m\eta \|_{5,2}^2)  \leqslant C^\kappa,\nonumber
 \end{align}
which, together with \eqref{202sdfas00110814002}, yields \begin{align}
 \sum_{0\leqslant i\leqslant 1}\left(\langle t\rangle^i(\|( \eta,u )\|_4^2+\|m\eta \|_{5,2}^2)
 + \int_0^t\langle  \tau\rangle^i( \|  u\|_4^2+ \| m\eta \|_{5,2}^2) \mm{d}\tau\right)
+ \int_0^t\langle  \tau\rangle \|\eta \|_{5,2}^2\mm{d}\tau\leqslant C^\kappa.
  \label{202911sfa9fsafasfasdfd52245}
\end{align}

Next, we introduce a global well-posedness result for the linear initial-value problem \eqref{01dsaf16asdfsfsdafafdsaasf} and
 a local well-posedness result for the nonlinear initial-value problem \eqref{01dsaf16asdfsfdsaasf}--\eqref{01dsafsfda16asdfasfsaf}.
\begin{pro}
\label{2020212240930}
Let $i\geqslant 1$ be an integer, $\kappa$ and $m$ be positive constants. If
$(\tilde{\eta}^0,\tilde{u}^0)\in (H^{i+1}_2\cap H^i_\sigma)\times H^i_\sigma $, then there exists a unique strong solution $(\eta^{\mm{L}},u^{\mm{L}}) \in\mathfrak{C}^0(\mathbb{R}^+,{H}^{i+1}_2 )\times {\mathfrak{U}}_\infty^i$
to the following linear initial-value problem:
 \begin{equation}
 \begin{cases}
\eta_t^{\mm{L}}=u^{\mm{L}}  , \\[1mm]
u_t^{\mm{L}}  -\kappa u^{\mm{L}} = m^2\partial_2^2\eta^{\mm{L}}  , \\[1mm]
\div u^{\mm{L}}  =0,\\
(\eta^{\mm{L}} ,u^{\mm{L}})|_{t=0}=(\tilde{\eta}^0, \tilde{u}^0).
\end{cases} \label{202012241017}
\end{equation}
\end{pro}
\begin{pf}
The proof of Proposition \ref{2020212240930} is trivial, and hence we omit it here.
\end{pf}
\begin{pro}\label{pro:0401nsadfaxsdfafdsaddfdxx}
Let $B^\kappa> 0$ and $(\eta^{\mm{L}},u^{\mm{L}})\in\mathfrak{C}^0( {\mathbb{R}^+},{H}^5_2 )\times {\mathfrak{U}}_\infty^4$ be
the unique global solution of \eqref{202012241017} with $(\tilde{\eta}^0,\tilde{u}^0)\in (H^{5}_2\cap H^4_\sigma)\times H^4_\sigma $.
Assume that $(  \eta^0,u^0)\in H^5_{2} \times H^4$,
$\|(u^0,\partial_2\eta^0)\|_4\leqslant B^\kappa$ and $\mm{div}_{\mathcal{A}^0}u^0=0$,
where $\mathcal{A}^0$ is defined by $\zeta^0$ and $ \zeta^0 = \eta^0+y$. Then, there is a constant $\delta_2^\kappa\in (0,1]$,
such that if, in addition,
\begin{align}
\|\nabla \eta^0\|_3\leqslant \delta_2^\kappa,
\label{201912saddsdaffafdads261426}
\end{align}
the initial value problem \eqref{01dsaf16asdfsfdsaasf}--\eqref{01dsafsfda16asdfasfsaf} possesses a unique local classical solution
$( \eta, u,q)\in \mathfrak{ {C}}^0( {I_T},{H}^5_2 )\times\mathfrak{U}_T\times ({C}^0(\overline{I_T},\underline{H}^3 )\cap\mathfrak{C}^0(I_T,H^4))$
 for some $T>0$ dependent of $B^\kappa$, $\kappa$, $m$ and $\delta_2^\kappa$. Moreover, $(\eta ,u)$ satisfies
\begin{eqnarray} &&
0<\inf_{(y,t)\in \mathbb{R}^2\times \overline{I_T}} \det(\nabla \eta+I),\qquad \sup_{t\in \overline{I_T}} \|\nabla \eta\|_{3}\leqslant 2\delta_2^\kappa ; \label{202012301432}  \\
&& \label{202012212151}
  \sup_{t\in \overline{I_T}}\| (\eta,  u,  m \partial_2 \eta )\|_{4} \leqslant c^\kappa_4 \|(\eta^0,u^0,m \partial_2 \eta^0 )\|_{4};
\end{eqnarray}
and
\begin{align}
  &\sup_{\tau\in \overline{I_t}}\| (\mm{curl}_{\mathcal{A}} \eta ^{\mm{d}}, \mm{curl}_{\mathcal{A}} u^{\mm{d}} , \mm{curl}_{\mathcal{A}}m \partial_2 \eta ^{\mm{d}} )(\tau)\|_{ 3 }^2 \nonumber  \\ &
\lesssim_\kappa  e^{ t  \|(\eta^0,u^0,m \partial_2 \eta^0 )\|_{4} } \|(\eta^0-\tilde{\eta}^0,u^0-\tilde{u}^0,m\partial_2(\eta^0-\tilde{\eta}^0)\|_{4}^2
+ t\Big( \|( u^0,m \partial_2 \eta^0, \tilde{u}^0,m \partial_2 \tilde{\eta}^0 )\|_{4}^3  \nonumber\\
& \qquad +  \|( u^0,m \partial_2 \eta^0 )\|_{4}^2  \|( u^0,m \partial_2 \eta^0, \tilde{u}^0,m \partial_2 \tilde{\eta}^0 )\|_{4}^2\Big),
 \quad \forall\, t\in \overline{I_T}
\label{20201221asfdasf2151}
\end{align}
where the constant $c_4^\kappa\geqslant 1$ depends on $\kappa$ at least, and $(\eta^{\mm{d}},u^{\mm{d}} ):=(\eta-\eta^{\mm{L}},u-u^{\mm{L}})$.
\end{pro}
\begin{pf}
We postpone the proof of Proposition \ref{pro:0401nsadfaxsdfafdsaddfdxx} to Section \ref{202008062143}.
\hfill $\Box$
\end{pf}
\begin{rem}
\label{2020102sdad52143}
If $(\eta^0,u^0)$ in Proposition \ref{pro:0401nsadfaxsdfafdsaddfdxx} further satisfies
$(\eta^0,u^0)\in \underline{H}^{5}_{1,2}\times \underline{H}^4$ and the odevity conditions \eqref{202007301500},
then for each fixed $t\in (0,T]$, $(\eta,u)(t)$ also belongs to $\underline{H}^{5}_{1,2}\times \underline{H}^4 $
and satisfies the odevity conditions \eqref{20208192055}.
\end{rem}
 \begin{rem}
\label{202010safa2sdad52143}
Since $(u,\partial_2\eta)$ may do not belong to $C^0(\overline{I_T},H^4)$, one needs the additional estimates \eqref{202012212151}
and \eqref{20201221asfdasf2151} to further establish the existence result of a global solution and its asymptotic behavior with respect to $m$,
repspectively.
\end{rem}

With the \emph{priori} estimate \eqref{202949}  and Propositions \ref{pro:0401nasfxdxx} and \ref{pro:0401nsadfaxsdfafdsaddfdxx} in hand, we can easily establish the existence and uniqueness of a global time-decay smooth solution stated in Theorem \ref{2019043019sdfa48}.
Next, we sketch the proof.

Let $m$ and $(\eta^0,u^0)\in (\underline{H}^{5}_{1,2}\cap H^4_*)\times \underline{H}^4$ satisfy the odevity conditions \eqref{202007301500}
and
\begin{align}
\max\{K_\kappa^{1/2},K^2_\kappa\}/m\leqslant\min\{\delta_1^\kappa,\delta_2^\kappa/c_{0},\delta_3/c_{0}\}=:c_3^\kappa\leqslant 1,
\label{201912261safd425}
\end{align}
where $K_\kappa$ is defined by \eqref{202008071646},  {and the constant $c_0$ is the same as in \eqref{2020008061647}. }
Thus, $\eta^0$ satisfies \eqref{201912saddsdaffafdads261426} by virtue of \eqref{2020008061647} and \eqref{201912261safd425}.
Moreover, by Proposition \ref{pro:0401nsadfaxsdfafdsaddfdxx} and Remark \ref{2020102sdad52143}, one sees that there is a unique local solution
$( \eta, u,q)$ of \eqref{01dsaf16asdfsfdsaasf}--\eqref{01dsafsfda16asdfasfsaf} defined on a maximal existence time interval $[0,T^{\max})$,
such that
\begin{itemize}
   \item for any $T\in  I_{T^{\max}}$, the solution $(\eta,u,q)$ belongs to
   $\underline{\mathfrak{H}}^{5,*}_{1,2,T}\times\underline{\mathfrak{U}}_T\times
   ({C}^0(\overline{I_T},\underline{H}^3)\cap\mathfrak{C}^0(I_T,H^4))$ and satisfies
$ \sup_{t\in \overline{I_{T}}}\|\nabla \eta\|_{3}\leqslant 2\delta_2^\kappa$;
   \item $\limsup_{t\to T^{\max} }\|\nabla \eta( t)\|_{3} > \delta_2^\kappa$, or $\limsup_{t\to T^{\max} }\|(u,\partial_2\eta)( t)\|_4=\infty$.
 \end{itemize}
In addition, the solution enjoys the odevity conditions \eqref{20208192055}.

Let
\begin{equation}
\nonumber
T^{*}=\sup\left\{ T \in I_{T^{\max}}~\left|~\|(\eta,u,m\partial_2\eta\|_{4}  \leqslant K_\kappa \mbox{ for any }t\leqslant T\right.\right\}.
\end{equation}
Recalling the definition of $K_\kappa$ and the condition $c_1^\kappa\geqslant 4$, we see that the definition of $T^*$ makes sense
and $T^*>0$. Moreover, $\eta(t)\in H^4_*$ for any $t\in I_{T^{\max}}$ by Proposition \ref{pro:0401nasfxdxx} and \eqref{2020008061647}.
To obtain the existence of a global solution, we prove $T^*=\infty$ by contradiction.

Let us assume $T^*<\infty$. Then, for any given $T^{**}\in I_{T^*}$, it holds that
\begin{equation}
\label{201911262sadf202}
\sup_{0\leqslant t\leqslant T^{**}}\|(\eta,u,m\partial_2\eta\|_{4}   \leqslant K_\kappa.
\end{equation}
Thanks to \eqref{202012212151} and \eqref{201911262sadf202}, the condition
``$\max\{K_\kappa^{1/2},K^2_\kappa\}/m\leqslant  \delta_1^\kappa$'' and the fact
  \begin{align}
  \sup_{t\in \overline{I_T}} \|f\|_0=\|f\|_{L^\infty(I_T,L^2)}\;\;\mbox{ for any }f\in C(\overline{I_T}, L^2_{\mm{weak}}),\label{20212232155}
\end{align}
 we can verify that the solution $(\eta,u )$ indeed satisfies the stability estimate \eqref{202949}
 by a standard regularity method. More precisely, we have
$$\sup_{0\leqslant t\leqslant T^*}\|(\eta,u,m\partial_2\eta)(t)\|_{4}  \leqslant K_\kappa/\sqrt{c_4^\kappa}.$$

Now, we take $(\eta(T^{**}),u(T^{**}))$ as an initial data. Noting that
 $$\|(u, \partial_2\eta)(T^{**})\|_{4} \leqslant  {B}^\kappa :=2\max\{{K}_\kappa/\sqrt{c_4^\kappa},{K}_\kappa/m\sqrt{c_4^\kappa}\}
\;\;\mbox{ and }\;\; \|\nabla \eta(T^{**})\|_{3}\leqslant \delta_2^\kappa ,$$
we apply Proposition \ref{pro:0401nxdxx} to see that there exists a unique local classical solution, denoted by $(\eta^*, u^*$, $q^*)$,
to the initial-value problem \eqref{01dsaf16asdfsfdsaasf}--\eqref{01dsafsfda16asdfasfsaf} with $(\eta(T^{**}),u(T^{**}))$
in place of $(\eta^0,u^0)$. Moreover,
\begin{align}
 \sup_{t\in [T^{**},T]}\| (\eta^*,u^*,m \partial_2 \eta^* )\|_{4} \leqslant K_\kappa \mbox{ and } \sup_{t\in [T^{**},T]}
 \|\nabla \eta^*\|_{3}\leqslant 2\delta_2^\kappa ,
\end{align}
where the local existence time $T>T^{**}$ depends possibly on $B^\kappa$, $\kappa$, $m$ and $\delta_2^\kappa$.

In view of
the uniqueness in Proposition \ref{pro:0401nsadfaxsdfafdsaddfdxx} and the fact that $T^{\max}$ is the maximal existence time,
we immediately see that $T^{\max}>T^*+T/2$ and
$ \sup_{t\in  {[0,T^*+T/2]}}\| (\eta,  u,  m \partial_2 \eta)\|_{4} \leqslant K_\kappa$.
This contradicts with the definition of $T^*$. Hence, $T^*=\infty$ and thus $T^{\max}=\infty$.

In addition, we also verify that the global solution $(\eta,u)$ indeed enjoys the time-decay estimates  \eqref{202008071638} and \eqref{202911sfa9fasdfd52245} by the \emph{a priori} estimates \eqref{20200safda071638} and \eqref{202911sfa9fsafasfasdfd52245}.
Finally, the uniqueness of the global solutions is obvious due to the uniqueness of the local solutions in Proposition \ref{pro:0401nsadfaxsdfafdsaddfdxx} and the fact $\sup_{t\geqslant 0}\|\nabla \eta\|_{3}\leqslant 2 \delta_2^\kappa$.
This completes the proof of the existence and uniqueness of global time-decay solutions stated in Theorem \ref{2019043019sdfa48}.

\subsection{ {Stability around the solutions of the linearized problem}}\label{20asdfas121}
We now turn to the proof of stability around the solutions of (\ref{01dsaf16asdfsfsdafafdsaasf}). Thanks to the regularity of
$(\eta^0,u^0)$ in Theorem \ref{2019043019sdfa48}, there exists a unique solution $(\eta^{\mm{r}},u^{\mm{r}},Q_1,Q_2)$
satisfying \eqref{203212051247} and \eqref{2032120512471}. Furthermore, $(\eta^{\mm{r}},u^{\mm{r}})$ satisfies \eqref{202011032122}
and \eqref{20asdff2011032123}.

Let $\tilde{\eta}^0=\eta^0+\eta^{\mm{r}}$ and $\tilde{u}^0=u^0+u^{\mm{r}}$. Then it is easy to see that $(\tilde{\eta}^0,\tilde{u}^0)$
is in $(H^5_2\cap\underline{H}^4_\sigma)\times\underline{H}^4_\sigma$ and enjoys the odevity conditions as $(\eta^0,u^0)$ does.
Hence, there is a unique global solution $(\eta^{\mm{L}},u^{\mm{L}})\in\mathfrak{C}^0(\overline{I_T},{H}^5_2 )\times\underline{\mathfrak{U}}_\infty$
to the problem \eqref{01dsaf16asdfsfsdafafdsaasf}. Furthermore, the solution enjoys the odevity conditions as $(\eta,u)$ does.

Let $(\eta^{\mm{d}},u^{\mm{d}})=(\eta-\eta^{\mm{L}},u-u^{\mm{L}})$, then the error function $(\eta^{\mm{d}},u^{\mm{d}})$
satisfies
\begin{equation}\label{01dsfdsaasfxx}
\begin{cases}
\eta_t^{\mm{d}}=u^{\mm{d}}, \\[1mm]
u_t^{\mm{d}}+\nabla_{\mathcal{A}} q+\kappa  u^{\mm{d}} = m^2 \partial_2^2\eta^{\mm{d}}, \\[1mm]
\div u^{\mm{d}}=-\mathrm{div}_{\tilde{\mathcal{A}}} {u},\\
(\eta^{\mm{d}},u^{\mm{d}})|_{t=0}=-(\eta^{\mm{r}},u^{\mm{r}}).
\end{cases}
\end{equation}
Moreover, $(u^{\mm{d}})_{\mathbb{T}^2}=(\eta^{\mm{d}})_{\mathbb{T}^2}=0$, $\mm{div}\eta^{\mm{d}}=\mm{div}\eta$, and $(\eta^{\mm{d}},u^{\mm{d}})$
also satisfies the odevity conditions as $(\eta,u)$ does.

Employing the arguments used for \eqref{211qwessebsaddfasfdadaiseqinM0846dfgssgsd} and \eqref{211qwessgsd}, together with a standard regularity method, we deduce from \eqref{01dsfdsaasfxx} that for a.e. $t>0$,
\begin{align}
& \frac{\mm{d}}{\mm{d}t}e^{ c_3^\kappa\min\{1,m\}  t} \mathcal{E}^\kappa_{\mm{d}} + c^\kappa e^{ c_3^\kappa\min\{1,m\} t}
(\|  u^{\mm{d}}\|_4^2+\|m\eta^{\mm{d}}   \|_{5,2}^2 ) \nonumber \\
&\quad \lesssim_\kappa   e^{ c_3^\kappa\min\{1,m\} t} \| (\eta ,u , m \partial_2\eta)\|_4^2 ( \mathcal{E}^\kappa_{\mm{d}}
+\|\nabla \eta\|_3^2\| (u,m\partial_2\eta)\|_4^2),  \label{211qwesdfgssgsd} \\
& \frac{\mm{d}}{\mm{d}t} (\langle t\rangle^i\mathcal{E}^\kappa_{\mm{d}}) + c^\kappa\langle t\rangle^i(\|  u^{\mm{d}}\|_4^2
+\|m\eta^{\mm{d}} \|_{5,2}^2 ) \nonumber \\
& \quad \lesssim_\kappa \langle t\rangle^{i-1} \mathcal{E}^\kappa_{\mm{d}}+ \langle t\rangle^i \| (\eta,u, m \partial_2\eta)\|_4^2
(\mathcal{E}^\kappa_{\mm{d}}  +\|\nabla \eta\|_3^2\| (u,m\partial_2\eta^2)\|_4^2),\quad i=1,2,
\label{211dsfgqwessgsd}
\end{align}
 where $\mathcal{E}^\kappa_{\mm{d}}$ is defined by \eqref{202012051924} with $(u^{\mm{d}},\eta^{\mm{d}})$ in place of $(u,\eta)$, in $W^{1,\infty}(\mathbb{R}^+)$, and satisfies that for any $t\geqslant 0$,
\begin{align}
 \|( \eta^{\mm{d}},u^{\mm{d}})\|_4^2+\|m\eta^{\mm{d}} \|_{5,2}^2  \lesssim_\kappa \mathcal{E}^\kappa_{\mm{d}} +\|\nabla \eta\|_3^2\| (u,m\partial_2\eta)\|_4^2,\quad  \mathcal{E}^\kappa_{\mm{d}}
\lesssim_\kappa \|( \eta^{\mm{d}},u^{\mm{d}},m\partial_2\eta^{\mm{d}}) \|_{4}^2.
\label{202012281}
\end{align}

In view of  \eqref{20201221asfdasf2151}, one sees that there is a $T^0>0$, such that for any $t<T^0$,
\begin{align}
 \sup_{\tau\in \overline{I_{t}}}\mathcal{E}^\kappa_{\mm{d}}  \lesssim_\kappa &(1+t) e^{ t  \|(\eta^0,u^0,m \partial_2 \eta^0 )\|_{4}}
 \|( {\eta}^{\mm{r}},u^{\mm{r}},m\partial_2 \eta^{\mm{r}})\|_{4}^2 \nonumber \\
& + t(\|( u^0,m \partial_2 \eta^0, \tilde{u}^0,m \partial_2 \tilde{\eta}^0 )\|_{4}^3 +
 t \|( u^0,m \partial_2 \eta^0 )\|_{4}^2  \|( u^0,m \partial_2 \eta^0, \tilde{u}^0,m \partial_2 \tilde{\eta}^0 )\|_{4}^2) ,
\nonumber
\end{align}
which, together with \eqref{20asdff2011032123}, yields that
\begin{align}
 \limsup_{t\to 0}\sup_{\tau\in \overline{I_{t}}}\mathcal{E}^\kappa_{\mm{d}}  \lesssim_\kappa &
 \|({\eta}^{\mm{r}},u^{\mm{r}},m\partial_2{\eta}^{\mm{r}})\|_{4}^2 \lesssim_\kappa
\|\partial_2\eta^0 \|_{4}^2\|(\eta^0, u^0,m\partial_2\eta^0)\|_{4}^2 .
\label{20220120852000}
\end{align}
Thus, if we apply Gronwall's lemma to \eqref{211qwesdfgssgsd}, and make use of \eqref{2019092sadf81832}, \eqref{202008071638},
\eqref{2020008061647}, \eqref{202012281} and \eqref{20220120852000}, we obtain
 \begin{align}
& \sup_{t\geqslant 0} (e^{  c_3^\kappa  \min\{1,m\} t}(\|(\eta^{\mm{d}}, u^{\mm{d}})\|_4^2+\| m\eta^{\mm{d}}\|_{5,2}^2))
+\int_0^t e^{ c_3^\kappa  \min\{1,m\} \tau }(\| u^{\mm{d}}\|_4^2 +\| m \eta^{\mm{d}}\|_{5,2}^2) \mm{d}\tau  \nonumber \\
&\qquad \lesssim_\kappa (1+\mathfrak{K}^4 ) \|(\eta^0, u^0,m\partial_2\eta^0)\|_{4}^4 e^{ c_\kappa \mathfrak{K}^2 } m^{-2}, \nonumber
\end{align}
which gives \eqref{2020212052001}.
Finally, following the same process as in the derivation of \eqref{202911sfa9fsafasfasdfd52245} with necessary modifications
in arguments, we obtain \eqref{202012052058}. This completes the proof of Theorem \ref{2019043019sdfa48}
\footnote{{It is easy to see from the proof of \eqref{2020212052001} and \eqref{202012052058} that
due to the special structure of energy inequalities of $(\eta^{\mm{d}},u^{\mm{d}})$, we can use \eqref{20220120852000} and
$\|m\nabla \eta\|_3^2\leqslant C^\kappa$ to obtain the convergence rate $m^{-2}$. For the viscous fluid case, however,
we mainly exploit $\|m\partial_2\eta\|_2\leqslant C$ to get the convergence rate $m^{-1}$. This consequently yields some nonlinear terms
of third-order derivatives, such as $\|\partial_2\eta\|_2\|\nabla \eta \|_2\|\nabla q\|_1$, on the right-hand side of the energy inequalities
for $(\eta^{\mm{d}},u^{\mm{d}})$ which can only provide the convergence rate $m^{-1}$. } }.

\section{Local well-posedness}\label{201912062021}

This section is devoted to the proof of the local well-posedness results in Propositions \ref{pro:0401nxdxx} and
\ref{pro:0401nsadfaxsdfafdsaddfdxx}, and is organized as follows.

First we establish the existence of both strong and classical solutions to the following linear initial-value problem
in Section \ref{2020010813027}:
 \begin{equation}
\label{201912060857}
  \begin{cases}
 u_t+\nabla_{\ml{A}}q -\nu \Delta_{\ml{A}}u+ \kappa u = f,    \\
\mm{div}_{\mathcal{A}}u=0 ,\\
 u|_{t=0}=u^0\mbox{ in }\mathbb{T}^2,
  \end{cases}
  \end{equation}
where $\nu>0$, $\kappa\geqslant 0$, $(\eta^0,u^0,w)$ are given,
\begin{align}
\ml{A}=(\nabla \zeta)^{-\mm{T}}\;\;\mbox{ and }\;\; \zeta=\int_0^t w\mm{d}y+\eta^0+y,
\label{202009151926}
 \end{align} see Propositions \ref{qwepro:0sadfa401nxdxx} and \ref{pro:04asfda01nxdxx} for the details.
Then we give the proof of Proposition \ref{pro:0401nxdxx} based on Proposition \ref{qwepro:0sadfa401nxdxx} by a standard iterative method
in Section \ref{202009111946}. Similarly, we also give the local existence of a unique classical solution to the following initial-value
problem:
 \begin{equation}\label{01dsaf16asdfsdasafasf}
                              \begin{cases}
\eta_t=u , \\[1mm]
 u_t+\nabla_{\ml{A}}q- \nu \Delta_{\ml{A}}u+\kappa u=  f, \\[1mm]
\div_\ml{A}u=0,\\
 (\eta,u)|_{t=0}=(\eta^0,u^0)\mbox{ in }\mathbb{T}^2,
\end{cases}
\end{equation}
see Proposition \ref{pro:04asfda01nxdxx} in Section \ref{202009111946} for the details.
Finally,  thanks to Proposition \ref{pro:04asfda01nxdxx}, we complete the proof of Proposition \ref{pro:0401nsadfaxsdfafdsaddfdxx}
by a standard method of vanishing viscosity limit in Section \ref{202008062143}.

Finally, we introduce new notations appearing in this section.
\begin{align}
& H^{-1}=\mbox{the dual space of }H,\
  H^{-1}_\sigma=\mbox{the dual space of }H^1_\sigma,\nonumber \\
& { <\cdot,\cdot>_{X^{-1},X} \mbox{ denotes dual product}, \
\mathcal{U}_T^2:=\{u\in\mathcal{U}_{2,T}~|~ u_{tt}\in L^2(I_T, H ) \}, } \nonumber \\
&  { \mathcal{G}_T:=\{f\in C^0(\overline{I_T}, L^2)~|~ (f,f_t)\in L^2(I_T,H \times H^{-1})\}, } \nonumber \\
& { \|v\|_{\mathcal{U}_{1,T}}:=  }
                            \sqrt{\sum_{0\leqslant j\leqslant 1}\left( \|\partial_t^jv\|_{C^0(\overline{I_T}, H^{2(1-j)})}^2+\|\partial_t^jv\|_{L^2(I_T,H^{2(1-j)+1} )}^2 \right)}, \nonumber  \\
 &  { \|v\|_{\mathcal{U}_{T}^2}:= } \sqrt{\sum_{0\leqslant j\leqslant 2} \|\partial_t^jv\|_{C^0(\overline{I_T}, H^{2(2-j)})}^2
 +  \sum_{0\leqslant j\leqslant 2} \|\partial_t^jv\|_{L^2(I_T,H^{2(2-j)+1} )}^2},    \nonumber\\
                          &  A\lesssim_{\mm{L}}B\mbox{ means }A\leqslant c^{\mm{L}} B, \nonumber   \end{align}
 where $\mathcal{U}_{i,T}$ for $i=1$, $2$ is defined by \eqref{202102111149}, and $c^{\mm{L}}$ denotes a generic positive constant depending on $\nu$, $\kappa$ and $m$, and may vary from one place to another (if not stated explicitly).

\subsection{Unique solvability of the linear initial-value problem \eqref{201912060857}}\label{2020010813027}
 {This section is devoted to establishing the existence and uniqueness of both strong and classical solutions to
\eqref{201912060857}. We start with the existence of a unique strong solution. }
\begin{pro}
\label{qwepro:0sadfa401nxdxx}
Let $B_1>0 $, $\delta>0$, $(\eta^0,u^0)\in H^3 \cap H^2$, $\mathcal{A}^0=(\nabla\eta^0+I)^{-\mm{T}}$, $w\in\mathcal{U}_{1,T}$,
$f\in \mathcal{G}_T$, $\mathcal{A}$ and $\zeta$ be defined by \eqref{202009151926}, $\eta=\zeta -y$ and
\begin{align}  {T:=\min\{1,(\delta /B_1)^4\} }.  \label{201912061028}
\end{align}
Assume that
\begin{align}
& \|\nabla \eta^0\|_{2,2} \leqslant \delta ,\ \mm{div}_{\mathcal{A}^0}u^0=0,\   w|_{t=0}=u^0,  \nonumber \\
& \sqrt{\|\nabla w\|_{C^0(\overline{I_T},H)}^2 +\| \nabla w\|_{L^2(I_T,H^2)}^2
+\|\nabla w_t\|_{L^2(I_T,L^2)}^2}\leqslant B_1,\label{2019122821231}
\end{align}
then there is a sufficiently small constant $\delta^{\mm{L}}_1 \in (0,1]$ independent of $m$ and $\nu$, such that for any
$\delta\leqslant\delta^{\mm{L}}_1$, there exists a unique local strong solution
 $(u,q)\in\mathcal{U}_{1,T}\times(C^0(\overline{I_T},\underline{H})\cap L^2(I_T,H^2) )$ to the initial-value  problem \eqref{201912060857}.
Moreover, the solution enjoys the following estimate:
\begin{align}
& 1\leqslant 2\det\zeta\leqslant {3},\quad \| \nabla \eta \|_{2,2} \leqslant  2\delta,\nonumber \\
& \|u\|_{\mathcal{U}_{1,T}}+  \|q\|_{C^0(\overline{I_T}, H)}+\|q\|_{L^2(I_T,H^2 )}  \lesssim_{\mm{L}} \sqrt{\mathfrak{B}_1(u^0,f )},
\label{202011102145}
\end{align}
where
\begin{align}
&\mathfrak{B}_1(u^0,f) := (1+\|\nabla \eta^0\|_2^2)  \tilde{\mathfrak{B}}_1(u^0,f),\nonumber \\
&  \tilde{\mathfrak{B}}_1(u^0,f):= \|u^0\|_2^2+ \| u^0\|_2^4 +  \|f\|_{C^0(\overline{I_T},L^2)}^2+ \| f\|_{L^2(I_T,H)}^2 \nonumber \\
&\qquad \qquad \qquad +\|f_t\|_{L^2(I_I,H^{-1})}^2+\|\nabla w\|_{C^0(\overline{I_T},H)}(1+ \|  u^0 \|_1 )( \|  u^0 \|_2^2
+ \|f\|_{L^2(Q_T)}^2).\nonumber \end{align}
Moreover, if $f=\partial_2^2\eta$, then
\begin{align}
&  \sqrt{\tilde{\mathfrak{B}}_1(u^0,\partial_2^2\eta)} \lesssim_0 1+ \| u^0\|_2^2 +\sqrt{\|\nabla w\|_{C^0(\overline{I_T},H)}}(1+ \| u^0\|_2^{3/2}), \label{202010192124} \\
&\Delta_{\mathcal{A}}  {q}  = m^2\mm{div}_{\mathcal{A}}\partial_2^2\eta+ \mm{div}_{\mathcal{A}_t}u\mbox{ holds for any }t\in \overline{I_T},
\label{20209292206}\\
& \|   q \|_{C^0(\overline{I_T} ,H^2)}\lesssim_{\mm{L}} 1+ \|\nabla w\|_{C^0(\overline{I_T},H)}\left(1+ \| u^0\|_2^2 +\sqrt{\|\nabla w\|_{C^0(\overline{I_T},H)}}\left(1+ \| u^0\|_2^{3/2}\right)\right),  \label{202010170831} \\
& \| q_t\|_{L^2(I_T,H )}\lesssim_{\mm{L}} (1+ \| u^0\|_2^3+\|\nabla w\|_{C^0(\overline{I_T},H)})
(1  +\|\nabla u^0\|_0+ \|\nabla (w,w_t)\|_{L^2(I_T,H^2\times L^2)}).
\label{202010092156asfda}
\end{align}
\end{pro}
\begin{rem}
For the case $f=\partial_2^2\eta$, from \eqref{202011102145} and \eqref{202010192124} we get
\begin{align}
& \|u\|_{\mathcal{U}_{1,T}}+\|q\|_{C^0(\overline{I_T},H)}+\|q\|_{L^2(I_T,H^2 )} \leqslant c^{\mm{L}} (1+ \| u^0\|_2^3)
 + {\|\nabla w\|_{C^0(\overline{I_T},H )}}/2.  \label{20201101111005}  \end{align}
\end{rem}
\begin{pf}
We shall break up the proof into three steps.

(1) \emph{Existence of local strong solutions}

Recalling that $\|\nabla \eta^0\|_{2,2}\leqslant \delta$, {the definition \eqref{201912061028} } and
\begin{align} \label{202010011513}
\eta=\int_0^tw\mm{d}\tau+\eta^0, \end{align}
we make use of the regularity of $w$, \eqref{201912061028} and \eqref{2019122821231} to find that $\eta\in C^0(\overline{I_T},H^3)$, and
\begin{equation}
\label{201912asfdsa061030}
\| \nabla \eta(t)\|_2\leqslant \|\nabla \eta^0\|_2+\sqrt{t} \|\nabla w\|_{L^2(I_T,H^2)}\leqslant 1+\|\nabla \eta^0\|_2  ,
\end{equation}
and
\begin{equation}
\label{201912061030}
\| \nabla \eta(t)\|_{2,2} \leqslant\delta +\sqrt{t} \|\nabla w\|_{L^2(I_T,H^2_2)}\leqslant 2\delta\;\;\;\mbox{ for all }t\in \overline{I_T}.
\nonumber
\end{equation}
By \eqref{202004221412} one has
\begin{align}
\label{201sdfa912061030}
\| \nabla \eta(t)\|_{L^\infty} \lesssim_0 \| \nabla \eta(t)\|_{2,2} \lesssim_0 \delta\;\;\;\mbox{ for any }t\in \overline{I_T}.
\end{align}
Thanks to the estimate \eqref{201sdfa912061030}, we have for sufficiently small $\delta$ that $1\leqslant 2J\leqslant {3}$,
where and in what follows, $J:=\det\nabla \zeta$. Therefore, $\mathcal{A}$ makes sense and is given by the following formula:
\begin{align}
 {\mathcal{A}}=J^{-1}\left(\begin{array}{cc}
          \partial_2\eta_2+1&  - \partial_1\eta_2\\
                - \partial_2\eta_1&   \partial_1\eta_1+1
                 \end{array}\right). \nonumber
\end{align}

\emph{We remark that the smallness of $\delta$ (independent of $\nu$) will be often used in the derivation of some estimates
and conclusions later, and we shall omit to mention it for the sake of simplicity.}

Inspired by the proof in \cite[Theorem 4.3]{GYTILW1}, we next solve the linear problem \eqref{201912060857} by applying the Galerkin method.
Let $\{\varphi^i\}_{i=1}^\infty$ be a countable orthogonal basis in $H^2_\sigma$.
For each $i\geqslant 1$ we define $\psi^i=\psi^i(t):=\nabla\zeta\varphi^i$. Let $\mathfrak{H}(t)= \{v\in H^2~|~\mm{div}_{\mathcal{A}}v=0\}$.
Then $\psi^i(t)\in \mathfrak{H}(t)$ and $\{\psi^i(t)\}_{i=1}^\infty$ is a basis of $\mathfrak{H}(t)$ for each $t\in \overline{I_T}$. Moreover,
\begin{equation}
\label{201912102051}
\psi_t^i =R \psi^i,
\end{equation}
where $R:=\nabla w\mathcal{A}^{\mm{T}}$.

For any integer $n\geqslant 1$, we define the finite-dimensional space
$\mathfrak{H}^n(t):=\mm{span}\{\psi^1,\ldots, \psi ^n\}\subset \mathfrak{H}(t)$,
and write $\mathcal{P}^n(t):\mathfrak{H}(t)\to \mathfrak{H}^n(t) $ for the $\mathfrak{H}$ orthogonal projection
onto $ \mathfrak{H}^n(t)$. Clearly, for each $v\in \mathfrak{H}(t)$,
$\mathcal{P}^n(t) v \rightarrow v $ as $n\rightarrow\infty$ and $\|\mathcal{P}^n(t) v \|_2\leqslant \|v\|_2$.

Now, we define an approximate solution
$$u^n(t)=a_j^n(t) \psi^j\;\;\;\mbox{ with }a_j^n:\ \overline{I_T}\rightarrow\mathbb{R}\mbox{ for }j=1,\ldots, n,$$
where $n\geqslant 1$ is given.
We want to choose the coefficients $a_j^n$, so that for any $1\leqslant i\leqslant n$,
\begin{equation}\label{appweaksolux}
 \int  u ^n_t\cdot \psi^i \mm{d} y +\nu\int
\nabla_{\mathcal{A}}  u ^n:\nabla_{\mathcal{A}} \psi^i \mm{d} y  +\kappa\int
   u ^n:  \psi^i \mm{d} y =    \int  f  \cdot \psi^i \mm{d} y
\end{equation}
with initial data
$ u ^n(0)=\mathcal{P}^n u _0\in \mathfrak{H}^n$.

Let
$$  \begin{aligned}
& X=(a_i^n)_{n\times 1},\;\;\mathfrak{N}=  \left(\int f \cdot \psi^i \mm{d} y\right)_{n\times 1},\;\;
\mathfrak{C}^1=\left(  \int\psi^i\cdot \psi^j \mm{d} y\right)_{n\times n}, \\
 & \mathfrak{C}^2=\left( \int R \psi^i\cdot \psi^j \mm{d} y+\nu  \int
\nabla_{\mathcal{A}}\psi^i:\nabla_{\mathcal{A}}\psi^j\mm{d} y+\kappa\int\psi^i\cdot\psi^j \mm{d} y\right)_{n\times n}.
\end{aligned}$$
Recalling the regularity of $w$, we easily verify that
\begin{align}
\label{202001200222}
\mathfrak{C}^1 \in C^{1,1/2}(\overline{I_T}),\ \mathfrak{C}^2\in C^{0,1/2}(\overline{I_T}),\ \mathfrak{N} \in C^{0}(\overline{I_T})\mbox{ and } \mathfrak{N}_t\in L^2(I_T).
\end{align}

 Noting that $\mathfrak{C}^1$ is invertible, we can rewrite \eqref{appweaksolux} as follows.
 \begin{equation}  \label{appweaksolu}
 X_{t}+(\mathfrak{C}^1)^{-1}( \mathfrak{C}^2 X- \mathfrak{N})=0
\end{equation}
with initial data
$$ X|_{t=0}=\left(\int \mathcal{P}^n u _0 \cdot \psi^i\mm{d}y \right)_{n\times 1}  ,$$
where $(\mathfrak{C}^1)^{-1}$ denotes the inverse matrix of $\mathfrak{C}^1$.
By virtue of the well-posedness theory of ODEs (see \cite[Section 6 in Chapter II]{WWODE148}), the equation \eqref{appweaksolu} has
exactly one solution $X\in C^{1}(\overline{I_T})$. Thus, one has established the existence of the approximate solution
$u^n(t)=a_j^n(t)\psi^j$.  Next, we derive uniform-in-$n$ estimates for $u^n$.

Due to \eqref{201912061030}, we easily get from \eqref{appweaksolux} with $u^n$ in place of $\psi$ that
for sufficiently small $\delta$,
\begin{equation}   \label{201912112038}
\frac{\mm{d}}{\mm{d}t} \| u ^n \|_0^2   + c^{\mm{L}} \|  u^n\|_1^2\lesssim_{\mm{L}} \|f\|_0^2.
\end{equation}

By \eqref{201912102051},
\begin{equation}
\label{201912122222}
u^n_t- R  u^n= \dot{a}_i^n  \psi^i.
\end{equation}
Obviously, we can replace $\psi$ by $\dot{a}_i^n\psi^i$ in \eqref{appweaksolux} and use \eqref{201912122222} to deduce that
\begin{align}
& \|u ^n_t\|_0^2 +\nu\int
\nabla_{\mathcal{A}}  u ^n:\nabla_{\mathcal{A}} u ^n_t \mm{d}y +\kappa\int
  u ^n: u^n_t \mm{d}y \nonumber  \\
&   =  \int ( u_t^n+\kappa u^n)\cdot(R u^n)\mm{d}y +\nu \int \nabla_{\mathcal{{A}}}u^n:\nabla_{\mathcal{A}} (R u^n)\mm{d}y
 + \int  f   \cdot   (u_t^n-R u^n) \mm{d}y. \label{appweakssafolux}
\end{align}

Thanks to \eqref{202004221saffad412}, \eqref{20200508} and \eqref{201912061030}, one can further obtain from \eqref{appweakssafolux} that
\begin{align}
& \frac{\mm{d}}{\mm{d}t} \left(\kappa \|u^n\|_0^2+  {\nu } \|\nabla_{\mathcal{A}}  u ^n\|_0^2\right) +  \|u ^n_t\|_0^2\nonumber \\
 &\lesssim_{\mm{L}}\|\nabla u^n\|_0(\|R \|_{L^\infty}\|\nabla u^n\|_0
+\|\mathcal{A}_t\|_{L^\infty} \|\nabla u^n\|_0 )+\int | \nabla u^n| |\nabla R ||u^n|\mm{d}y \nonumber \\
&\quad + \|(f, R u^n) \|_{0}^2+\|R\|_{L^\infty} \|  u^n \|_{0}^2 \nonumber \\
&  \lesssim_{\mm{L}}   \|\nabla w\|_{2,2}  (\|u^n\|_0^2+\|\nabla_{\mathcal{A}}u^n\|_0^2)
+\|\nabla w\|_{1}\|\nabla w\|_{2,2} \|u^n\|_0^2+\|f\|_0^2.  \label{appwfolux}
\end{align}

With the help of Gronwall's lemma, \eqref{201912061028} and \eqref{201912061030}, we infer from \eqref{201912112038} and \eqref{appwfolux}
  that for any $t\in\overline{I_T}$,
\begin{align}
& \| u ^n
\|_1^2+\int_0^t\|u_\tau\|_0^2\mm{d}\tau \lesssim_{\mm{L}}  \left( \|\mathcal{P}^nu^0\|_1^2+\int_0^t\|f\|_0^2\mm{d}\tau\right)
 e^{\int_0^t c (\|\nabla w\|_2+\|\nabla w\|_{1}\|\nabla w\|_{2,2})\mm{d}\tau}\nonumber \\
& \qquad \lesssim_{\mm{L}} \|u^0\|_2^2+ \|f\|_{L^2(Q_t)}^2 .
\label{201912241359}
\end{align}

Recall (see \cite[Theorem 1.67]{NASII04})
$$ \int  f(\tau)\cdot \psi^i(\tau)\mm{d}y\bigg|_{\tau=0}^{\tau=t}=\int_s^t \left(<f_\tau,\psi^i>_{H^{-1},{H}  }+ \int f\cdot \psi^i_\tau\mm{d}y\right)\mm{d}\tau,$$
which gives
\begin{align}
\frac{\mm{d}}{\mm{d}t}\int f(t)\cdot \psi^i(t)\mm{d}y=<f_\tau,\psi^i>_{H^{-1},{H}  }+ \int f\cdot \psi^i_\tau\mm{d}y\;\;\mbox{ for a.e. } t\in I_T.
\label{202011091820}
\end{align}
Hence, $\mathfrak{C}^2 _t\in L^2(I_T)$. In view of \eqref{202001200222} and \eqref{appweaksolu}, we have $\ddot{a}_j^n(t)\in L^2({I_T})$.
This means that $u_{tt}^n$ makes sense. So, with the help of \eqref{201912102051} and \eqref{202011091820},
we get from \eqref{appweaksolux} that
 \begin{align}
&  \int  u ^n_{tt}\cdot \psi^i \mm{d} y+\nu\int \nabla_{\mathcal{A}}  u ^n_t:\nabla_{\mathcal{A}} \psi^i \mm{d} y +\kappa\int
  u ^n_t\cdot \psi^i \mm{d} y \nonumber \\
& = <  f_t \cdot \psi^i>_{H^{-1},{H}  } + \int  (f- u ^n_{t}-\kappa u^n)\cdot (R \psi^i)  \mm{d} y \nonumber \\
&\quad -\nu\int (\nabla_{\mathcal{A}_t} u^n:\nabla_{\mathcal{A}} \psi^i
+ \nabla_{\mathcal{A}} u^n:(\nabla_{\mathcal{A}_t} \psi^i
+ \nabla_{\mathcal{A}}(R\psi^i ))\mm{d} y\;\;\;\mbox{  a.e. in }I_T. \label{201912102242}
\end{align}

Noting that (also see \cite[Theorem 1.67]{NASII04})
$$ \begin{aligned}
& \frac{1}{2}\|u ^n_t\|_0^2 - \int  u_{t}^n\cdot (R u^n)\mm{d}y- \left.\left(\frac{1}{2}\|u ^n_t\|_0^2
 - \int  u_{t}^n\cdot ( R u^n)\mm{d}y\right)\right|_{t=0} \\
 & = \int_0^t \left(\int  u ^n_{\tau\tau}\cdot  ( u^n_{\tau}-Ru^n) \mm{d} y
 -\int  u ^n_{\tau}\cdot  (Ru^n)_{\tau} \mm{d} y\right)\mm{d}\tau
 \end{aligned} $$
and
$$  \begin{aligned}
&\int f(\tau)\cdot (Ru^n)(\tau) \mm{d}y\bigg|_{\tau=0}^{\tau=t}  =\int_0^t\left( <f_\tau,  R u^n>_{H^{-1},{H}  }
+  \int f\cdot( R u^n )_\tau\mm{d} y\right)\mm{d}\tau ,
 \end{aligned} $$
we utilize \eqref{201912122222} and the above two identities to infer from \eqref{201912102242} with $\psi^i$ replaced
by $(u^n_t-Ru^n)$ that
\begin{align}
&\frac{1}{2 }\|u ^n_t\|_0^2  - \int   u_{t}^n\cdot (R u^n)\mm{d}y+ \int f\cdot (Ru^n) \mm{d}y
+\int_0^t (\kappa\|u_{\tau}^n\|_0^2+\nu\|\nabla_{\mathcal{A}}  u ^n_\tau\|_0^2)\mm{d}\tau\nonumber \\
& = \left. \left(\frac{1}{2}\|u ^n_t\|_0^2 - \int u_{t}^n\cdot ( Ru^n)\mm{d}y
+\int f \cdot (Ru^n) \mm{d}y\right)\right|_{t=0}+I_7,
\label{201912241249}
\end{align}
where
$$    \begin{aligned}
I_7 :=&\int_0^t\bigg( <f_{\tau},  u_{\tau}^n >_{H^{-1},{H}  }
+\int \Big( f\cdot (2 Ru^n_{\tau} +R_{\tau}u^n -R^2 u^n)+ (u_\tau^n+\kappa  u^n )\cdot (R (u_\tau^n\nonumber \\
&\quad - R u^n))\Big)\mm{d}y +\int (\kappa u_{\tau}^n(Ru^n)- u_{\tau}^n\cdot (R u^n)_{\tau})\mm{d}y
-\nu\int\big(  \nabla_{\mathcal{A}} u^n: (\nabla_{\mathcal{A}_{\tau}} (u_{\tau}^n-R u^n) \nonumber \\
&\quad +\nabla_{\mathcal{A}}(R (u_{\tau}^n -R u^n))\big) +\nabla_{\mathcal{A}_{\tau}} u^n:\nabla_{\mathcal{A}} (u_{\tau}^n-R u^n)
- \nabla_{\mathcal{A}}  u^n_\tau :\nabla_{\mathcal{A}} (R u^n)) \mm{d} y\bigg)\mm{d}\tau .
\end{aligned}
$$

Keeping in mind that
\begin{align}
  \|\nabla w\|_0\leqslant \int_0^t\|\nabla w_{\tau}\|_0\mm{d}\tau+\|\nabla u^0\|_0, \label{202009121026}
\end{align}
we get from \eqref{201912241249} that
 \begin{align}
& \|u ^n_t\|_0^2  +\int_0^t(\kappa\|u_\tau\|_0^2+ \nu\| \nabla u ^n_{\tau}\|_0^2)\mm{d}\tau\nonumber \\
& \lesssim_{\mm{L}}\|\nabla w\|_{0}\|\nabla w\|_1\|u^n\|_1^2+\|u^0\|_2^4 +\|f\|_{C^0(\overline{I_T},L^2)}^2+\|u^n_t|_{t=0}\|_0^2+ I_9 \nonumber\\
& \lesssim_{\mm{L}}  \|\nabla w\|_1(1+ \|u^0 \|_1 )\left(\|u^0 \|_2^2+\|f\|_{L^2(Q_t)}^2\right)+ \|u^0\|_2^4+ \|f\|_{C^0(\overline{I_T},L^2)}^2
+\|u ^n_t|_{t=0}\|_0^2+ I_9,
\label{201912242022}
\end{align}
where we have used \eqref{20200508} in the first inequality, \eqref{201912061028}, \eqref{201912241359} and \eqref{202009121026}
in the second inequality. Below, we shall bound the the last two terms in \eqref{201912242022}.

Replacing $\psi^i$ by $(u_t^n-Ru^n)$ in \eqref{appweaksolux}, one sees that
\begin{align}
 \| u^n_t \|_0^2 = &  \int f \cdot (u ^n_t-R u^n) \mm{d} y
 +\nu\int \Delta_{\mathcal{A}}u^n:(u_t^n-R u^n)\mm{d}y\nonumber\\
&+ \int u_t^n\cdot (R u^n)\mm{d}y-\kappa\int u^n\cdot (u_t^n-Ru^n)\mm{d}y ,
\end{align}
which implies
$$ \|u^n_t \|_0^2\lesssim_{\mm{L}}\| f \|_0^2+\|  u^n\|_2^2 +\|\nabla w\|_1^2\|u^n\|_2^2,\quad\forall\, t\in [0,T). $$
In particular, \begin{equation}
\label{201912112039}
\| u ^n_t|_{t=0} \|_0^2\lesssim_{\mm{L}}  \|u^0\|_2^2+ \|u^0\|_2^4+\|f^0\|_0^2 .
\end{equation}

 Thus, the last term on the right-hand side of \eqref{201912241249} can be estimated as follows, using
the first inequality in \eqref{20200508} and \eqref{201912061028}.
\begin{align}
 I_7\lesssim_{\mm{L}} & \int_0^t\Bigg(\|f_{\tau}\|_{-1} \|u^n_{\tau}\|_{1}
+\|u^n\|_1\left( \| f\|_0\|\nabla w\|_{L^\infty}\|\nabla w\|_1+ \sqrt{\| f\|_0\|f\|_1}\|\nabla w_{\tau}\|_0\right) \nonumber \\
&\quad+ \|f\|_1\| \nabla w\|_1\|u_{\tau}^n\|_0+ \| u^n  \|_1^2 \|\nabla w\|_{L^\infty}\left( \|\nabla w\|_{L^\infty}
+\sqrt{\|\nabla  w\|_1\|\nabla  w\|_{2,2}}  \right) \nonumber \\
 &\quad+ \| u^n_\tau \|_0(\|\nabla w\|_2\| u^n_\tau \|_0 + (1+\|\nabla w\|_{L^\infty})\|\nabla w\|_1 \|u^n \|_1)
+ \|\nabla w_\tau\|_0 \|u^n \|_1\sqrt{\| u^n_\tau \|_0 \| u^n_\tau \|_1}\nonumber \\
&\quad  +\|u^n \|_1\| u^n_\tau \|_1\left(   \|\nabla w\|_{L^\infty}+  \sqrt{\|\nabla  w\|_1\|\nabla  w\|_{2,2}}  \right)
\Bigg)\mm{d}\tau \nonumber \\
 \lesssim_{\mm{L}}  & \|u^0\|_2^2+\|f\|_{C^0(\overline{I_t},L^2)}^2 +\|  f\|_{L^2(I_t,H)}^2
+  \int_0^t \|f_{\tau}\|_{H^{-1}}  \| u_{\tau}^n\|_1 \mm{d}\tau\nonumber \\
& + \left(  \|u^0\|_2   + \| f\|_{L^2(I_t, H )} \right)( \|u_t^n\|_{{C^0(\overline{I_t},L^2)} } + \| \nabla u_\tau^n\|_{L^2(I_t,L^2)} )
 +\delta  \|u_t^n\|_{{C^0(\overline{I_t},L^2)} }^2. \label{202001221640}
\end{align}

Substituting \eqref{201912112039} and \eqref{202001221640} into \eqref{201912242022}, and applying Young's inequality, we arrive at
\begin{align}
&  \| u ^n_t \|_{{C^0(\overline{I_T},L^2)} }^2 +   \| u ^n_{t}\|_{L^2(I_T,H)}^2  \lesssim_{\mm{L}}  \tilde{\mathfrak{B}}_1(u^0,f).\label{appwasdfasfsadfeaksolux}
\end{align}
Summing up \eqref{201912241359} and \eqref{appwasdfasfsadfeaksolux}, we conclude
\begin{align}
 \| (u^n ,u ^n_t) \|_{C^0(\overline{I_T},H\times L^2)}^2+ \|(u ^n,u _t^n)\|_{L^2(I_T,H)}^2 \lesssim_{\mm{L}}    \tilde{\mathfrak{B}}_1(u^0,f). \label{sumestimeat}
\end{align}

In view of \eqref{sumestimeat}, the Banach--Alaoglu and Arzel\`{a}--Ascoli theorems, up to the extraction of a subsequence
(still labelled by $u^n$), we have, as $n\to\infty$, that
$$\begin{aligned}
&(u^n,u^n_t)\rightarrow (u,u_t)\;\; \mbox{ weakly-$*$ in }L^\infty(I_T,H\times L^2 ) ,\\
&(u^n,u^n_t)\rightarrow (u,u_t)\;\; \mbox{ weakly  in }L^2(I_T,H\times H),\\
& u_n\rightarrow u \;\;\mbox{ strongly in }C^0(\overline{I_T},L^2),\\
&\mm{div}_{\mathcal{A}} u =0 \;\; \mbox{ a.e. in }\Omega_T  \mbox{ and }\;u(0)=u_0,
\end{aligned}$$
where $u$ and $u_t$ are measurable functions defined on $\Omega_T$.
Moreover,
\begin{align}
&  \|(u,u_t) \|_{L^\infty(I_T,L^2)} + \|(u,u _t)\|_{L^2(I_T,H)} \lesssim_{\mm{L}}\sqrt{\tilde{\mathfrak{B}}_1(u^0,f)} .  \label{2019121asdfa12157}
\end{align}
Therefore, we can take to the limit in \eqref{appweaksolux} as $n\to\infty$, and obtain
\begin{equation}\label{0507n1}
  \int u _t\cdot\zeta \mm{d} y +\nu\int\nabla_{\mathcal{A}} u:\nabla_{\mathcal{A}} \zeta\mm{d} y
+\kappa\int u\cdot\zeta\mm{d} y =\int f\cdot\zeta\mm{d} y\;\;\;\mbox{ a.e. in }I_T,\quad\forall\,\zeta\in\mathfrak{H}.
\end{equation}

Now, we begin to show spatial regularity of $u$. Let us further assume that $\delta $ is so small that
$\eta$ satisfies \eqref{2312018031adsadfa21601xx} by virtue of Proposition \ref{pro:0401nasfxdxx}.
Denoting $F:=f-\kappa u-u_t$, $\tilde{F}:=F(\zeta^{-1},t)$ and $\tilde{J}:=J(\zeta^{-1},t)$, we see
that $\tilde{F}$ has the same regularity as that of $F$, i.e.,
\begin{align}
 \label{20218181407} \| \tilde{F} \|_{L^\infty(I_T,L^2)}+  \|  \tilde{F} \|_{L^2(I_T,H)}  <\infty.\end{align}
Moreover, \begin{align}
\int F(y,t)\mm{d}y= \int \tilde{F}\tilde{J}^{-1}\mm{d}x . \nonumber  \end{align}

Applying the regularity theory of the Stokes problem, we see that there is a unique strong solution
 $\alpha \in  {L^\infty( {I_T},\underline{H}^2)} \cap  L^2(I_T,{H}^3)$ with a unique associated function
 $p \in {L^\infty( {I_T},\underline{H})} \cap  L^2(I_T,{H}^2)$, such that
\begin{equation}
\label{appstokes}
\begin{cases}
\nabla  p -\nu \Delta \alpha =\tilde{F}  , \\
\mm{div}\,v =0.
\end{cases}
\end{equation}
 Let $\varpi=\alpha(\zeta,t)$ and $q=p(\zeta,t)- (p(\zeta,t))_{\mathbb{T}^2} $, then $(\varpi,q)\in ({L^\infty( {I_T},H^2)} \cap  L^2(I_T,{H}^3))\times ({L^\infty( {I_T},\underline{H})} \cap  L^2(I_T,{H}^2)) $  satisfies the following system:
\begin{equation}
\label{appsdfstokes}
\begin{cases}
\nabla_{\mathcal{A}}  {q} -\nu \Delta_{\mathcal{A}}\varpi= {F},  \\
\mm{div}_{\mathcal{A}} w =0\;\;\quad\mbox{ for a.e. }t\in I_T.
\end{cases}
\end{equation}
By a density argument, the identity \eqref{0507n1} also holds for $\zeta \in {H}$ with $\mm{div}_{\mathcal{A}}\zeta=0$.
This fact, together with \eqref{appsdfstokes}, implies $\nabla u=\nabla\varpi$.

 {Following the derivation of \eqref{2017020614181721} with slight modification,}
we can get from \eqref{appsdfstokes} that for a.e. $t\in I_T$,
\begin{align}
&\|q\|_{i+1} \lesssim_0 \|\nabla q\|_i \lesssim_0  \|(f,  u_t)\|_i \quad \hbox{ for }i=0,\ 1.
 \label{2020119092203}
\end{align}
Similarly to the derivation of \eqref{2020119092203} with $i=1$, we can  derive from \eqref{appsdfstokes} that
\begin{align}
&\| \varpi\|_2 \lesssim_0
  \|(u,  u_t,\nabla q,f)\|_0,
\label{2020119safas092203}\\
 &\|   \varpi\|_{3,2}  \lesssim_0
  \|(u,  u_t,\nabla q,f)\|_{1}.
\label{20201fas092203} \\
&\| \varpi\|_3  \lesssim_0
  \|(u,  u_t,\nabla q,f)\|_{1}+\|\nabla \eta\|_2\| \varpi\|_{3,2}.
\label{20201fassfdads092203}
\end{align}

So, it follows from \eqref{201912asfdsa061030}, \eqref{2019121asdfa12157} and  \eqref{2020119092203}--\eqref{20201fassfdads092203} that
\begin{align}
  \| (  u,u_t,q)\|_{L^\infty(I_T,H^2\times L^2\times  H)} \lesssim_0
  \|(u,  u_t,f)\|_{L^\infty(I_T,L^2)}\lesssim_{\mm{L}} \sqrt{\tilde{\mathfrak{B}}_1(u^0,f)}
\label{201912221611x}
\end{align}
and
\begin{align}
  \|(  u,u_t, q)\|_{L^2(I_T,H^2\times H \times H^2)} \lesssim_0
 (1+\|\nabla \eta^0\|_2)    \|(u,  u_t,f)\|_{L^2(I_T,H)} \lesssim_{\mm{L}} \sqrt{\mathfrak{B}_1(u^0,f)}.
\label{2019asdfasf12221611x}
\end{align}
Combining \eqref{201912221611x} with \eqref{2019asdfasf12221611x}, one obtains
\begin{align}
 \|  (u,u_t,q) \|_{L^\infty(I_T,H^{2 }\times L^2\times H)}  + \|  (u, u_t,q) \|_{L^2(I_T,H^3\times H\times H^2)}
\lesssim_{\mm{L}}  \sqrt{\mathfrak{B}_1(u^0,f) }.\label{201912221611}
 \end{align}
This completes the existence of local strong solutions. Moreover, a strong solution, which enjoys the regularity of $(\eta,u)$
constructed above, is obviously unique.

(2) \emph{Strong continuity of $(u_t,u,q)$ on $\overline{I}_T$ with values in $L^2\times H^2\times H$}.

For any given $\varphi\in {H}  $, let $\psi =\varphi(\zeta(y,t)) $. Noting $J_t=J\mm{div}_{\mathcal{A}}w$, we can derive
from \eqref{appsdfstokes} with $v$ in place of $w$ that for any $\phi\in C_0^\infty(I_T)$,
$$ -\int_0^t \phi_t\int  u \cdot \psi J\mm{d} y \mm{d}\tau
 = \int_0^t\phi \int ( f +\nu \Delta_{\mathcal{A}}  u -  \kappa   u - \nabla_{\mathcal{A}}  q)
\psi  J \mm{d} y \mm{d}\tau- \int_0^t \phi\int w\cdot \nabla_{\mathcal{A}}u\cdot\psi J \mm{d} y \mm{d}\tau .$$

Let ${v}=u(\zeta^{-1}(x,t),t)$, $\tilde{w}=w(\zeta^{-1}(x,t),t)$ and $ g=f(\zeta^{-1}(x,t),t)$, then from
the above identity we can get that
\begin{align}
 -\int_0^t \phi_t\int  v \cdot \varphi   \mm{d} x \mm{d}\tau
& =\int_0^t\phi \int  (  f+\nu \Delta_{\mathcal{A}}  u -  \kappa   u- \nabla_{\mathcal{A}} q
+ w\cdot\nabla_{\mathcal{A}}u )\cdot  \psi   J  \mm{d} y \mm{d}\tau \nonumber \\
& = \int_0^t\phi\int (g +\nu\Delta v -\kappa v -\nabla p -\tilde{w}\cdot\nabla v)\cdot\varphi\mm{d}x\mm{d}\tau,
\label{appweakdfsafssfaolux}
\end{align}
which immediately results in
\begin{align}
&\label{appweaksssaffaolux}
 v_t= g +\nu \Delta v -  \kappa  v - \nabla  p - \tilde{w} \cdot \nabla v \in  L^\infty(I_T,L^2)
\end{align}
and
\begin{align}
&  \int  v_t \cdot \varphi  \mm{d} x =  \int  ( f+\nu
\Delta_{\mathcal{A}}  u-  \kappa   u-   \nabla_{\mathcal{A}}  q - w\cdot\nabla_{\mathcal{A}}    u  ) \cdot \psi J  \mm{d} y  .\label{appweaksadfasssfaolux}
\end{align}
We find by  \eqref{appsdfstokes}$_1$ and \eqref{appweaksssaffaolux} that
\begin{align}
u_t=(v_t+\tilde{w} \cdot \nabla v  ) |_{x=\zeta}=v_t|_{x=\zeta}+{w} \cdot \nabla_{\mathcal{A}} u .
\label{202010042131}
\end{align}

Now let us further assume $\varphi\in H^1_\sigma$, then $\mm{div}_{\mathcal{A}}\psi=0$.
Recalling $\partial_j(J\mathcal{A}_{ij} )=0$, the identity \eqref{appweaksadfasssfaolux} implies
\begin{align}
 \frac{\mm{d}}{\mm{d}t}\int  v_t \cdot \varphi  \mm{d} y =& <f_t,\psi J>_{H^{-1},H} + \int
(f(\psi J)_t- \partial_t (((\kappa   u +w\cdot\nabla_{\mathcal{A}}    u) \cdot \psi \nonumber\\
& +\nu \nabla_{\mathcal{A} }  u :\nabla_{\mathcal{A}} \psi )J))\mm{d}y =:<\chi,\varphi>_{H^{-1}_\sigma,H_\sigma}. \label{appwasfdaeakssfaolux}
\end{align}

Noting that $\|\psi\|_1\lesssim_0 \|\varphi\|_1\lesssim_0 \|\psi\|_1$ for any $t\in I_T$, and recalling
the definition of $<\chi,\varphi>_{H^{-1}_\sigma,H_\sigma}$, we have $\chi\in {L^2(I_T,H^{-1}_\sigma)}$.
Therefore, $v_{tt}=\chi$ and $v_t\in C^0(\overline{I_T},L^2)$. Consequently, the identity \eqref{202010042131} implies
$u_t\in  C^0(\overline{I_T},L^2)$. Since $u\in L^2(I_T,H^3)$ and $ u_t\in L^2(I_T, {H}  )$, $u\in C^0(\overline{I_T}, H^2)$.
Hence, $u\in \mathcal{U}_{1,T}$.
In addition, we can derive from \eqref{appsdfstokes}$_1$ that $q\in C^0(\overline{I_T},\underline{H}^1)$ for sufficiently small $\delta$.
Thanks to the strong continuity of $(u_t,u, q)$ on $\overline{I}_T$ with values in $L^2\times H^2\times H$, we immediately
get \eqref{202011102145} from \eqref{201912221611}.

(3) \emph{More regularities of $q$ under the case ``$f=\partial_2^2\eta$''.}

Obviously, \eqref{202010192124} holds for $f=\partial_2^2\eta$. Keeping in mind that $p$ is in $C^0(\overline{I_T},H)$,
we obtain from \eqref{appstokes} that
\begin{equation}\label{appstsafokes}
\int \nabla  p \cdot \nabla \varphi\mm{d}x= \int (g-\tilde{w}\cdot \nabla  v)\cdot \nabla \varphi\mm{d}x
\;\;\;\mbox{ for any }t\in \overline{I_T},
\end{equation}
which implies that $p(t)\in \underline{H}^2$ ($t\in\overline{I_T}$) satisfies
\begin{align}
\Delta  p =\mm{div}g- \nabla\tilde{w}: \nabla  v^{\mm{T}}\;\;\mbox{  for any }t\in \overline{I_T}.
\label{202010071249}
\end{align}
Then, $q(t)\in H^2$ satisfies \eqref{20209292206}. Thus, the estimate \eqref{202010170831} follows from
\eqref{202011102145}--\eqref{20209292206}.

Recalling the derivation of \eqref{2017020614181721} with $i=1$ and the fact
$\|u(\tau)|_{\tau=s}^{\tau=t}\|_1\leqslant \int_s^t\|u_\tau\|_1\mm{d}\tau$,
we see that $q\in AC^0(\overline{I_T},H)$, i.e., $ q$ is absolutely contiguous in $\overline{I_T}$ with respect to the norm $H$.
Next, we proceed to show that $q_t$ exists and enjoys the estimate \eqref{202010092156asfda}.

By virtue of the Riesz representation theorem, it is easy to show that there is a unique function
$\chi\in {L}^2(I_T,\underline{H}^1)$, such that
\begin{align}
-\int \nabla_{\mathcal{A}}\chi\cdot \nabla_{\mathcal{A}}\varsigma\mm{d}y
=\int\big(\partial_t(m^2 \mathcal{A}^{\mm{T}}\partial_2^2\eta+\mathcal{A}_t^{\mm{T}} u)
+\mathcal{A}^{\mm{T}}\nabla_{\mathcal{A}_t} q +\mathcal{A}^{\mm{T}}_t\nabla_{\mathcal{A}} q\big)\cdot\nabla\varsigma\mm{d}y,\;\;\;
\varsigma\in \underline{H}^1.
\label{202011202136}
\end{align}
Moreover, $w$ enjoys the following estimate
\begin{align}
\|\chi\|_{L^2(I_T,H)}\lesssim_{\mm{L}} & \|(u,u_t,q)\|_{C^0(\overline{I_T},H^2\times L^2\times H)}\big( 1+\|\nabla u^0\|_0\nonumber \\
&  + \|\nabla (w,w_t)\|_{L^2(I_T,H^2\times L^2)}\big) +\|\nabla w\|_{C^0(\overline{I_T},H)}.
\label{202010092156}
\end{align}

Let $t\in I_T$ and $D_h\vartheta=\big(\vartheta(y,t+h)-\vartheta(y,t)\big)/h$ where $t+h\in I_T$.
Multiplying \eqref{20209292206} by $\varsigma$ in $L^2$ and applying then $D_h$ to the resulting equation, we get
\begin{align}
-\int \nabla_{\mathcal{A}} D_hq\cdot \nabla_{\mathcal{A}}\varsigma\mm{d}y =&\int(D_h(m^2 \mathcal{A}^{\mm{T}}\partial_2^2\eta
+\mathcal{A}_t^{\mm{T}} u) + \mathcal{A}^{\mm{T}}\nabla_{D_h\mathcal{A} } q(t+h)\nonumber \\
& +D_h(\mathcal{A}^{\mm{T}} ) \nabla_{\mathcal{A}(y,t+h) } q(y,t+h)) \cdot \nabla\varsigma \mm{d}y. \label{2020112sdfa02136}
\end{align}
Subtracting \eqref{202011202136} from \eqref{2020112sdfa02136} and denoting $\varsigma =D_hq -\chi$, we have
\begin{align}
\|  D_hq-\chi \|_1\lesssim_{\mm{L}} & \|(D_h-\partial_t)(m^2 \mathcal{A}^{\mm{T}}\partial_2^2\eta+\mathcal{A}_t^{\mm{T}} u)\|_0
+ \|\mathcal{A}^{\mm{T}}\nabla_{D_h\mathcal{A} - \mathcal{A}_t} q(y,t+h)\|_0 \nonumber \\
  & +\|\mathcal{A}^{\mm{T}}\nabla_{\mathcal{A}_t } (q(t+h)-q(t))\|_0 +
\|(D_h-\partial_t)\mathcal{A}^{\mm{T}} \nabla_{\mathcal{A}(y,t+h)} q(y,t+h)\|_0 \nonumber \\
& +\| \mathcal{A}^{\mm{T}}_t(\nabla_{\mathcal{A}(y,t+h) } q(y,t+h) -\nabla_{\mathcal{A}} q)\|_0 =:\Theta(t)\;\;\mbox{ for a.e. }t\in I_T.
\nonumber
\end{align}

 Noting that the generalized derivative with respect to $t$ is automatically strong derivative, we easily see that
$\Theta(t)\to 0$ for a.e. $t\in I_T$. So, $\|(D_hq-\chi)\|_1^2\to 0$ as $h\to 0$ for a.e. $t\in I_T$.
 {This means that the strong derivative of $ q$ with respect to $t$ is equal to that of $\chi$.}
Because of $q\in AC^0(\overline{I_T},H)$, $q_t=\chi$, where $q_t$ denotes the generalized derivative of $q$.
Hence, $q_t\in L^2(I_T,H)$ satisfies the estimate \eqref{202010092156asfda} by \eqref{202010192124}, \eqref{201912221611x}
and \eqref{202010092156}. This completes the proof of Proposition \ref{qwepro:0sadfa401nxdxx}. \hfill $\Box$
\end{pf}

Now, we turn to establishing the existence and uniqueness of classical solutions
to the initial-value problem \eqref{01dsaf16asdfsdasafasf}.
\begin{pro} \label{pro:04asfda01nxdxx}
Let $\delta>0$. Under the assumptions of Proposition \ref{qwepro:0sadfa401nxdxx} with $f=\partial_2^2\eta$,
assume further that $(\eta^0,u^0)\in H^5\times H^4$, $\|\nabla\eta^0\|_{4}\leqslant\delta$ and
\begin{align}\label{202011111628}
& \sqrt{\|\nabla ( w,w_t)\|_{C^0(\overline{I_T},H^3\times H)}^2
+\|w_{tt}\|_{C^0(\overline{I_T},L^2)}^2+\|\nabla (w_t,w_{tt})\|_{L^2(I_T,H^2\times L^2)}^2}\leqslant B_1.
\end{align}
Then, there is a sufficiently small constant $\delta^{\mm{L}}_2\leqslant \delta^{\mm{L}}_1$,
such that for any $\delta\leqslant \delta^{\mm{L}}_2$, the solution $(u,q)$ constructed by Proposition \ref{qwepro:0sadfa401nxdxx}
belongs to $\mathcal{U}_T^2\times C^0(\overline{I_T},H^4)$, where $\delta^{\mm{L}}_2$ is independent of $m$ and $\nu$. Moreover,
  \begin{eqnarray}
 && \|u\|_{\mathcal{U}_{2,T}} + \| q_t  \|_{C^0(\overline{I_T},{H}^2)} \lesssim_\kappa \left(1+\|\nabla \eta^0\|_4\right)
 (1+\|u^0\|_4^6+\|u^0\|_2^4\|\nabla w_t|_{t=0} \|_0), \label{20191asfas2242055} \\
 && \|q\|_{C^0(\overline{I_T},{H}^4)} \lesssim_\kappa\|\nabla\eta^0\|_4+ \|w\|_{C^0(\overline{I_T},{H}^4)}\|u\|_{C^0(\overline{I_T},{H}^4)}, \label{20191asfas2242055xx}  \\
 && \|\nabla \eta\|_{4} \leqslant 2\delta,\;\;\;\mbox{ and }\;\;\;\|u_t|_{t=0}\|_2\lesssim 1+\|u^0\|_4. \label{20191asfasfasas2242055}
\end{eqnarray}
\end{pro}
\begin{pf}
Let $(u,q)$ be constructed in Proposition \ref{qwepro:0sadfa401nxdxx} and
 $\|\nabla \eta^0\|_{4,2}\leqslant \delta\leqslant \delta^{\mm{L}}_1$. Recalling
\begin{align}
T=\min\{1,(\delta /B_1)^4\}, \label{2011221940}
\end{align}
we see that
\begin{align}
\|\nabla \eta\|_{4 }\leqslant 2\delta\quad\mbox{ for any }t\in \overline{I_T}.
\label{202009301500}
\end{align}
 We remark here that the smallness of $\delta$ will be used in the derivation of some estimates later. From now on, we denote
  $D_t\sigma:=\sigma_t-R\sigma$.

By \eqref{20209292206} we see that $q^0\in\underline{H}^2$ satisfies
$$\Delta_{\mathcal{A}^0}{ q^0}=\mm{div}_{\mathcal{A}^0}\partial_2^2\eta^0+\mm{div}_{\mathcal{A}_t|_{t=0}}u^0.$$
Thanks to \eqref{202009301500}, we can apply a standard difference quotient method to deduce that $q^0\in \underline{H}^3$.
Moreover,
\begin{align}   \label{202010191943}
\|\nabla q^0\|_2\lesssim_0 \|\partial_2^2\eta^0\|_2+\|u^0\|_3^2 \lesssim_0 1+\|u^0\|_3^2.
\end{align}

Recalling $u_t^0+\nabla_{\ml{A}^0}q^0 -\nu \Delta_{\ml{A}^0}u^0+ \kappa u^0 =\partial_2^2\eta^0$, we have
\begin{align}
\|u_t|_{t=0}\|_2\lesssim 1+ \| \nabla q^0 \|_2+\|u^0\|_4\lesssim 1+\|u^0\|_4,
\label{2020212271533}
\end{align}
which gives
\begin{align}
\|D_t u|_{t=0}\|_2\lesssim 1 +\|u^0\|_4^2.
\label{202009301532}
\end{align}

Noting that
\begin{align}
\label{202011222006}
\|\nabla w\|_2\leqslant \int_0^t\|\nabla w_t\|_2\mm{d}\tau +\| u^0\|_3,
\end{align}
we make use of \eqref{202010192124}, \eqref{202011111628}, \eqref{2011221940} and \eqref{202011222006} to
deduce from the definition of $\mathfrak{B}_1(u^0,\partial_2^2\eta )$ that
\begin{align}
\mathfrak{B}_1(u^0,\partial_2^2\eta )\lesssim_{\mm{L}} 1+ \| u^0\|_2^4.
\label{202012091547}
\end{align}
Similarly, we can obtain by using \eqref{202011102145}, \eqref{202010170831} and \eqref{202010092156asfda} that
\begin{align}
\|u\|_{\mathcal{U}_{1,T}} + \| q \|_{C^0(\overline{I_T},H^2)} +\|q_t\|_{L^2(I_T,H)}\lesssim_{\mm{L}} 1 +\|u^0\|_3^4.
\label{202010092156safdaasfda}
\end{align}

Let $\mathcal{F}=\nu(\mm{div}_{\mathcal{A}_t}\nabla_{\mathcal{A}}u^{\mm{T}}
+\mm{div}_{\mathcal{A}}\nabla_{\mathcal{A}_t}u^{\mm{T}} +\Delta_{\mathcal{A}}(Ru)-R\Delta_{\mathcal{A}}u )+R\nabla_{\mathcal{A}} q-\nabla_{\mathcal{A}_t}q-R_tu$.
One may use \eqref{202011111628}, \eqref{2011221940}, and \eqref{202009301500}, \eqref{202009301532} and \eqref{202011222006} to verify that
\begin{align}
 {\mathfrak{B}_1(D_tu|_{t=0},D_t\partial_2^2\eta+\mathcal{F})}
\lesssim 1 & +  \|u^0\|^8_4+  \|\mathcal{F} \|^2_{C^0(\overline{I_T},L^2)} + \| \mathcal{F}  \|_{L^2(I_T,H)}^2 +\| \mathcal{F}_t  \|_{L^2(I_T,H^{-1})}^2\nonumber \\
& + (1+\|u^0\|_4^3  ) \left(1+\|u^0\|_4^4  +\|\mathcal{F}\|_{L^2(Q_T)}^2\right) .
\label{202010111353}
\end{align}

Thanks to \eqref{202010092156asfda}, \eqref{202011111628}, \eqref{2011221940}, \eqref{202011222006}  and \eqref{202010092156safdaasfda}, we have the following upper bounds
for $\mathcal{F}$ and $\mathcal{F}_t$:
\begin{align}
& \|\mathcal{F}\|_{L^2(I_T,H)}
  \lesssim_{\mm{L}} T \sup_{t\in \overline{I_T}}(\|u\|_{2}(\|\nabla w\|_{3}+\|\nabla w\|_{2}^2+\|\nabla w_t\|_{1}))  \nonumber \\
 &\qquad +  \|\nabla w \|_2 ( \|u\|_{L^2(I_T,H^3)}+\|\nabla q\|_{L^2(I_T,H)}) \lesssim_{\mm{L}}   1+ \|  u^0\|_4^4   \nonumber
\end{align}
and
\begin{align}
\! \|\mathcal{F}_t\|_{L^2(I_T,H^{-1})}\lesssim_{\mm{L}}&  T\sup_{t\in \overline{I_T}}\big( (\|\nabla w\|_3+\|\nabla w\|_2^2+\|\nabla w_t\|_1
+\|\nabla w\|_2\|\nabla w_t\|_0)(\|u\|_2+\|u_t\|_0+\|\nabla q\|_1)\big)\nonumber \\
 & +\sup_{t\in \overline{I_T}}   \|\nabla w\|_2  (\|u_t\|_{L^2(I_T,H)}+ \|\nabla q_t\|_{L^2(I_T,L^2)} )\nonumber \\
&   +
\int_0^T\Big(\sqrt{\|w_{tt}\|_0\|w_{tt}\|_1}\|u\|_2+ \sqrt{\|\nabla w_{t}\|_1 \|\nabla w_{t}\|_2 }(\|u\|_2+\|u_t\|_0)\Big)\mm{d}t
\nonumber \\
\lesssim_{\mm{L}}& 1+ \| u^0\|_4^5. \nonumber
\end{align}
Moreover,
\begin{align}
\|\mathcal{F}\|_{C^0(\overline{I_T},L^2)}\lesssim_0 \|\mathcal{F}\|_{L^2({I_T},H)}+\|\mathcal{F}\|_{L^2({I_T},H^{-1})}\lesssim_{\mm{L}}
1+ \|u^0\|_4^5. \nonumber
\end{align}
Substitution of the above three estimates into \eqref{202010111353} yields
\begin{align}
 \sqrt{\mathfrak{B}_1(D_tu|_{t=0},D_t\partial_2^2\eta+\mathcal{F})} \lesssim_{\mm{L}} 1+
\|u^0\|_4^5. \label{202012091528}
\end{align}

Now, let us consider the problem
\begin{align}
\label{20191206fda0857}
  \begin{cases}
U_t+\nabla_{\ml{A}}Q -\nu \Delta_{\ml{A}}U+ \kappa U=
D_t\partial_2^2\eta+\mathcal{F},    \\
\mm{div}_{\mathcal{A}}U=0,\\
w|_{t=0}=D_tu|_{t=0}.
  \end{cases}
\end{align}
Recalling $\mathrm{div}_{\mathcal{A}^0} (D_tu)|_{t=0}=0$, and using \eqref{202009301532}, \eqref{202011222006} and \eqref{202012091528}, we can apply Proposition \ref{qwepro:0sadfa401nxdxx}
to \eqref{20191206fda0857} to see that the initial-value problem \eqref{20191206fda0857} admits a unique strong solution
$(U,Q)\in\mathcal{U}_{1,T}\times(C^0(\overline{I_T},\underline{H})\cap L^2(I_T,H^2) )$ , which satisfies
\begin{align}
 \|U\|_{\mathcal{U}_T } + \| Q\|_{C^0(\overline{I_T},H^2)}\lesssim_{\mm{L}} 1 +\|u^0\|_4^5 .\label{20191224assdfadf2055}
\end{align}

Let $\tilde{U}=U(\zeta^{-1}(x,t),t)$ and $\tilde{U}(x,0)=(D_tu|_{y=\zeta^{-1}(x,t)})|_{t=0}$.
Similarly to \eqref{appweakdfsafssfaolux}, we can obtain from \eqref{20191206fda0857} that
\begin{align}
& -\int_0^T\phi\int  \tilde{U}\cdot \xi \mm{d}x\mm{d}t  +\int \phi\int( \nu \nabla \tilde{U}:\nabla \xi +\kappa\tilde{U} \cdot\xi)\mm{d}x\mm{d}t\nonumber \\
&\quad =\int_0^T\phi\int ((D_t\partial_2^2\eta+\mathcal{F})|_{y=\zeta^{-1}(x,t)}-\tilde{w}\cdot \nabla \tilde{U})\cdot\xi\mm{d}x\mm{d}t,\quad\;
\phi\in C_0^\infty (I_T),\; \xi\in H_\sigma^1, \nonumber
\end{align}
which implies
\begin{align}
&\frac{\mm{d}}{\mm{d}t}\int\tilde{U}\cdot\xi\mm{d}x +\int( \nu\nabla\tilde{U}:\nabla \xi +\kappa\tilde{U}\cdot\xi)\mm{d}x\nonumber \\
& \quad  =\int ((D_t\partial_2^2\eta+\mathcal{F})|_{y=\zeta^{-1}(x,t)} -\tilde{w}\cdot \nabla v)\cdot\xi\mm{d}x \;\;\; \mbox{ for a.e. }t\in I_T.
\label{2020sa854}
\end{align}

Since $u$ solves \eqref{01dsaf16asdfsdasafasf} with $f=\partial_2^2\eta$, we utilize the regularity of $u$ to find that
\begin{align}
&\frac{\mm{d}}{\mm{d}t}\int  \widetilde{Du_t}\cdot \xi \mm{d}x  +\int( \nu \nabla \widetilde{Du_t}:\nabla \xi +\kappa\widetilde{Du_t}\cdot\xi)\mm{d}x\nonumber \\
&\quad =\int ((D_t\partial_2^2\eta+\mathcal{F})|_{y=\zeta^{-1}(x,t)}-\tilde{w}\cdot \nabla v)\cdot\xi\mm{d}x\;\;\;\mbox{ for a.e. }t\in I_T,
\label{2020sa8sfda54}
\end{align}
where $\widetilde{Du_t}:= Du_t|_{x=\zeta^{-1}(y,t)}$ and $\widetilde{Du_t}|_{t=0}=\tilde{U}(x,0)$.  Hence, from \eqref{2020sa854}
and \eqref{2020sa8sfda54} it follows that $\tilde{U}=\widetilde{Du_t}$, which givest $U=D_t u$.  Thus, in view of
\eqref{202011102145}, \eqref{202012091547} and \eqref{20191224assdfadf2055}, we have
\begin{align}
& \|u_t\|_{\mathcal{U}_T} \lesssim_{\mm{L}} 1+\|u^0\|_4^6+\|u^0\|_2^4\|\nabla w_t|_{t=0}\|_0. \label{20191224asdfasassdfadf2055}
\end{align}

Taking into account the regularity of $u$, we get from \eqref{01dsaf16asdfsdasafasf} that
\begin{align}
\label{20191206fda0saf857}
(D_t u)_t+\nabla_{\ml{A}} q_t -\nu \Delta_{\ml{A}}D_t u+ \kappa D_t u = D_t\partial_2^2\eta+\mathcal{F}.
\end{align}
Thus, $Q=q_t$ by virtue of \eqref{20191206fda0857} and \eqref{20191206fda0saf857}.

Similarly to the derivation of \eqref{202011102145}, we can obtain
\begin{align}
&  \|  u \|_{C^0(\overline{I_T},H^4)} +\|  u \|_{L^2({I_T},H^5)}  \nonumber \\
 & \lesssim_{\mm{L}} \left(1+\|\nabla \eta^0\|_4\right)(1+ \|u^0\|_4^6+\|u^0\|_2^4\|\nabla w_t|_{t=0} \|_0),\label{2019122216sdaf11}
 \end{align}
which, together with \eqref{20191224assdfadf2055} with $q_t$ in place of $Q$ and \eqref{20191224asdfasassdfadf2055}, yields \eqref{20191asfas2242055}.
Employing the same arguments as in the proof of \eqref{202010170831}, one gets \eqref{20191asfas2242055xx}. Finally,
the two estimates in \eqref{20191asfasfasas2242055} are obvious to get by using \eqref{202009301500} and \eqref{2020212271533}.
This completes the proof of Proposition \ref{pro:04asfda01nxdxx}. \hfill $\Box$
\end{pf}

\subsection{Proof of Proposition \ref{pro:0401nxdxx}}\label{202009111946}

Now we are in a position to show Proposition \ref{pro:0401nxdxx}. To start with, let $(\eta^0,u^0)$ satisfy all the assumptions
in Proposition \ref{pro:0401nxdxx} and $\|\nabla\eta^0\|_{2,2}^2\leqslant\delta\leqslant\delta^{\mm{L}}_1$.
We should remark here that the smallness of $\delta$ (independent of $m$ and $\nu$) will be frequently used
in the calculations that follow.

Denote
\begin{align}
B_1:= 2c^{\mm{L}}(  1+B^3),
\label{202011111007}
\end{align}
 where $B$ comes from Proposition \ref{pro:0401nxdxx} and  the constant $c^{\mm{L}}$ is the same as in  \eqref{20201101111005}.
  By Proposition \ref{qwepro:0sadfa401nxdxx} with $B_1$ defined by \eqref{202011111007} and Remark \ref{20201101111005}, one can
 easily construct a function sequence $\{u^k,{q}^{k}\}_{k=1}^\infty$ defined on $\Omega_T$ with $T$ satisfying \eqref{201912061028}.
 Moreover,
\begin{itemize}
  \item for $k\geqslant 1$, $({u}^{k+1},q^{k+1})\in  \mathcal{U}_{1,T}\times C^0(\overline{I_T} ,\underline{H}^2) $ and
 \begin{equation}\label{iteratiequat}\begin{cases}
\eta^k=\int_0^tu^k\mm{d}\tau+\eta^0,\\  u_t^{k+1}+\nabla_{\mathcal{A}^k} {q}^{k+1}-\nu \Delta_{\mathcal{A}^k} {u}^{k+1}=m^2\partial_2^2\eta^k ,\\[1mm]
\mathrm{div}_{\mathcal{A}^k} {u}^{k+1}=0\end{cases}  \end{equation}
with initial condition ${u}^{k+1} |_{t=0}=u^0$,
where $\mathcal{A}^k$ is defined by $\zeta^k:=\eta^k+y$;
  \item $({u}^1,q^1)$ is constructed by Proposition \ref{qwepro:0sadfa401nxdxx} with $w=0$ and $\partial_2^2\eta^0$ in place of $f$;
\item
the solution sequence $\{  {u}^k,q^k\}_{k=1}^\infty$ satisfies the following uniform estimates: for all $k\geqslant 1$,
\begin{align}
\label{2020100sfas11506}
& 1\leqslant 2\det(\nabla \eta^k+I)\leqslant 3,\  \|\nabla \eta^k\|_{2,2}\leqslant 2\delta\;\;\mbox{ for all }t\in\overline{I_T}, \\
 & \|u^k \|_{\mathcal{U}_{1,T}}\leqslant B_1\;\;\mbox{ and }\|\nabla \eta^k\|_{2} +\|q^k \|_{C^0(\overline{{I}_T},H^2)}
 \lesssim_{\mm{L}} 1+B_1^2\label{n053112}.\end{align}
\end{itemize}

In order to take limits in \eqref{iteratiequat} as $k\to\infty$, we have to show that $\{{u}^k,q^k\}_{k=1}^\infty$ is a Cauchy sequence.
To this end, we define for $k\geqslant 2$,
$$( \bar{\eta}^{k}, \bar{ {u}}^{k+1}, \bar{\mathcal{A}}^k, \bar{q}^{k+1}):= (\eta^{k}- {\eta}^{k-1},{u}^{k+1}- {u}^k,\tilde{\mathcal{A}}^k-\tilde{\mathcal{A}}^{k-1},{q}^{k+1} -q^{k}),$$
which satisfies
 \begin{equation}\label{difeequion} \begin{cases}
\bar{\eta}^k=\int_0^t\bar{u}^k\mm{d}\tau,\\ \Delta  \bar{q}^{k+1} =\mathcal{M}_k  , \\[1mm]
 \bar{ { u}}_t^{k+1}+\nabla \bar{q}^{k+1}-\nu \Delta  \bar{ {u}}^{k+1}- m^2 \partial_2^2 \bar{\eta}^k=\mathcal{N}_k,\\[1mm]
\mathrm{div} \bar{u}^{k+1}= -(\mathrm{div}_{\bar{\mathcal{A}}^k } {u}^{k+1}+\mathrm{div}_{\tilde{\mathcal{A}}^{k-1}} \bar{u}^{k+1}), \\
\bar{ {u}}^{k+1} |_{t=0}=0, \end{cases}  \end{equation}
where
$$
\begin{aligned}
\mathcal{M}_k:=&m^2\partial_2^2(\mm{div}\bar{\eta}^k+ \mm{div}_{\bar{\mathcal{A}}^k}  {\eta}^{k}+\mm{div}_{ \tilde{\mathcal{A}}^{k-1}}  \bar{\eta}^{k})\nonumber  \\
& +  \mm{div}_{\bar{\mathcal{A}}^k_t}{u}^{k+1} + \mm{div}_{\mathcal{A}_t^{k-1}}\bar{u}^{k}
-(\mm{div}_{\bar{\mathcal{A}}^k}\nabla_{\ml{A}^k}{q}^{k+1} +
\mm{div}_{\tilde{\ml{A}}^{k-1}}\nabla_{\bar{\mathcal{A}}^k} {q}^{k+1}\nonumber &\\
&+ \mm{div}_{\tilde{\ml{A}}^{k-1}}\nabla_{\tilde{\ml{A}}^{k-1}}\bar{q}^{k+1}
+ \mm{div} (\nabla_{\bar{\mathcal{A}}^k}q^{k+1}+ \nabla_{\tilde{\ml{A}}^{k-1}}\bar{q}^{k+1})), \\
\mathcal{N}_k:=& \nu ( \mm{div}_{\bar{\mathcal{A}}^k}\nabla_{\ml{A}^k}u^{k+1}
+ \mm{div}_{ \tilde{\ml{A}}^{k-1}}\nabla_{\bar{\mathcal{A}}^k}u^{k+1}+\mm{div}_{ \tilde{\ml{A}}^{k-1}}\nabla_{\ml{A}^{k-1}}\bar{u}^{k+1}\\
&+ \mm{div}\nabla_{\bar{\mathcal{A}}^k}u^{k+1}+ \mm{div} \nabla_{\tilde{\ml{A}}^{k-1}}\bar{u}^{k+1})
-(\nabla_{\bar{\mathcal{A}}^k}q^{k+1}+\nabla_{\tilde{\ml{A}}^{k-1}}\bar{q}^{k+1}).
\end{aligned}$$

Keeping in mind that
$$
\begin{aligned}
&\|\bar{\mathcal{A}}^k \|_{2}\lesssim_0 (1+B_1^2)\|\nabla \bar{\eta}^k\|_{2}\lesssim_0 T^{1/2} (1+B_1^2)\| \nabla \bar{u}^k\|_{L^2(I_T,H^2)}, \\
&\|\bar{\mathcal{A}}^k_t  \|_i \lesssim_0 \|\nabla \bar{u}^k\|_i +B_1\|\nabla \bar{\eta}^k\|_{2,2}, \quad i=0,\ 1, \\
&\int|\bar{\mathcal{A}}^k_t||\nabla u^{k+1}||\Delta \bar{q}|\mm{d}y \lesssim_0
 \sqrt{\| \bar{\mathcal{A}}^k_t\|_0 \| \bar{\mathcal{A}}^k_t\|_1} \|u^{k+1}\|_2 \| \Delta \bar{q}\|_{0} ,\\
&\int|{\mathcal{A}}^{k-1}_t||\nabla \bar{u}^{k}||\Delta \bar{q}|\mm{d}y \lesssim_0
 \|\mathcal{A}^k_t\|_1\sqrt{\| \nabla \bar{u}^k\|_0 \| \nabla \bar{u}^k\|_1} \| \Delta \bar{q}\|_{0} ,\\
&\|\nabla \bar{u}^{k}\|_0\leqslant T^{1/2}\| \nabla \bar{u}^k_t\|_{L^2(I_T,L^2)},
\end{aligned}
$$
we make use of \eqref{202004221saffad412},  \eqref{20200508}, \eqref{2020100sfas11506}, \eqref{n053112} and the above five estimates
to deduce from \eqref{difeequion}$_2$--\eqref{difeequion}$_4$ that
\begin{align}
&\|\nabla \bar{q}^{k+1}\|_{C^0(\overline{I_T},H)}
\lesssim_{\mm{L}}  T^{1/4}(1+B_1^2) (\|\nabla\bar{u}^k\|_{L^2(I_T,H^2)}+\| \nabla \bar{u}^k_t\|_{L^2(I_T,L^2)}+ \|\nabla\bar{u}^k\|_{C^0(\overline{I_T},H )}),\label{202001231907}\\[1mm]
&\|\nabla^2  \bar{u}^{k+1}\|_{C^0(\overline{I_T},L^2)}^2+ \|\nabla^3\bar{u}^{k+1} \|_{L^2( {I_T}, L^2)}^2
\lesssim_{\mm{L}} {T} (1+B_1^8)\|\nabla \bar{u}^k \|_{L^2(I_T,H^2)}^2 \nonumber \\
 & \qquad +  (1+B_1^4)( \|\nabla  \bar{q}^{k+1}\|_{C^0(\overline{I_T},H)}^2
 +T^{1/2}(\|\nabla^2  \bar{u}^{k+1}\|_{C^0(\overline{I_T},L^2)}^2+ \|\nabla^3\bar{u}^{k+1} \|_{L^2( {I_T}), L^2)}^2)),
\label{2020012319071}\\[1mm]
& \| \bar{u}_t^{k+1}\|_{C^0(\overline{I_T},L^2)} + \| \bar{u}_t^{k+1}\|_{L^2(I_T,H)}
   \lesssim_{\mm{L}}(1+T^{1/4 }(1+B_1^2))  \|\nabla^2  \bar{u}^{k+1}\|_{C^0(\overline{I_T},L^2)}  \nonumber \\
&\qquad   +(1+T^{1/4 }(1+B_1^2) )  \|\nabla  \bar{u}^{k+1}\|_{C^0(\overline{I_T},H^2)}   \nonumber \\
& \qquad  + \|\nabla  \bar{q}^{k+1}\|_{C^0(\overline{I_T},H)} +   T^{1/2}(1+B_1^5)
 \| \nabla \bar{u}^k \|_{L^2(I_T,H^2)}.
\label{202001232104}
\end{align}
Recalling $( \bar{u}^{k+1})_{\mathbb{T}^2}=0$, we put \eqref{202001231907}--\eqref{202001232104}
together to conclude that for sufficiently small $T$ (depending possibly on $B_1$, $\nu$ and $m$),
$$\|\bar{u}^{k+1}\|_{\mathcal{U}_{1,T}} +\|\bar{q}^{k+1}\|_{C^0(\overline{I_T},H^2)}\leqslant\|\bar{u}^{k}\|_{\mathcal{U}_{1,T}}/{2}
\quad\mbox{ for any }\; k\geqslant 1,$$
which implies
\begin{align}
\sum_{k=2}^\infty(\|\bar{u}^k\|_{\mathcal{U}_T}+\|\bar{q}^{k}\|_{C^0(\overline{I_T},H^2)})<\infty. \nonumber
\end{align}
Hence, $\{ {u}^k,q^k\}_{k=1}^\infty$ is a Cauchy sequence in ${\mathcal{U}_{1,T}}\times C^0(\overline{I_T},\underline{H}^2)$ and
\begin{equation}
\label{strongconvegneuN}
(\eta^k,u^k, q^k)\rightarrow (\eta ,u, q)\quad\mbox{strongly in }C^0(\overline{I_T},H^3)\times
{\mathcal{U}_{1,T}}\times C^0(\overline{I_T}, \underline{H}^2),
\end{equation}
where
\begin{align}
\eta :=\int_0^tu\mm{d}\tau +\eta^0.
\label{20201111413111}
\end{align}

Remembering that \eqref{20201111413111} implies $\eta_t=u$, we infer
from \eqref{iteratiequat} and \eqref{strongconvegneuN} that the limit $(\eta,{u},q)$ is a solution to the initial-value
problem \eqref{01dsaf16asdfasf}--\eqref{01dsaf16asdfasfsaf}.
The uniqueness of solutions to \eqref{01dsaf16asdfasf}--\eqref{01dsaf16asdfasfsaf}
in the function class $C^0(\overline{I_T},H^3)\times {\mathcal{U}_{1,T}}\times C^0(\overline{I_T}, \underline{H}^2)$
is easily verified by a standard energy method, and its proof will be omitted here.
The proof of Proposition \ref{pro:0401nxdxx} is complete.

Similar to Proposition \ref{pro:0401nxdxx}, we can use Proposition \ref{pro:04asfda01nxdxx} to establish the following existence
and uniqueness of a classical solution to the problem \eqref{01dsaf16asdfsdasafasf} with an additional damping term.
\begin{pro}
\label{qwepro:0sadassdfafdasfa401nxdxx}
Let $\delta>0$, $(\eta^0,u^0)\in H^5\times H^4$  and $\|\nabla \eta^0\|_3 \leqslant \delta$. Then, there is
a sufficiently small constant $\delta^{\mm{a}}\in (0,\delta_2^{\mm{L}} ]$, independent of $m$ and $\nu$, such that
for any $\delta\leqslant \delta^{\mm{a}}$, the initial-value problem \eqref{201912060857} admits a unique local strong solution
$(\eta^\nu,u^\nu,q^\nu)\in C^0(\overline{I_T}, H^5)\times\mathcal{U}_{T}^2\times C^0(\overline{I_T},\underline{H}^4)$,
where $T:=\min\{( 2c^\kappa\|( u^0, m\partial_2\eta^0,\sqrt{\nu}\nabla\eta^0)\|_4^2)^{-1},\delta(\sqrt{2 c^\kappa }\|
 (u^0, m\partial_2 \eta^0,\sqrt{\nu} \nabla \eta^0) \|_4)^{-1}\}$, and the constant $c^\kappa$ is the same as in the definition of $T$.
Moreover,
\begin{align}
& 1\leqslant 2\det(\nabla \eta+I)\leqslant 3,\ \|\nabla \eta^\nu\|_3\leqslant  2 \delta , \label{202010211335} \\
& \| \eta^\nu  \|_{4} \lesssim_\kappa    \| \eta^0\|_{4} +t\|(u^0,m \partial_2 \eta^0,\sqrt{\nu} \nabla \eta^0)\|_{4}
\;\;\;\mbox{ for any }t\in \overline{I_T}, \\
& \|(u^\nu,m\partial_2 \eta^\nu )\|_{C^0(\overline{I_T},H^4)} +\sqrt{\nu}\|u^\nu\|_{L^2(I_T,H^5)}
\lesssim_\kappa\|(u^0,m \partial_2 \eta^0,\sqrt{\nu} \nabla \eta^0)\|_{4} , \\
& \|( u^\nu_t , q^\nu ,q^\nu_t) \|_{C^0(\overline{I_T},H^2 \times H^4\times H^2)} \lesssim_\kappa I_0,  \label{20201111sfaf2113}
\end{align}
where
$$
\begin{aligned}
I_0:=(1 & +\|(u^0,m \partial_2 \eta^0,\sqrt{\nu} \nabla \eta^0)\|_{4}  )(\|( u^0,m\partial_2 \eta^0, \sqrt{\nu}  \nabla \eta^0 )\|_{4}^2\\
&+(1+\nu+m^2)\|( u^0,m\partial_2 \eta^0, \sqrt{\nu}  \nabla \eta^0 )\|_{4}) .
\end{aligned}
$$
\end{pro}
\begin{pf}
We divided the proof into three steps.

(1) Let  $\delta\in (0,\delta_2^{\mm{L}} ]$, $(\eta^0,u^0)\in H^5\times H^4$ and $\|\nabla \eta^0\|_3 \leqslant \delta$.  Let $ {B}^\kappa  >0$
be an undetermined constant that satisfies ${B}^\kappa\geqslant \|(\nabla \eta^0,u^0) \|_4$ and will be defined in \eqref{202010211938}.
Thanks to Proposition \ref{pro:04asfda01nxdxx}, we can follow the same arguments as in Section \ref{202009111946} to deduce
 that there is a constant $\delta_1^{\mm{a}}$ independent of $\nu$ and $B^\kappa$, such that for any $\delta\leqslant \delta_1^{\mm{a}}$,
 \begin{itemize}
     \item there are a function sequence $\{ \eta^k, u^k,{q}^{k} \}_{k=1}^\infty\in C^0(\overline{I_{T^\nu}}, H^5)
     \times\mathcal{U}^2_{T^\nu} \times C^0(\overline{I_{T^\nu}}, \underline{H}^4) $ and a limit function $(\eta^\nu,u^\nu,q^\nu )$,
      such that as $k\to\infty$,
\begin{align}\label{strongcsafaonvegneuN}
& (\eta^k,u^k, q^k,q^k_t)\rightarrow (\eta^\nu,u^\nu, q^\nu,q_t^\nu)  \nonumber \\
 & \qquad \mbox{ in }C^0(\overline{I_{T^\nu}}, H^5) \times  \mathcal{U}_{T^\nu}^2 \times  C^0(\overline{I_{T^\nu}}, H^4)
\times C^0( {I_{I_{T^\nu}}}, \underline{H}^2) ,  \\
& \label{01dsaf16asdfsadfsasdasafasf}
                              \begin{cases}
\eta^\nu_t=u^\nu, \\[1mm]
 u_t^\nu+\nabla_{\ml{A}}q^\nu- \nu \Delta_{\ml{A}}u^\nu+\kappa u^\nu=
 m^2\partial_2^2\eta^\nu, \\[1mm]
\div_{\ml{A}^\nu}u^\nu=0,\\
 (\eta^\nu,u^\nu)|_{t=0}=(\eta^0,u^0),
\end{cases} \\
&1\leqslant 2\det(\nabla \eta^\nu +I)\leqslant 3,\quad \sup_{t\in \overline{I_{T^\nu}}} \|\nabla \eta^\nu\|_3
\leqslant \|\nabla \eta^0\|_3+\delta <\delta':=4\delta, \label{2020111141603}
\end{align}
where the local existence time $T^\nu_\delta\in (0,1]$ may depend on ${B}^\kappa$, $\nu$, $m$ and $\delta$.
       \item
the function $(\eta^\nu,u^\nu,q^\nu )$ is just the unique solution of \eqref{01dsaf16asdfsadfsasdasafasf}, i.e.,
if there is another solution $( \tilde{\eta}^\nu, \tilde{u}^\nu,\tilde{q}^\nu)$ in
$C^0(\overline{I_{T^\nu}}, H^5) \times  \mathcal{U}_{T^\nu}^2 \times  C^0(\overline{I_{T^\nu}}, H^4)$ satisfying
$0<\inf_{(y,t)\in \mathbb{R}^2\times \overline{I_T}} \det(\nabla \tilde{\eta}^\nu+I)$, then $(\eta^\nu,u^\nu,q^\nu)=(\tilde{\eta}^\nu,\tilde{u}^\nu,\tilde{q}^\nu)$ by using the smallness condition
``$\sup_{t\in \overline{I_T}}\|\nabla \eta^\nu\|_{3}\leqslant 4\delta $''.
                                             \end{itemize}
From now on, we further take $\delta\leqslant \delta_1^{\mm{a}}/2$, then the definition
$T^\nu_{\min}:=\min\{ T^\nu_{\delta},T^\nu_{2\delta} \}$ makes sense.

(2) Noting that \eqref{202008061631} holds for any $\chi\in L^2$, and
$$\mathrm{div}\eta^\nu= \partial_1\eta_2^\nu\partial_2 \eta_1^\nu
-\partial_1\eta_1^\nu\partial_2 \eta_2^\nu +\det(\nabla \eta+I)|_{t=0}-1,$$
we have the inequality:
\begin{align}
&\|(u^\nu,m \partial_2 \eta^\nu)\|_{4}^2 \lesssim_0
 \|(u^\nu,\nabla^3  \mm{curl}_{\mathcal{A}^\nu}(u^\nu,m \partial_2 \eta^\nu  ))\|_0^2
+ \|\eta^\nu \|_4^2\|   (u^\nu,m \partial_2 \eta^\nu  ) \|_4^2\nonumber \\
& \quad + \| \eta^\nu\|_4^2\| u^\nu\|_4^2+ \|\eta^\nu\|_4^2 \|m\partial_2 \eta^\nu\|_4^2 +  \| m\partial_2\eta^0\|_4^2  .
\label{202012262107}
\end{align}
On the one hand, remembering that $\sup_{t\in \overline{I_{T^\nu_{\min}}}} \|\nabla \eta\|_3\leqslant \delta'$ in \eqref{2020111141603},
we can use \eqref{01dsaf16asdfsadfsasdasafasf}$_1$, \eqref{01dsaf16asdfsadfsasdasafasf}$_2$ and \eqref{202012262107},
 and follow the same process (under slight modifications) as in the derivation of \eqref{211qwessebsadd}, to
 deduce that there is a constant $\delta_2^{\mm{a}}$ independent of $\nu$ and $B^\kappa$, such that
 for any $\delta'\leqslant \delta_2^{\mm{a}}$, the solution $(\eta^\nu,u^\nu )$ satisfies
\begin{align}
&  \frac{\mm{d}}{\mm{d}t}   \|(u^\nu,m\partial_2\eta^\nu,\nabla^3  \mm{curl}_{\mathcal{A}^\nu}(u^\nu,m \partial_2 \eta^\nu  ))\|_0^2  + \nu \|   u^\nu \|_5^2+\kappa \|u^\nu\|_4^2/2  \nonumber  \\
&\qquad \lesssim_\kappa   \| ( u^\nu,m\partial_2\eta^\nu)\|_4^4+\nu   \|u^\nu\|_3^2\|\nabla \eta^\nu\|_4^2 ,\label{2020111130859}
 \end{align}
where
\begin{eqnarray*} &&
\|(u^\nu,m \partial_2 \eta^\nu)\|_{4}^2  \lesssim_0 \|(u^\nu,\nabla^3  \mm{curl}_{\mathcal{A}^\nu}(u^\nu,m \partial_2 \eta^\nu  ))\|_0^2
+\| m\partial_2\eta^0\|_4^2, \\[1mm]
&& \|(u^\nu,\nabla^3  \mm{curl}_{\mathcal{A}^\nu}(u^\nu,m \partial_2 \eta^\nu  ))\|_0^2\lesssim_0 \|(u^\nu,m \partial_2 \eta^\nu)\|_{4}^2.
\end{eqnarray*}

On the other hand, by \eqref{01dsaf16asdfsadfsasdasafasf}$_1$ we find that
\begin{align}
\label{202010171448}
\|\nabla\eta^\nu(t)\|_4\leqslant \|u^\nu\|_{L^2(I_t,H^5)}+\|\nabla\eta^0\|_4.
\end{align}
Thus, one concludes from \eqref{2020111130859}--\eqref{202010171448} that for any $t\in (0,T^\nu_{\min}]$,
\begin{align}
\|(u^\nu,m \partial_2 \eta^\nu)\|_{4}^2 +\nu\int_0^t  \| u^\nu\|_5^2  \mm{d}\tau \leqslant \frac{c^\kappa \|(u^0,m \partial_2 \eta^0,\sqrt{\nu} \nabla \eta^0)\|_{4}^2 }{1- c^\kappa t\|(u^0,m \partial_2 \eta^0,\sqrt{\nu} \nabla \eta^0)\|_{4}^2  }. \label{202sdfaf008061949}
\end{align}
In particular, taking $\tilde{T}=1/2 c^\kappa \|  (u^0, m  \partial_2 \eta^0,\sqrt{\nu} \nabla \eta^0) \|_4^2$ , we derive from \eqref{202sdfaf008061949} that for any $t\leqslant \min\{\tilde{T},T^\nu_{\min}\}$,
\begin{align}
  \| (  u^\nu , m \partial_2 \eta^\nu )\|_{4} + \sqrt{\nu}  \| u^\nu   \|_{L^2(I_t,H^5)} \leqslant  \sqrt{2c^\kappa} \|(u^0,m \partial_2 \eta^0,\sqrt{\nu} \nabla \eta^0)\|_{4} =: {B}^\kappa_1.\label{202sdfafsadf008061949}
\end{align}

Let $T^{\min}=\min\{\tilde{T},T^\nu_{\min} ,\delta/{B}^\kappa_1\}$.With the help of \eqref{202sdfaf008061949}--\eqref{202sdfafsadf008061949},
we then infer from \eqref{01dsaf16asdfsadfsasdasafasf}$_1$ that for any $t\leqslant T^{\min}$,
\begin{align}
&   \|  \eta^\nu\|_4 \leqslant  \| \eta^0\|_4+ t {B}^\kappa_1 , \label{201asfdsa7asdfsaf}\\
&   \|\nabla \eta^\nu\|_{3} \leqslant  \|\nabla \eta^0\|_{3}+ t {B}^\kappa_1\leqslant 2\delta  , \label{2017asdfsaf}\\
& \|\nabla \eta^\nu\|_4 \lesssim_\kappa \|\nabla \eta^0\|_4+ \|(u^0,m \partial_2 \eta^0,\sqrt{\nu} \nabla \eta^0)\|_{4}^2\sqrt{2 c^\kappa/ {\nu}} =: {B}^\kappa_2. \label{2017asdfsaf2assfdasf1}
\end{align}
Furthermore, by \eqref{2017asdfsaf}, \eqref{01dsaf16asdfsadfsasdasafasf}$_2$ and \eqref{01dsaf16asdfsadfsasdasafasf}$_3$,
we find that there is a constant $\delta_3^{\mm{a}}$ independent of any parameters, such that for any $\delta\leqslant \delta_3^{\mm{a}}$,
\begin{align}
 \|(u_t^\nu,q^\nu,q^\nu_t )\|_{C^0(\overline{I_{T^{\min}}}, H^2\times H^4\times H^2)}\lesssim_\kappa I_0. \label{202sdfafsadf0080asf61949}
\end{align}

(3) Now, if we take
\begin{align}
B^\kappa:=\max\{{B}^\kappa_1/2,B^\kappa_2/2,\|(\nabla \eta^0,u^0) \|_4\},
\label{202010211938}
\end{align}
then we have by \eqref{202sdfafsadf008061949} and \eqref{2017asdfsaf2assfdasf1} that
\begin{align}
  \| ( \nabla \eta,u^\nu)\|_4\leqslant B^\kappa\;\;\;\mbox{ for any }t\leqslant T^{\min}.
\label{202011141926}
\end{align}
Denote $\tilde{T}^{\min} =\{\tilde{T},\delta/{B}^\kappa_1\}$. If $\tilde{T}^{\min}\leqslant T^\nu$, we easily see that the conclusion
in Proposition \ref{qwepro:0sadassdfafdasfa401nxdxx} holds for any given $\delta\in (0,\delta^{\min}]$, where $\delta^{\min}:=\min\{\delta_1^{\mm{a}}, \delta_2^{\mm{a}}/4, \delta_3^{\mm{a}}\}$.
Next, we consider the case $\tilde{T}^{\min}> T^\nu_{\min}$.

For the case $\tilde{T}^{\min}>T^\nu_{\min}$, we can take $(\eta(T^\nu_{\min}),u(T^\nu_{\min}))$ as initial data. By \eqref{2017asdfsaf}, \eqref{202011141926} and Step (1), we see that for any $\delta\leqslant  \delta^{\min}/2$,
there exists a unique local solution $(\eta^*,u^*,q^*)$ defined on $\mathbb{T}^2\times (T^\nu_{\min}, 2T^\nu_{\min} ]$ of
 the initial-value problem:
$$\begin{cases}
\eta^*_t=u^*_t, \\[1mm]
 u_t^*+\nabla_{\ml{A}}q^*- \nu \Delta_{\ml{A}}u^*+\kappa u^*=
 m^2\partial_2^2\eta^*, \\[1mm]
\div_{\ml{A}^*}u^*=0,\\
 (\eta^*,u^*)|_{t=T^\nu}=(\eta^\nu(T^\nu),u^\nu(T^\nu)).
\end{cases}$$
Moreover, $1\leqslant 2\det(\nabla \eta^*+I)\leqslant 3$ and $  \|\nabla \eta^*(t)\|_{3} \leqslant \delta'$
for any $t\in ( T^\nu_{\min}, 2T^\nu_{\min}]$.  {Due to the uniqueness,
we can get a new local solution} defined on $\mathbb{T}^2\times (0, 2T^\nu_{\min} ]$, still denoted this new solution
by $(\eta^\nu,u^\nu,q^\nu)$. By Step (2) we find that for any $\delta\leqslant\delta^{\min}/2$, the new solution satisfies \eqref{202sdfafsadf008061949}--\eqref{2017asdfsaf}, \eqref{202sdfafsadf0080asf61949} and \eqref{202011141926} for any
$t\in (0,\min\{\tilde{T}^{\min},2T^\nu_{\min}\}]$.
Let $n=[(\tilde{T}^{\min}-T^\nu_{\min})/T^\nu_{\min}]+1$, where $[\cdot]$ means the integer part.
Therefore, by performing $n$-times extension with respect to time, we obtain the desired conclusion
in Proposition \ref{qwepro:0sadassdfafdasfa401nxdxx}.
\hfill $\Box$
\end{pf}

\subsection{Proof of Proposition \ref{pro:0401nsadfaxsdfafdsaddfdxx}} \label{202008062143}
  With the help of Proposition \ref{qwepro:0sadassdfafdasfa401nxdxx}, we are able to prove Proposition \ref{pro:0401nsadfaxsdfafdsaddfdxx}
  by the method of vanishing viscosity limit.

Let $B^\kappa>0$, $\delta\in (0,\delta^a]$, $(\eta^0,u^0)$ satisfy the assumptions in Proposition \ref{pro:0401nsadfaxsdfafdsaddfdxx}
and $\|\nabla \eta^0\|_3\leqslant \delta $, where $\delta^a$ is the same constant as in Proposition \ref{qwepro:0sadassdfafdasfa401nxdxx}.

Let $\varepsilon\in (0,1)$ and $S_\varepsilon$ be a standard mollifier (or a regularizing operator, see \cite[Section 1.3.4.4]{NASII04}
for definition).
It is well-known that $S_\varepsilon(\vartheta)\in C^\infty(\mathbb{T}^2)$, and
$\|S_\varepsilon(\vartheta)\|_{i} \leqslant \tilde{c}_i\|\vartheta\|_{i}$  for $\vartheta\in H^i$ and $i \geqslant 0$,
where the positive constant $\tilde{c}_i$ depends on $i$ only.  Let  $(\eta^0_\varepsilon,u^0_\varepsilon,\tilde{\eta}^0_\varepsilon,\tilde{u}^0_\varepsilon)=S_\varepsilon( \eta^0,u^0,\tilde{\eta}^0,\tilde{u}^0)$.

We now fix $\varepsilon>0$. By virtue of Proposition \ref{qwepro:0sadassdfafdasfa401nxdxx}, there is a sufficiently small $\nu^0$
 (depending possibly on $\|S_\varepsilon(\eta^0)\|_5$, $\|( u^0,m \partial_2 \eta^0 )\|_{4} $ and $\|\nabla\eta^0\|_5$),
 such that for any $\nu<\nu^0$, there exists a unique solution
 $(\eta^\nu,u^\nu,q^\nu)\in C^0(\overline{I_{T}}, H^5) \times \mathcal{U}_{T}^2 \times  C^0(\overline{I_{T}}, \underline{H}^4)$
 to the initial-value problem:
\begin{equation}
\label{20191206safas0857}
\begin{cases}
\eta^\nu=\int_0^tu ^\nu\mm{d}\tau+\eta^0_\varepsilon ,\\
u_t^\nu+\nabla_{\ml{A}^\nu}q^\nu+ \kappa u^\nu-\nu \Delta_{\ml{A}^\nu}u^\nu=
m^2\partial_2^2\eta^\nu ,    \\
\mm{div}_{\mathcal{A}^\nu}u^\nu=0 ,\\
u^\nu|_{t=0}=u^0_\varepsilon\mbox{ in }\mathbb{T}^2,
\end{cases}
\end{equation}
where $\mathcal{A}^\nu=(\nabla \eta^\nu+I)^{-\mm{T}}$,
$T:=\min\{( 1+2c^\kappa \|  (u^0, m  \partial_2 \eta^0) \|_4^2)^{-1}, \delta(1+ \sqrt{2 c^\kappa }\|(u^0, m \partial_2 \eta^0 )\|_4)^{-1}\}$.
Moreover, the solution satisfies the uniform estimates:
\begin{align}
& 1\leqslant 2\det(\nabla \eta^{\nu}+I)\leqslant 3,\ \|\nabla \eta^\nu\|_3\lesssim_0 \delta  ,\label{20201safd0211335} \\
 &\| \eta^\nu\|_{4} \lesssim_\kappa\|( \eta^0,u^0,m\partial_2 \eta^0,\sqrt{\nu} \nabla\eta^0)\|_{4}\;\;\mbox{ for each }t\in \overline{I_T}, \\
& \|( u^\nu,m\partial_2\eta^\nu ) \|_{C^0(\overline{I_{T}},H^4  )}+\sqrt{\nu}  \| u^\nu   \|_{L^2(I_T,H^5)} \lesssim_\kappa\|( u^0,m \partial_2 \eta^0 )\|_{4} , \label{202012262132} \\   & \|( u^\nu_t ,q^\nu ,q^\nu_t) \|_{C^0(\overline{I_{T}},  H^2
 \times H^4\times H^2)} \lesssim_\kappa ((1+ m^2)(1+\|( u^0,m \partial_2 \eta^0 )\|_{4}))^3,  \label{202011112113} \\
& \nu\|\eta^\nu   \|_5^2\lesssim_0\nu \|\eta^0_\varepsilon\|_5^2+t\nu\int_0^t\|u^\nu   \|_5^2\mm{d}\tau
\lesssim_0 \nu\|\eta^0_\varepsilon\|_5^2+t \|( u^0,m \partial_2 \eta^0 )\|_{4}^2.
\label{2020212272138} \end{align}
From now on, we take $\nu=\nu_n=1/n$ with $n\geqslant 1/\nu^0$, and renew to define $(\eta^{\nu_n},u^{\nu_n},q^{\nu_n})$
and $\mathcal{A}^{\nu_n}$ by $(\eta^n,u^n,q^n)$ and $\mathcal{A}^n$, respectively.

Let $({\eta}^{\mm{L}},
{u}^{\mm{L}}) \in\mathfrak{C}^0(\mathbb{R}^+,{H}^{5}_2 )\times {\mathfrak{U}}_\infty^4$ be the unique global solution of the linear
initial-value problem \eqref{202012241017}. In view of Proposition \ref{2020212240930}, we see that there is a unique strong solution
$( {\eta}^{\varepsilon,\mm{L}},{u}^{\varepsilon,\mm{L}})$ to \eqref{202012241017} with $(\tilde{\eta}^0_\varepsilon, \tilde{u}^0_\varepsilon)$
in place of $(\tilde{\eta}^0,\tilde{u}^0)$. Moreover,
\begin{itemize}
  \item for $i\geqslant 0$, $( {\eta}^{\varepsilon,\mm{L}},{u}^{\varepsilon,\mm{L}}) \in\mathfrak{C}^0(\mathbb{R}^+,{H}^{i+1}_2 )
  \times {\mathfrak{U}}_\infty^i$,  and for any $t\geqslant 0$,
\begin{align}
& \|(  u^{\varepsilon,\mm{L}},m\partial_2\eta^{\varepsilon,\mm{L}} )\|_i^2+ \int_0^t \|u^{\varepsilon,\mm{L}}\|_i^2\mm{d}\tau
\lesssim_\kappa \|( \tilde{u}^0_\varepsilon,m\partial_2\tilde{\eta}^0_\varepsilon)\|_i^2 ; \label{202012272141}
\end{align}
  \item  for any $t\in I_T$,
\begin{align}
({\eta}^{\varepsilon,\mm{L}},{u}^{\varepsilon,\mm{L}},\partial_2\eta^{\varepsilon,\mm{L}} )
\to (\eta^{\mm{L}},u^{\mm{L}},\partial_2\eta^{\mm{L}} )\;\;\mbox{ weakly-$*$ in }  {L^\infty ({I_t},H^4 )}\;\;\mbox{ as }\varepsilon\to 0.
\label{2020212271422xx}
\end{align}
\end{itemize}

Defining $({\eta}^{n,\varepsilon,\mm{d}},{u}^{n,\varepsilon,\mm{d}}):=(\eta^n,u^n)-(\eta^{\varepsilon,\mm{L}},u^{\varepsilon,\mm{L}})$,
then
\begin{equation}\label{01dsfsafdsaasfxx}
\begin{cases}
\eta_t^{n,\varepsilon,\mm{d}}=u^{n,\varepsilon,\mm{d}}, \\[1mm]
u_t^{n,\varepsilon,\mm{d}}+\nabla_{\mathcal{A}^n} q^n+ \kappa u^{n,\varepsilon,\mm{d}} =    m^2
\partial_2^2\eta^{n,\varepsilon,\mm{d}}+\nu\Delta_{\mathcal{A}^n} u^\nu, \\[1mm]
\div u^{n,\varepsilon,\mm{d}}=-\mathrm{div}_{\tilde{\mathcal{A}^n}} {u},\\
(\eta^{n,\varepsilon,\mm{d}},u^{n,\varepsilon,\mm{d}})|_{t=0}=(\eta^0_\varepsilon - \tilde{\eta}^0_\varepsilon,  u^0-\tilde{u}^0_\varepsilon).
\end{cases}
\end{equation}
It is easy to see from \eqref{20191206safas0857}$_3$, \eqref{20201safd0211335} and \eqref{01dsfsafdsaasfxx} that for any $t\in I_T$,
\begin{align}\nonumber
 &\| \mm{curl}_{\mathcal{A}^n}(  u^{n,\varepsilon,\mm{d}}, m \partial_2 \eta^{n,\varepsilon,\mm{d}} )\|_{3}^2\\
&\lesssim_0  e^{t  \|( u^n,m\partial_2\eta^n) \|_{C^0(\overline{I_{t}},  H^4\times H^4 )}   }\Bigg(\| ( u^{n,\varepsilon,\mm{d}},m \partial_2 \eta^{n,\varepsilon,\mm{d}})|_{t=0}\|_4^2+ t  \sup_{\tau\in \overline{I_t}}( \|( u^n,m\partial_2\eta^n, m\partial_2\eta^{\varepsilon,\mm{L}}\|_4^3 )\nonumber \\
&\qquad+\nu t  \sup_{\tau\in \overline{I_t}}(\| \eta^n\|_5^2   \|(u^{\varepsilon,\mm{L}},u^{n})\|_4^2  +\| \eta^n\|_5 (  \|u^{\varepsilon,\mm{L}}\|_5^2+ \|u^n\|_4^2))  +  \nu \int_0^t \|u^{\varepsilon,\mm{L}}\|_5^2 \mm{d}\tau\Bigg).  \nonumber
\end{align}
and
\begin{align}\nonumber
 \| \mm{curl}_{\mathcal{A}^n} \eta^{n,\varepsilon,\mm{d}} \|_{3} \lesssim_0\| \mm{curl}_{\mathcal{A}^n} \eta^{n,\varepsilon,\mm{d}} |_{t=0}\|_3+t\sup_{\tau\in \overline{I_t}}\| \mm{curl}_{\mathcal{A}^n} u^{n,\varepsilon,\mm{d}} \|_{3}.
\end{align}
If we making use of \eqref{202012262132}, \eqref{2020212272138}, \eqref{202012272141} and the above two estimates, we further obtain
\begin{align}\nonumber
&\sup_{\tau\in \overline{I_t}}\| \mm{curl}_{\mathcal{A}^n}( \eta^{n,\varepsilon,\mm{d}},  u^{n,\varepsilon,\mm{d}}, m \partial_2 \eta^{n,\varepsilon,\mm{d}} )(\tau)\|_{3}^2 \\ &
\lesssim_\kappa  e^{  t  \|(\eta^0,u^0,m \partial_2 \eta^0 )\|_{4} } (  \|(\eta^0-\tilde{\eta}^0,u^0-\tilde{u}^0,m \partial_2 (\eta^0-\tilde{\eta}^0) )\|_{4}^2\nonumber \\
&\qquad+t(\|( u^0,m \partial_2 \eta^0, \tilde{u}^0,m \partial_2 \tilde{\eta}^0 )\|_{4}^3+
 (t \|( u^0,m \partial_2 \eta^0 )\|_{4}^2 \nonumber \\
&\qquad +\nu   \|\eta^0_\varepsilon\|_5^2)\|( u^0,m \partial_2 \eta^0, \tilde{u}^0,m \partial_2 \tilde{\eta}^0 )\|_{4}^2+\sqrt{\nu} (\sqrt{t} \|( u^0,m \partial_2 \eta^0 )\|_{4} \nonumber \\
 &\qquad + \sqrt{\nu}   \|\eta^0_\varepsilon\|_5)(\|( \tilde{u}^0_\varepsilon,m\partial_2\tilde{\eta}^0_\varepsilon)\|_5^2 +\|( u^0,m \partial_2 \eta^0 )\|_{4}^2)+\nu  \|(  \tilde{u}^0_\varepsilon,m\partial_2\tilde{\eta}^0_\varepsilon)\|_5^2 ).
\label{20212272147}
\end{align}

Thanks to the uniform estimates \eqref{20201safd0211335}--\eqref{202011112113} and \eqref{20212272147}, we can choose a sequence of $\{\eta^n,u^n,q^n\}$ (still labelled by $\{\eta^n,u^n,q^n\}$ for the sake of simplicity), such that as $n\to\infty $,
\begin{align}
&(\eta^n, \partial_2\eta^n,  u^n,q^n )
\to (\eta^\varepsilon,\partial_2\eta, u^\varepsilon,q^\varepsilon) \;\mbox{ strongly in }{C^0(\overline{I_{T}},H^3\times H^3\times  H^3\times \underline{H}^3 )},\nonumber \\
&(\eta^n, u^n,u_t^n, q^n, q^n_t ) \to (\eta^\varepsilon ,u^\varepsilon,u_t^\varepsilon, q^\varepsilon , q^\varepsilon_t)
\; \mbox{ weakly-$*$ in }  L^\infty(I_{T},H^5_2  \times H^4  \times H^2\times H^4\times H^2)  ,\nonumber\\
& ( \mm{curl}_{\mathcal{A}^n}\eta^{n,\varepsilon,\mm{d}},  \mm{curl}_{\mathcal{A}^n}u^{n,\varepsilon,\mm{d}},
\mm{curl}_{\mathcal{A}^n}\partial_2 \eta^{n,\varepsilon,\mm{d}})\to (\mm{curl}_{\mathcal{A}^\varepsilon} \eta^{\varepsilon,\mm{d}},  \mm{curl}_{\mathcal{A}^\varepsilon} u^{\varepsilon,\mm{d}},\mm{curl}_{\mathcal{A}^\varepsilon}\partial_2 \eta^{\varepsilon,\mm{d}} )\nonumber \\
&\quad \mbox{ weakly-$*$ in } L^\infty(I_t,H^3)\; \mbox{ for any }t\in (0,T] ,\nonumber
\end{align}
where
$(\eta^{\varepsilon,\mm{d}},u^{\varepsilon,\mm{d}}):=(\eta^{\varepsilon}- {\eta}^{\varepsilon,\mm{L}},u^{\varepsilon} -u^{\varepsilon,\mm{L}})$.
Moreover, the limit functions $\eta^\varepsilon$, $u^\varepsilon$ and $q^\varepsilon$ solve the initial-value problem
\begin{equation}
\label{20191206safasadfs0857}
  \begin{cases}
 \eta_t^\varepsilon=u^\varepsilon   ,\\
u_t^\varepsilon +\nabla_{\ml{A}}q^\varepsilon + \kappa u^\varepsilon   =
m^2\partial_2^2\eta^\varepsilon,    \\
\mm{div}_{\mathcal{A}}u^\varepsilon=0,  \\
u^\varepsilon|_{t=0}=u^0 ,
  \end{cases}
  \end{equation}
and satisfies
\begin{align}
& 1\leqslant 2\det(\nabla \eta^\varepsilon+I)\leqslant 3,\  \|\nabla \eta^\varepsilon\|_{L^\infty({I_{T}},H^3)}\leqslant  2 \delta   ,\label{2020011121337} \\
&\mm{ess}\sup_{t\in I_T} \|(\eta^\varepsilon,u^\varepsilon,m\partial_2\eta^\varepsilon ) \|_{4} \lesssim_\kappa  \|(\eta^0,u^0,m \partial_2 \eta^0 )\|_{4} ,\label{2020011121337xx}  \\   & \|( u^\varepsilon_t ,q^\varepsilon ,q^\varepsilon_t) \|_{L^\infty( {I_{T}}, H^2\times H^4
 \times H^2 )} \lesssim_\kappa ((1+ m^2)(1+\|(\eta^0,u^0,m \partial_2 \eta^0 )\|_{4}))^3,   \label{20221safsa13}\\
  &   \mm{ess}\sup_{\tau\in I_t}\| (\mm{curl}_{\mathcal{A}} \eta^{\varepsilon, d}, \mm{curl}_{\mathcal{A}} u^{\varepsilon, d},  \mm{curl}_{\mathcal{A}}m \partial_2 \eta^{\varepsilon, d} )(\tau)\|_{3}^2 \nonumber  \\ &
\lesssim_0   e^{ 2 t  \|(\eta^0,u^0,m \partial_2 \eta^0 )\|_{4} } \|(\eta^0-\tilde{\eta}^0,u^0-\tilde{u}^0,m\partial_2(\eta^0-\tilde{\eta}^0)\|_{4}^2 \nonumber \\
& \quad + t(\|( u^0,m \partial_2 \eta^0, \tilde{u}^0,m \partial_2 \tilde{\eta}^0 )\|_{4}^3+
  \|( u^0,m \partial_2 \eta^0 )\|_{4}^2  \|( u^0,m \partial_2 \eta^0, \tilde{u}^0,m \partial_2 \tilde{\eta}^0 )\|_{4}^2). \label{202012271254}
\end{align} In addition, by $u^\varepsilon\in L^\infty(I_{T},H^4)$, one has
\begin{align}
\label{2020120292103}
 \eta^\varepsilon \in C^0(\overline{I_{T}},H^4).
\end{align}

Let $\alpha$ and $\beta$ satisfy $|\alpha|=4$ and $|\beta|=3$. With the help of the uniform estimates
\eqref{20201safd0211335}--\eqref{202011112113},we can use \eqref{20191206safas0857}$_1$ and \eqref{20191206safas0857}$_2$ to deduce
that the sequences $\{\partial^\alpha u^n\}$, $\{\partial^\beta(\mm{curl}_{\mathcal{A}^n}u^n)\}$, $\{\partial^\alpha \partial_2\eta^n \}$ and
$\{\partial^\beta\mm{curl}_{\mathcal{A}^n}\partial_2\eta^n\}$ are uniformly continuous in $H^{-1}$ on $\overline{I}$. Besides,
they are also uniformly bounded in $L^2$. Therefore, there is a sequence of $\{\eta^n,u^n,q^n\}$ (still denoted by $\{\eta^n,u^n,q^n\}$),
 such that (see \cite[Lemma 6.2]{NASII04} for example),
\begin{align}
&\partial^\alpha (  u^n ,  \partial_2\eta^n )\to \partial^\alpha( u^\varepsilon ,\partial_2\eta^\varepsilon)
\mbox{ in }C^0(\overline{I_{T}}, L^2_{\mm{weak}}),\nonumber
\end{align}
where we have relabelled $\partial^\alpha( u^\varepsilon ,\partial_2\eta^\varepsilon)$ on a set of zero-measure in $\overline{I_T}$.
Thus, by \eqref{20212232155}, \eqref{2020120292103} and the above result of weak continuity, the notation ``$\mm{ess}$'' can be removed
in \eqref{2020011121337xx} and \eqref{202012271254}.

Thanks to the uniform estimates \eqref{2020011121337}--\eqref{20221safsa13}
and the limit behavior \eqref{2020212271422xx},
 we can again take to the limit as $\varepsilon\to 0$ by employ the same arguments as used in obtaining
 $(\eta^\varepsilon,u^\varepsilon,q^\varepsilon)$, and thus obtain a limit function $(\eta,u,q)$
 which is a solution of the initial-value problem \eqref{01dsaf16asdfsfdsaasf}--\eqref{01dsafsfda16asdfasfsaf} and
 satisfies the estimates \eqref{20221safsa13} with $(\eta,u,q)$ in place of $(\eta^\varepsilon,u^\varepsilon,q^\varepsilon)$, \eqref{202012301432}--\eqref{20201221asfdasf2151}, and the same regularity as thst of $(\eta^\varepsilon,u^\varepsilon,q^\varepsilon)$.
 Moreover, the obtained solution $(\eta,u,q)$ is unique, provided that $\delta$ is sufficiently small.  To complete the proof of Proposition \ref{qwepro:0sadassdfafdasfa401nxdxx}, obviously, it suffices to show the strong continuity
 of $\partial^\alpha (u,q)$ in $L^2$ with respect to time for any $\alpha$ satisfying $|\alpha|=4$. Next we will verify this fact.

To begin with, we easily see by a straightforward calculation that
\begin{align}
&\partial_t S_\varepsilon(\partial^\beta  \mm{curl}_{\mathcal{A}}u  )+
   \kappa S_\varepsilon( \partial^\beta  \mm{curl}_{\mathcal{A}} u ) -m^2\partial_2
 S_\varepsilon( \partial^\beta \mm{curl}_{\mathcal{A}}\partial_2 \eta )\nonumber
\\
&  =   S_\varepsilon(\partial^\beta \mm{curl}_{\mathcal{A}_t}u   ) -m^2
 S_\varepsilon( \partial^\beta \mm{curl}_{\partial_2 \mathcal{A}}\partial_2 \eta ), \label{2021022041446}
\end{align}
and
\begin{align}
& \partial_tS_\varepsilon(\partial^\beta  \mm{curl}_{ \mathcal{A}}\partial_ 2\eta) =  \partial_ 2 S_\varepsilon(\partial^\beta  \mm{curl}_{\mathcal{A}}u  )+ S_\varepsilon(\partial^\beta \mm{curl}_{\mathcal{A}_t}\partial_2 \eta )
- S_\varepsilon(\partial^\beta \mm{curl}_{\partial_ 2 \mathcal{A}}u). \nonumber
\end{align}
Thanks to the above two identities, we further get that for a.e. $t\in I_T$,
\begin{align}
&\int  \partial_t S_\varepsilon(\partial^\beta  \mm{curl}_{\mathcal{A}}u)   S_\varepsilon(\partial^\beta  \mm{curl}_{\mathcal{A}}u ) \mm{d}y
=\frac{1}{2}\frac{\mm{d}}{\mm{d}t}\int |S_\varepsilon(\partial^\beta  \mm{curl}_{\mathcal{A}}u )|^2\mm{d}y,\nonumber \\
&  - \int \partial_2  S_\varepsilon( \partial^\beta \mm{curl}_{\mathcal{A}}\partial_2 \eta )
 S_\varepsilon(\partial^\beta  \mm{curl}_{\mathcal{A}}u  )\mm{d}y\nonumber \\
 & =\frac{1}{2}\frac{\mm{d}}{\mm{d}t}\int| S_\varepsilon(\partial^\beta  \mm{curl}_{\mathcal{A}}\partial_2\eta)|^2\mm{d}y
 +\int S_\varepsilon( \partial^\beta \mm{curl}_{\mathcal{A}}\partial_2 \eta ) (
   S_\varepsilon(\partial^\beta  \mm{curl}_{\partial_ 2 \mathcal{A}}u)-
S_\varepsilon(\partial^\beta  \mm{curl}_{\mathcal{A}_t}\partial_2 \eta ))\mm{d}y, \nonumber
\end{align}
where $\| S_\varepsilon(\partial^\beta  \mm{curl}_{\mathcal{A}}u )\|_0^2 $  and $\| S_\varepsilon( \partial^\beta \mm{curl}_{\mathcal{A}}\partial_2\eta )\|_0^2\in {AC}^0(I_{T}\setminus \tilde{Z})$ for some zero-measurable set $\tilde{Z}$.

For $\phi\in C_0^\infty(I_{T})$, we multiply \eqref{2021022041446} by $S_\varepsilon(\partial^\beta\mm{curl}_{\mathcal{A}}u)\phi$ in $L^2(\Omega_T)$,
and take then to the limits as $\varepsilon\to 0$ to obtain
\begin{align}
& -  \frac{1}{2}\int_0^T \int (|\partial^\beta  \mm{curl}_{\mathcal{A}}u  |^2
  + m^2 \partial^\beta \mm{curl}_{\mathcal{A}}\partial_2\eta |^2)\mm{d}y\phi_\tau\mm{d}\tau +
  \kappa  \int_0^T\int |\partial^\beta  \mm{curl}_{\mathcal{A}} u  |^2\mm{d}y\phi\mm{d}\tau \nonumber \\
&=\int_0^T    \int (\partial^\beta  \mm{curl}_{\mathcal{A}_t}u \partial^\beta  \mm{curl}_{\mathcal{A}}u \mm{d}y +m^2(
 \partial^\beta \mm{curl}_{\mathcal{A}}\partial_2\eta(
\partial^\beta \mm{curl}_{\mathcal{A}_t}\partial_2\eta-\partial_2\eta
\partial^\beta \mm{curl}_{\partial_2\mathcal{A}}u ) \nonumber
 \\
&\quad\quad\quad\quad-\partial^\beta \mm{curl}_{\partial_2\mathcal{A}}\partial_2\eta
\partial^\beta \mm{curl}_{\mathcal{A}}u ))
\mm{d}y \phi\mm{d}\tau, \nonumber
\end{align}
which implies that for any $\beta$  satisfying $|\beta|=3$,
 \begin{align}
&  \frac{1}{2} \frac{\mm{d}}{\mm{d}t}\|\partial^\beta (\mm{curl}_{\mathcal{A}}u , m \partial^\beta \mm{curl}_{\mathcal{A}}\partial_2\eta )\|_0^2 +
 \kappa   \|  \partial^\beta  \mm{curl}_{\mathcal{A}} u  \|_0^2 \nonumber \\
&=    \int (\partial^\beta  \mm{curl}_{\mathcal{A}_t}u \partial^\beta  \mm{curl}_{\mathcal{A}}u \mm{d}y +m^2
 (\partial^\beta \mm{curl}_{\mathcal{A}}\partial_2\eta
\partial^\beta \mm{curl}_{\mathcal{A}_t}\partial_2\eta \nonumber
 \\
&\quad-\partial^\beta \mm{curl}_{\partial_2\mathcal{A}}\partial_2\eta
\partial^\beta \mm{curl}_{\mathcal{A}}u -\partial^\beta \mm{curl}_{\mathcal{A}}\partial_2\eta
\partial^\beta \mm{curl}_{\partial_2\mathcal{A}}u )) \mm{d}y \quad\mbox{ for a.e. }t\in I_{T},\nonumber
\end{align}
whence,
\begin{align}
\label{202010141017}
\|\partial^\beta   (\mm{curl}_{\mathcal{A}}u , m \mm{curl}_{\mathcal{A}}\partial_2\eta )\|_0^2\in {AC}^0(I_{T}\setminus Z)
\;\;\mbox{ for some zero-measurable set }Z.\end{align}
By \eqref{202010141017} and the fact ``$\partial^\beta( \mm{curl}_{\mathcal{A}}  u  ,\mm{curl}_{\mathcal{A}}\partial_2\eta  )\in C^0(\overline{I_{T}}, L^2_{\mm{weak}})$'',  we immediately  get
\begin{align}
\label{20201111212143}
\partial^\alpha(\mm{curl}_{\mathcal{A}} u ,\mm{curl}_{\mathcal{A}} \partial_2\eta ) \in C(I_{T}\setminus Z, L^2).
\end{align}

Since
$\partial^\alpha (u,\partial_2\eta)\in C^0(\overline{I_{T}}, L^2_{\mm{weak}})$, one has by \eqref{20212232155} that
\begin{align}
\sup_{t\in \overline{I_{T}}}\|\partial^\alpha  (  u,\partial_2\eta)  \|_0 = \mm{ess}\sup_{t\in I_{T}} \|\partial^\alpha  ( u ,\partial^\alpha  \partial_2\eta)  \|_0.
\nonumber
\end{align}
Keeping in mind that for any $t$ and $t_0\in \overline{I_{T}}$,
\begin{align}
 \| \partial^\beta \nabla(u(t)-u(t_0)) \|_0= &\sqrt{\| \partial^\beta \mm{curl}  (u(t)-u(t_0))\|_0^2+\|  \partial^\beta \mm{div}  (u(t)-u(t_0))\|_0^2} \nonumber   \\
\lesssim_0 &\| \partial^\beta (\mm{curl}_{\mathcal{A}}  u(t)-\mm{curl}_{\mathcal{A}^0} u(t_0))\|_0+\| \partial^\beta (\mm{curl}_{\mathcal{A}}u(t_0) -  \mm{curl}_{\mathcal{A}^0}u(t_0))\|_0\nonumber  \\
  &+  \| \partial^\beta (\mm{div}_{\tilde{\mathcal{A}}}  u(t)-\mm{div}_{\tilde{\mathcal{A}}} u(t_0))\|_0+\| \partial^\beta (\mm{div}_{\tilde{\mathcal{A}}}u(t_0)-\mm{div}_{\tilde{\mathcal{A}}^0}u(t_0))\|_0\nonumber
\end{align}
and
\begin{align}
& \| \partial^\beta \nabla \partial_2(\eta(t)-\eta(t_0)) \|_0^2 \nonumber \\
&= \| \partial^\beta \mm{curl} \partial_2 (\eta(t)-\eta(t_0))\|_0^2+\|  \partial^\beta \mm{div}\partial_2  (\eta(t)-\eta(t_0))\|_0^2  \nonumber \\
&\lesssim_0  \| \partial^\beta (\mm{curl}_{\mathcal{A}}  \partial_2\eta(t)-\mm{curl}_{\mathcal{A}^0} \partial_2\eta(t_0))\|_0+\| \partial^\beta (\mm{curl}_{\mathcal{A}}\partial_2\eta(t_0) -  \mm{curl}_{\mathcal{A}^0}\partial_2\eta(t_0))\|_0\nonumber  \\
&\quad +  \| \partial^\beta\partial_2 ((\partial_1\eta_2\partial_2 \eta_1
-\partial_1\eta_1\partial_2 \eta_2)(t)-(\partial_1\eta_2\partial_2 \eta_1
-\partial_1\eta_1\partial_2 \eta_2)(t_0))\|_0 , \nonumber
\end{align}
we employ the regularity $(\eta ,u )\in C^0(\overline{I_{T}},H^4\times H^3)$ to arrive at
$\partial^\alpha( u,\partial_2\eta)\in C(I_{T}\setminus Z, L^2)$. Consequently, $\partial^\alpha q\in C(I_{T}\setminus Z, L^2)$.
The proof of Proposition \ref{qwepro:0sadassdfafdasfa401nxdxx} is complete.

\vspace{4mm} \noindent\textbf{Acknowledgements.}
The research of Fei Jiang was supported by NSFC (Grant Nos. 11671086 and 12022102) and the Natural Science Foundation of Fujian Province of China (2020J02013), and the research of Song Jiang by National Key R\&D Program (2020YFA0712200), National Key Project (GJXM92579), and
NSFC (Grant No. 11631008), the Sino-German Science Center (Grant No. GZ 1465) and the ISF-NSFC joint research program (Grant No. 11761141008).
\renewcommand\refname{References}
\renewenvironment{thebibliography}[1]{ %
\section*{\refname}
\list{{\arabic{enumi}}}{\def\makelabel##1{\hss{##1}}\topsep=0mm
\parsep=0mm
\partopsep=0mm\itemsep=0mm
\labelsep=1ex\itemindent=0mm
\settowidth\labelwidth{\small[#1]}%
\leftmargin\labelwidth \advance\leftmargin\labelsep
\advance\leftmargin -\itemindent
\usecounter{enumi}}\small
\def\newblock{\ }
\sloppy\clubpenalty4000\widowpenalty4000
\sfcode`\.=1000\relax}{\endlist}
\bibliographystyle{model1b-num-names}

\end{document}